\documentclass[reqno, 12pt]{smfart}

\usepackage{srcltx}
\usepackage[utf8]{inputenc}

\usepackage{a4wide}
\usepackage{amsfonts, amsthm, amsmath, amssymb}
\usepackage{amsmath, amsfonts, amssymb, amsthm}
\usepackage{graphicx, subfigure}
\usepackage{fancyhdr}
\usepackage{stmaryrd}
\usepackage{enumerate}
\usepackage{pstricks}
\usepackage{pstricks-add}
\usepackage{shorttoc}
\usepackage{xypic}
\input cyracc.def
\font\tencyr=wncyr10          % tailles de 5 à 10
\def\cyr{\tencyr\cyracc}

\theoremstyle{plain}
\newtheorem{theor}{Théorème}[section]
\newtheorem{lemme}{Lemme}[section]

\newcommand{\Mod}[1]{\left({\rm{mod}}\ #1\right)}
\renewcommand{\theequation}{\thesection.\arabic{equation}}

\author{Kevin Destagnol}
\address{Institut de Mathématiques de Jussieu-Paris Rive Gauche\\ UMR 7586\\ Université Paris Diderot-Paris 7\\ Case postale 6052\\ Bâtiment Sophie Germain\\75205 Paris Cedex 13\\ France} 
\email{kevin.destagnol@imj-prg.fr}
\urladdr{webusers.imj-prg.fr/~kevin.destagnol/}
\date{27 décembre 2015}

\begin{document}
\renewcommand{\contentsname}{Sommaire}
\renewcommand{\abstractname}{Abstract}
\setcounter{section}{0}
\rhead{}
\chead{}
\lhead{}
\fancyhead[CO]{\footnotesize{ \textsc{La conjecture de Manin pour certaines surfaces de Châtelet}}}
\fancyhead[CE]{\footnotesize{\textsc{K. Destagnol}}}

\title{\textbf{La conjecture de Manin pour\\ certaines surfaces de Châtelet}}
\maketitle
\begin{abstract} 
Following the line of attack from La Bretèche, Browning and Peyre, we prove Manin's conjecture in its strong form conjectured by Peyre for a family of Châtelet surfaces which are defined as minimal proper smooth models of affine surfaces of the form
$$
Y^2-aZ^2=F(X,1),
$$
where $a=-1$, $F \in \mathbb{Z}[x_1,x_2]$ is a polynomial of degree 4 whose factorisation into irreducibles contains two non proportional linear factors and a quadratic factor which is irreducible over $\mathbb{Q}[i]$. This result deals with the last remaining case of Manin's conjecture for Châtelet surfaces with $a=-1$ and essentially settles Manin's conjecture for Châtelet surfaces with $a<0$.
\end{abstract}

 %\addtosections{toc}{\bfseries}
\setcounter{tocdepth}{2}

\begin{center}
\tableofcontents
\end{center}

\section{Introduction}
On définit la surface de Châtelet $S_{a,F}$ comme modèle minimal propre et 
lisse de variétés affines de $\mathbb{A}_{\mathbb{Q}}^3$ 
de la forme
$$X^2-aY^2=F(X,1)$$
où $F$ est une forme binaire à coefficients entiers de degré 4, de discriminant non nul et $a$ un entier qui n'est pas un carré. Ces surfaces introduites par Châtelet (\cite{Ch1} et \cite{Ch2}) sont arithmétiquement très riches puisqu'elles correspondent à l'exemple historique donné en 1971 par Swinnerton-Dyer de variétés ne satisfaisant ni le principe de Hasse ni l'approximation faible. Elles sont également des désingularisations minimales de surfaces de del Pezzo de degré 4 de type de singularité $2A_1$ conjuguées \cite[remarque 2.3]{35}.  Les surfaces de del Pezzo ont été parmi les premières à être considérées dans le traitement de la conjecture de Manin (\cite{20}, \cite{De}, \cite{DT} ou encore \cite{L} entre autres) puisque les plus simples dans la classification birationnelle d'Iskovskikh~\cite{I}.\par
L'objet de cet article est de décrire la répartition des points rationnels sur les surfaces de Châtelet dans les cas où $S_{a,F}:=S$ avec $F$ produit de deux formes linéaires $F_1$ et $F_2$ non proportionnelles
par une forme quadratique $F_3$ irréductible sur $\mathbb{Q}[i]$ et $a=-1$. Tous les éléments quant à la géométrie de ces surfaces dont nous auront besoin ont été démontrés dans \cite{CoSSWD1} et \cite{CoSSWD2} et sont résumés dans \cite{Brow} et \cite{35} dans le cas où $F$ est scindé. Ils s'adaptent au cas présent sans difficulté. En particulier, le système linéaire anticanonique $| \omega_S^{-1}|$ est sans point base, ce qui donne lieu à un morphisme $\psi:S \rightarrow \mathbb{P}^4$. On s'intéresse alors à la quantité
$$
N(B)=\#\{x\in S(\mathbb{Q}) \quad| \quad (H_4 \circ \psi)(x) \leqslant B\},
$$
où $H_4$ est une hauteur sur $\mathbb{P}^4(\mathbb{Q})$ qui sera définie en section 6.1. Le rang du groupe de Picard de $S$ vaut ici $2$ \cite{Brow} si bien que la conjecture de Manin prend la forme $$N(B) \sim c_S B \log(B)$$ pour une certaine constante $c_S>0$. Peyre \cite{P95} dans un premier temps, puis Batyrev et Tschinkel \cite{BatTsch} dans un cadre plus général, ont proposé une expression conjecturale de la constante $c_S$ s'exprimant en termes d'une mesure de Tamagawa sur l'espace adélique associé à la surface~$S$ (voir aussi la formule empirique donnée par Peyre dans \cite[formule 5.1]{P03}).\par
D'après \cite{CoSSWD1} et \cite{CoSSWD2}, l'obstruction de Brauer-Manin est la seule obstruction au principe de Hasse et à l'approximation faible pour les surfaces de Châtelet et en particulier les surfaces considérées spécifiquement dans cet article ne satisfont pas nécessairement l'approximation faible mais satisfont le principe de Hasse. Il s'agit ici du troisième cas, après \cite{35} et \cite{39}, pour lequel la conjecture de Manin dans sa forme forte conjecturée par Peyre dans \cite[formule 5.1]{P03} est établie par des méthodes de descente sur des torseurs pour une classe de variétés ne satisfaisant pas l'approximation faible.\par
Les travaux de La Bretèche, Browning et Peyre \cite{35} utilisant l'approche par les torseurs versels (introduits par Colliot-Thélène et Sansuc dans \cite{CoS1}, \cite{CoS2} et \cite{CoS3} et utilisés dans ce cadre pour la première fois par Salberger dans \cite{Sal}) ont permis dans un premier temps d'établir la conjecture de Manin dans le cas de factorisations de type un produit de formes linéaires non proportionnelles \cite{35}, puis les travaux de La Bretèche et Browning ont fourni le cas d'un produit d'une forme linéaire et d'une forme cubique irréductible sur $\mathbb{Q}$ \cite{36}, et enfin les travaux de La Bretèche et Tenenbaum ceux d'une forme irréductible sur $\mathbb{Q}[i]$ de degré 4 ou d'un produit de deux formes quadratiques irréductibles sur $\mathbb{Q}[i]$ \cite{39}.\par
On achève dans cet article la preuve de la conjecture de Manin dans le cas des surfaces de 
Châtelet avec $a=-1$. Le théorème suivant améliore notamment la borne $N(B) \ll B\log(B)$ obtenue par Browning dans \cite{Brow}.
\begin{theor}
Lorsque $a=-1$ et $F=F_1F_2F_3$ avec $F_1$ et $F_2$ deux formes linéaires non proportionnelles et $F_3$ une forme quadratique irréductible sur $\mathbb{Q}[i]$, on a
$$
N(B) \underset{B \rightarrow +\infty}{\sim} c_0 B \log(B),
$$ où la constante $c_0=c_S$ est la constante conjecturée par Peyre.
\label{theor1}
\end{theor}
\noindent
\textit{Remarque:} On peut se convaincre assez aisément que les méthodes développées ici permettraient de généraliser le résultat à tous les $a<0$ mais au prix d'un certain nombre de complications techniques comme expliqué dans \cite[Introduction]{39}. Le cas $a>0$ semble nécessiter une approche différente du fait du nombre infini d'unités du corps quadratique réel $\mathbb{Q}\left(\sqrt{a}\right)$.\\
\newline
\indent
La méthode utilisée ici suit celle développée par La Bretèche, Browning et Peyre dans le cas précédem\-ment cité. En particulier, la résolution de ce problème repose essentiellement sur l'estimation asymptotique de sommes du type
\begin{eqnarray}
S(X)=\sum_{\mathbf{x} \in \mathbb{Z}^2 \cap X\mathcal{R}} r\big(F_1(\mathbf{x})\big)r\big(F_2(\mathbf{x})\big)r\big(F_3(\mathbf{x})\big),
\label{sx1}
\end{eqnarray}
pour une région $\mathcal{R} \subset \mathbb{R}^2$ convenable et où $r(n)$ désigne le nombre de représentations d'un entier $n$ comme somme de deux carrés. On rappelle qu'on a l'expression
$$
r(n)=4 \sum_{d|n} \chi(d),
$$
où $\chi$ désigne le caractère de Dirichlet non principal modulo 4. La méthode pour estimer ces sommes est inspirée de l'article \cite{30} où le même genre de sommes est étudiée mais pour la fonction $\tau$ nombre de diviseurs à la place de la fonction $r$. Pour la vérification de la conjecture de Peyre, on s'inspire ici de \cite{39} et de \cite{35}. Cependant, à la différence du cas traité dans~\cite{35}, on établit que le traitement de la conjecture de Manin se fait au moyen d'un passage sur des torseurs distincts des torseurs versels dont on donne une description précise en section 7.\\
\newline
\noindent
\textbf{Remerciements.}-- L'auteur tient ici à exprimer toute sa gratitude à son directeur de thèse Régis de la Bretèche pour ses conseils, son soutien et ses relectures tout au long de ce travail ainsi qu'à Tim Browning,  Cyril Demarche, Ulrich Derenthal, Dan Loughran, Emmanuel Peyre et Marta Pieropan pour de nombreuses discussions éclairantes.

\section{Estimations de sommes liées à la fonction $r$}

Plus précisément, on considère $\mathbf{x}=(x_1,x_2)$, $F_1$ et $F_2$ deux formes linéaires binaires dans~$\mathbb{Z}[x_1,x_2]$, $F_3$ une forme quadratique dans $\mathbb{Z}[x_1,x_2]$ et $\mathcal{R}$ est une région de $\mathbb{R}^2$ vérifiant les hypothèses suivantes que l'on notera \textbf{NH}:
\begin{enumerate}[i)]
 \item 
 Les formes $F_1$ et $F_2$ ne sont pas proportionnelles;
 \item
 La forme quadratique $F_3$ est irréductible sur $\mathbb{Q}[i]$;
 \item
 $\forall \mathbf{x} \in \mathcal{R}, \quad F_i(\mathbf{x})>0;$
\hfill (2.0)
 \item
 La région $\mathcal{R}$ est convexe, bornée avec une frontière continûment différentiable.
\end{enumerate}
On définit
$$
X \mathcal{R}=\{ X \mathbf{x} \quad | \quad \mathbf{x} \in \mathcal{R} \}
$$
et on pose
$$
F_1(\mathbf{x})=a_1x_1+b_1x_2, \quad F_2(\mathbf{x})=a_2x_1+b_2x_2 \quad \mbox{et} \quad F_3(\mathbf{x})=a_3x_1^2+b_3x_2^2+c_3x_1x_2,
$$
avec les $a_i$, $b_i$ et $c_3$ entiers. Lorsque $J, L \in \mathbb{Z}[x_1,x_2]$ sont deux formes binaires homogènes, on note $\mbox{disc}(J)$ le discriminant de~$J(x,1)$ et $\mbox{Res}(J,L)$ le résultant de $J(x,1)$ et de $L(x,1)$. On considère alors
\begin{eqnarray}
\Delta=\mbox{disc}(F_3)=c_3^2-4a_3b_3 \neq 0, \quad \Delta_{12}=\mbox{Res}(F_1,F_2)=a_1b_2-a_2b_1 \neq 0,
\label{disc}
\end{eqnarray}
et
$$
\Delta_{i3}=\mbox{Res}(F_i,F_3)=a_3b_i^2+b_3a_i^2-c_3a_ib_i=F_3(-b_i,a_i) \neq 0,
$$
pour $i \in \{1,2\}$
et, pour $\mathbf{d}=(d_1,d_2,d_3) \in \mathbb{N}^3$,
$$
\Lambda(\mathbf{d})=\left\{ \mathbf{x} \in \mathbb{Z}^2 \quad | \quad d_i |F_i(\mathbf{x}) \right\},
$$
où l'on note systématiquement (sauf mention contraire) $d_i|F_i(\mathbf{x})$ pour $d_1|F_1(\mathbf{x})$, $d_2|F_2(\mathbf{x})$ et~$d_3 | F_3(\mathbf{x})$. On introduit alors la quantité essentielle
$$
 \rho(\mathbf{d})=\rho(\mathbf{d},F_1,F_2,F_3)=\#\left( \Lambda(\mathbf{d}) \cap [0,d[^2 \right)
$$
où l'on note systématiquement $d=d_1d_2d_3.$\\
\newline
On pose également
$$
 L_{\infty}=L_{\infty}(F_1,F_2,F_3)=\max\{\parallel F_i\parallel \},
$$
où $\parallel F_i \parallel$ désigne le plus grand coefficient en valeur absolue de la forme $F_i$,
$$
 r_{\infty}=r_{\infty}\left(\mathcal{R}\right)= \sup_{\mathbf{x} \in \mathcal{R}} \max \{|x_1|,|x_2|\}
$$
et
$$
 r'=r'(F_1,F_2,F_3,\mathcal{R})=\sup_{\mathbf{x} \in \mathcal{R}} \max\left\{|F_1(\mathbf{x})|,|F_2(\mathbf{x})|,\sqrt{|F_3(\mathbf{x})|}\right\}.
$$
On note enfin
$$
\mathcal{E}=\Big\{ n \in \mathbb{Z} \quad | \quad \exists \ell \in \mathbb{N}, \quad n \equiv 2^{\ell} \Mod{2^{\ell+2}} \Big\}
$$
ainsi que $\mathcal{E}_{2^n}$ sa projection modulo $2^n$
$$
\mathcal{E}_{2^n}=\left\{ k \in \mathbb{Z}/2^n\mathbb{Z} \quad | \quad \exists \ell \in \mathbb{N}, \quad k \equiv 2^{\ell} \Mod{2^{\min\{\ell+2,n\}}} \right\},
$$
et
\begin{eqnarray}
\eta=1-\frac{1+\log\log(2)}{\log(2)}=0.086071...
\label{eta}
\end{eqnarray}
La borne
$
S(X) \ll  X^2
$
fournie par \cite{B}, où la quantité $S(X)$ a été définie en (\ref{sx1}), est alors améliorée par le résultat suivant.
\begin{theor}
On suppose que les formes $F_1$, $F_2$ et $F_3$ et la région $\mathcal{R}$ vérifient les hypothèses \textbf{NH} (2.0) et soient $\varepsilon>0$ et $X \geqslant 1$ tels que $r'X^{1-\varepsilon} \geqslant 1$. On a alors
$$
S(X)=\pi^3 {\rm{vol}}(\mathcal{R}) X^2 \prod_{p} \sigma_p+O_{\varepsilon}\left( \frac{L_{\infty}^{\varepsilon}(r_{\infty}r'+r_{\infty}^2)X^2}{\left(\log(X)\right)^{\eta-\varepsilon}}\right),
$$
où
$$
 \sigma_p=\left( 1-\frac{\chi(p)}{p} \right)^3 \sum_{\boldsymbol{\nu} \in \mathbb{N}^3} \frac{\chi(p)^{\nu_1+\nu_2+\nu_3}\rho\left( p^{\nu_1},p^{\nu_2},p^{\nu_3} \right)}{p^{2(\nu_1+\nu_2+\nu_3)}}
$$
lorsque $p$ est impair et
$$
 \sigma_2=8\lim_{n \rightarrow +\infty} 2^{-2n} \#\left\{ 
 \mathbf{x} \in \left(\mathbb{Z}/2^n\mathbb{Z}\right)^2 \quad \bigg| \quad 
  F_i(\mathbf{x}) \in \mathcal{E}_{2^n}\\
 \right\}.
$$
Le produit $\underset{p>2}{\prod} \sigma_p$ est bien absolument convergent.
\label{theor2}
\end{theor}
On introduit alors l'ensemble
$$
 \mathfrak{D}=\{ (\mathbf{d},\mathbf{D}) \in \mathbb{N}^6 \quad | \quad d_i|D_i, \quad 2 \nmid D_i \}
$$
et pour $(\mathbf{d},\mathbf{D}) \in \mathfrak{D}$, $X \geqslant 1$, on pose
\begin{eqnarray}
 S(X,\mathbf{d},\mathbf{D})=S(X,\mathbf{d},\mathbf{D};\mathcal{R},F_1,F_2,F_3)=\sum_{\mathbf{x} \in \Lambda(\mathbf{D}) \cap X \mathcal{R}} r\left(\frac{F_1(\mathbf{x})}{d_1}\right)r\left(\frac{F_2(\mathbf{x})}{d_2}\right)r\left(\frac{F_3(\mathbf{x})}{d_3}\right).
\label{Ss}
\end{eqnarray}
Cette somme est directement liée à un problème de comptage sur la variété affine de $\mathbb{A}^8$ d'équations
\begin{eqnarray}
F_i(\mathbf{x})=d_i(s_i^2+t_i^2), \qquad (i=1,2,3)
\label{var}
\end{eqnarray}
où les $(x_1,x_2)$ sont restreints à une certaine région.\\
 \par
On introduit ensuite les entiers $\ell_1$, $\ell_2$, $\ell_3$ et les formes primitives $F_1^{\ast}$, $F_2^{\ast}$ et $F_3^{\ast}$ telles que
\begin{eqnarray}
 F_i=\ell_iF_i^{\ast}.
 \label{not}
\end{eqnarray}
On considère alors 
$$
 \mathbf{D'}=\left( \frac{D_1}{(D_1,\ell_1)},\frac{D_2}{(D_2,\ell_2)},\frac{D_3}{(D_3,\ell_3)} \right),
$$
\begin{eqnarray}
 a(\mathbf{D},\mathbf{\Delta})=(D_1,\Delta_{12}\Delta_{13})(D_2,\Delta_{12}\Delta_{23})(D_3,\Delta(\Delta_{13},\Delta_{23}))
 \label{a}
\end{eqnarray}
avec $\mathbf{\Delta}=(\Delta,\Delta_{12},\Delta_{13},\Delta_{23})$ et 
$$
 a'(\mathbf{D},\mathbf{\Delta})=a(\mathbf{D}',\mathbf{\Delta}')=(D'_1,\Delta'_{12}\Delta'_{13})(D'_2,\Delta'_{12}\Delta'_{23})(D'_3,\Delta'(\Delta'_{13},\Delta'_{23}))
$$
où 
$$
\mathbf{\Delta}'=(\Delta',\Delta'_{12},\Delta'_{13},\Delta'_{23})=\left(\frac{\Delta}{\ell_3^2},\frac{\Delta_{12}}{\ell_1\ell_2},\frac{\Delta_{13}}{\ell_1^2\ell_3},\frac{\Delta_{23}}{\ell_2^2\ell_3}\right).
$$
On a alors le résultat suivant, crucial en vue de l'obtention de la conjecture de Manin.
\begin{theor}
Soient $\varepsilon>0$ et $X \geqslant 1$ tels que $r'X^{1-\varepsilon} \geqslant 1$. Si on suppose que les formes $F_1$, $F_2$ et $F_3$ et la région $\mathcal{R}$ vérifient les hypothèses \textbf{NH} (2.0) et que $(\mathbf{d},\mathbf{D}) \in \mathfrak{D}$, alors
$$
S(X,\mathbf{d},\mathbf{D})= \pi^3 {\rm{vol}}(\mathcal{R}) X^2 \prod_{p} \sigma_p(\mathbf{d},\mathbf{D})+O_{\varepsilon}\left( \frac{L_{\infty}^{\varepsilon}D^{\varepsilon}(r_{\infty}r'+r_{\infty}^2)a'(\mathbf{d},\mathbf{\Delta})X^2}{\left(\log(X)\right)^{\eta-\varepsilon}}\right),
$$
où
$$
 \sigma_p(\mathbf{d},\mathbf{D})=\left( 1-\frac{\chi(p)}{p} \right)^3 \sum_{\boldsymbol{\nu} \in \mathbb{N}^3} \frac{\chi(p)^{\nu_1+\nu_2+\nu_3}\rho\left( p^{N_1},p^{N_2},p^{N_3} \right)}{p^{2(N_1+N_2+N_3)}}
$$
avec $N_i=\max\{\nu_p(D_i),\nu_i+\nu_p(d_i)\}$ lorsque $p>2$ et
\begin{eqnarray}
 \sigma_2(\mathbf{d},\mathbf{D})=\sigma_2(\mathbf{d})=8\lim_{n \rightarrow +\infty} 2^{-2n} \#\left\{ 
 \mathbf{x} \in \left(\mathbb{Z}/2^n\mathbb{Z}\right)^2 \quad \bigg| \quad 
  F_i(\mathbf{x}) \in d_i\mathcal{E}_{2^n}
 \right\}.
\label{sigmadeux}
\end{eqnarray}
De plus, on a 
$$
 \prod_{p} \sigma_p(\mathbf{d},\mathbf{D}) \ll L_{\infty}^{\varepsilon}D^{\varepsilon}a'(\mathbf{d},\mathbf{\Delta}).
$$
\label{theor3}
\end{theor}
\noindent
\textit{Remarque}: Le même raisonnement que celui mené dans la section 6 de \cite{30} permet, \textit{mutatis mutandis}, de se ramener au Théorème \ref{theor2} grâce à un changement de variables. Le fait que~$\Lambda(\mathbf{D})$ ne soit pas un réseau (contrairement au cas de \cite{35}) est pallié en le reliant à une réunion de réseaux de la même manière que dans \cite{D99}. Ainsi, on ne détaille pas ici la démonstration du Théorème \ref{theor3}.\\
\newline
\indent
Pour terminer cette section, dans l'optique de vérifier la conjecture de Peyre, il est bon de réinterpréter la constante obtenue dans le Théorème \ref{theor3}. On définit pour $\boldsymbol{\lambda}\in \mathbb{N}^3$, $\boldsymbol{\mu} \in \mathbb{N}^3$ et $p$ premier différent de 2:
\begin{eqnarray}
 N_{\boldsymbol{\lambda},\boldsymbol{\mu}}(p^n)=\#\left\{  (\mathbf{x},\mathbf{s},\mathbf{t}) \in \left(\mathbb{Z}/p^n\mathbb{Z}\right)^8 \quad \Bigg| \quad 
 \begin{array}{l}
  F_i(\mathbf{x}) \equiv p^{\lambda_i}(s_i^2+t_i^2)\Mod{p^n},\\
  p^{\mu_i}| F_i(\mathbf{x})\\
 \end{array}
 \right\}
\label{N2}
\end{eqnarray}
et
$$
 \omega_{\boldsymbol{\lambda},\boldsymbol{\mu}}(p)=\lim_{n \rightarrow +\infty}p^{-5n-\lambda_1-\lambda_2-\lambda_3}N_{\boldsymbol{\lambda},\boldsymbol{\mu}}(p^n).
$$
Ces quantités sont bien définies et correspondent aux densités $p$-adiques associées à la variété définie par les équations (\ref{var}). Pour le cas $p=2$, on introduit
$$
N_{\mathbf{d}}(2^n)=\#\left\{  (\mathbf{x},\mathbf{s},\mathbf{t}) \in \left(\mathbb{Z}/2^n\mathbb{Z}\right)^8 \quad \bigg| \quad 
  F_i(\mathbf{x}) \equiv d_i(s_i^2+t_i^2)\Mod{2^n}
 \right\}
$$
et
\begin{eqnarray}
 \omega_{\mathbf{d}}(2)=\lim_{n \rightarrow +\infty}2^{-5n}N_{\mathbf{d}}(2^n).
\label{omegadeux}
\end{eqnarray}
Posant enfin $\omega_{\infty}(\mathcal{R})$ la densité archimédienne associée au problème de comptage (\ref{var}), on obtient le résultat suivant, qui traduit en réalité le fait que le principe de Hasse est vérifié sur les torseurs décrits en section 7.
\begin{theor}
Supposons que les formes $F_1$, $F_2$ et $F_3$ et la région $\mathcal{R}$ vérifient les hypothèses \textbf{NH} (2.0) et que $(\mathbf{d},\mathbf{D}) \in \mathfrak{D}$. Pour tout premier $p>2$, on a $\omega_{\boldsymbol{\lambda},\boldsymbol{\mu}}(p)=\sigma_p(\mathbf{d},\mathbf{D})$ avec~$\boldsymbol{\lambda}=(\nu_p(d_1),\nu_p(d_2),\nu_p(d_3))$ et~$\boldsymbol{\mu}=(\nu_p(D_1),\nu_p(D_2),\nu_p(D_3))$, $\omega_{\mathbf{d}}(2)=\sigma_2(\mathbf{d})$ et $\omega_{\infty}(\mathcal{R})= \pi^3 \rm{vol}(\mathcal{R})$.
\label{theor4}
\end{theor}

\section{Propriétés des fonctions $\rho$}

 La fonction $\rho$ a été abondamment étudiée, notamment dans \cite{D99}, \cite{M1} et \cite{M2} ou \cite{30}. Il résulte du théorème chinois que la fonction $\rho$ est multiplicative dans le sens où si l'on se donne deux triplets $\mathbf{d}=(d_1,d_2,d_3)$ et $\mathbf{d'}=(d'_1,d'_2,d'_3)$ tels que $(d_1d_2d_3,d'_1d'_2d'_3)=1$, alors on a
$$
\rho(d_1d'_1,d_2d'_2,d_3d'_3)=\rho(\mathbf{d})\rho(\mathbf{d'}).
$$
Il suffit donc de l'étudier sur les triplets de nombres premiers.  Le lemme immédiat suivant tiré de \cite[lemme 1]{30} permet de se ramener à ce cas.
\begin{lemme}
Soient $F_1$, $F_2 \in \mathbb{Z}[x_1,x_2]$ des formes linéaires non proportionnelles et $F_3\in \mathbb{Z}[x_1,x_2]$ une forme quadratique. Avec les notations de (\ref{not}), on a pour $\mathbf{d} \in \mathbb{N}^3$,
$$
 \frac{\rho(\mathbf{d},F_1,F_2,F_3)}{(d_1d_2d_3)^2}=\frac{\rho(\mathbf{d'},F_1^{\ast},F_2^{\ast},F_3^{\ast})}{(d'_1d'_2d'_3)^2}
$$
où $\mathbf{d'}=\left( \frac{d_1}{(d_1,\ell_1)},\frac{d_2}{(d_2,\ell_2)},\frac{d_3}{(d_3,\ell_3)} \right).$
\label{lemme1}
\end{lemme}
\noindent
On a alors le résultat suivant.
\begin{lemme}
Soient $F_1$, $F_2$ des formes linéaires primitives non proportionnelles et $F_3$ une forme quadratique primitive irréductible sur $\mathbb{Q}$.
 \begin{enumerate}[a)]
  \item 
 Pour tout nombre premier $p$, on a
 $$
 \rho(p^{\nu},1,1)=\rho(1,p^{\nu},1)=p^{\nu}.
 $$
 \item
 Si $p \nmid 2{\rm{disc}}(F_3)$, alors
  $$
  \rho(1,1,p^{\nu})=\varphi(p^{\nu})\left( 1+\left( \frac{{\rm{disc}}(F_3)}{p} \right) \right)\left\lceil \frac{\nu}{2} \right\rceil+p^{2(\nu-\lceil \nu/2\rceil)} \leqslant (\nu+1)p^{\nu}.
  $$
  De plus, si $p$ est un facteur impair de ${\rm{disc}}(F_3)$, on a 
  $$
  \rho(1,1,p^{\nu}) \ll (\nu+1)p^{\nu+\min\{\left\lfloor\nu_p(\footnotesize{{\rm{disc}}}(F_3)/2)\right\rfloor,\left\lfloor\nu/2\right\rfloor\}}
  $$
  et si $p=2$,
  $$
  \rho(1,1,2^{\nu}) \ll (\nu+1)2^{\nu}.
  $$
  \item
  Lorsque $\max\{\nu_1,\nu_2\} \leqslant \left \lceil \frac{\nu_3}{2} \right\rceil$ et $p$ impair, on a
  $$
  \rho(p^{\nu_1},p^{\nu_2},p^{\nu_3}) \ll (\nu_3+1)p^{2(\nu_1+\nu_2)+\nu_3}p^{\min\{\left\lfloor\nu_p(\footnotesize{{\rm{disc}}}(F_3)/2)\right\rfloor,\left\lfloor\nu_3/2\right\rfloor\}}.
  $$
  De plus, lorsque $\max\{\nu_1,\nu_2\}=\nu_3=1$ et $p\nmid\Delta_{12}\Delta_{13}\Delta_{23}$, on a
  $$
  \rho(p^{\nu_1},p^{\nu_2},p^{\nu_3}) \ll p^{2(\nu_1+\nu_2)}.
  $$
  Enfin, lorsque $\max\{\nu_1,\nu_2\} \leqslant \left \lceil \frac{\nu_3}{2} \right\rceil$, on a
  $$
  \rho(2^{\nu_1},2^{\nu_2},2^{\nu_3}) \ll 2^{2(\nu_1+\nu_2)+\nu_3}.
  $$
  \item
  Lorsque $\left \lceil \frac{\nu_3}{2} \right\rceil \leqslant \min\{\nu_1,\nu_2\}$ et pour tout nombre premier, on a
  $$
   \rho(p^{\nu_1},p^{\nu_2},p^{\nu_3}) \ll p^{\nu_1+\nu_2+2\nu_3+\min\{\nu_1,\nu_2,\nu_p(\Delta_{12})\}}.
  $$
  \item
  Lorsque $\nu_j \leqslant \left \lceil \frac{\nu_3}{2} \right\rceil \leqslant \nu_3 \leqslant \nu_i$ avec $\{i,j\}=\{1,2\}$ et pour tout nombre premier $p$, on a
  $$
  \rho(p^{\nu_1},p^{\nu_2},p^{\nu_3}) \ll \nu_3 p^{2 \nu_j+\nu_i+\nu_3+\rfloor\nu_3/2\rfloor}\left( p^{\min\{\lceil \nu_3/2 \rceil,\lceil \nu_p(\Delta_{i3})/2 \rceil\}}+p^{r_p} \right),
  $$
  où
  $$
  r_p=\min\left\{\nu_i-\lceil \nu_3/2 \rceil,\lceil \nu_p(\Delta_{i3})/2 \rceil\right\}+\min\left\{\left\lceil\nu_3/2\right\rceil,\left\lceil\nu_p({\rm{disc}}(F_3))/2\right\rceil\right\}.
  $$
    \item
  On a pour tout nombre premier $p$, la majoration
  $$
  \rho(p^{\nu_1},p^{\nu_2},p^{\nu_3}) \ll \min\left\{ p^{2(\nu_2+\nu_3)+\nu_1},p^{2(\nu_1+\nu_3)+\nu_2},(\nu_3+1)p^{2(\nu_1+\nu_2)+3\frac{\nu_3}{2}} \right\}.
  $$
 \end{enumerate}
\label{lemme2}
\end{lemme}
\noindent
\textit{Démonstration.}-- La preuve des points a) à e) se trouve dans le lemme 3 de \cite{30}. Démontrons donc le point~f). En négligeant les deux conditions 
$$
p^{\nu_1} | F_1(\mathbf{x}) \quad \mbox{et} \quad p^{\nu_3} | F_3(\mathbf{x}),
$$
on obtient les majorations
$$
\rho(p^{\nu_1},p^{\nu_2},p^{\nu_3}) \leqslant p^{2(\nu_1+\nu_3)}\rho(1,p^{\nu_2},1)=p^{2(\nu_1+\nu_3)+\nu_2}.
$$
On obtient aussi la majoration
$$
\rho(p^{\nu_1},p^{\nu_2},p^{\nu_3}) \leqslant p^{2(\nu_2+\nu_3)+\nu_1}
$$
en inversant les rôles de $F_1$ et $F_2$. En oubliant les deux conditions sur les formes linéaires, on obtient
$$
\rho(p^{\nu_1},p^{\nu_2},p^{\nu_3}) \leqslant p^{2(\nu_1+\nu_2)}\rho(1,1,p^{\nu_3})
$$
\nopagebreak[0]
et on conclut grâce au point b).
\hfill
$\square$\\
\indent
Suivant toujours la section 3 de \cite{30}, on peut déduire de ces deux lemmes que la fonction
$$
f(\mathbf{d})=\frac{\rho(\mathbf{d})}{d_1d_2d_3}
$$
est proche, au sens de la convolution, de la fonction $R:\mathbf{d} \mapsto R(\mathbf{d})=r_{\Delta}(d_3)$ avec
\begin{eqnarray}
r_{\Delta}(\ell)=\sum_{k|\ell} \chi_{\Delta}(k)
\label{rrrr}
\end{eqnarray}
où $\chi_{\Delta}(n)=\left(\frac{\Delta}{n} \right)$ est le symbole de Kronecker et $\Delta$ a été défini en (\ref{disc}).  On pose également~$h$ la fonction arithmétique satisfaisant
$$
f(\mathbf{d})=(h \ast R) (\mathbf{d})=\sum_{\mathbf{k} \in \mathbb{N}^3 \atop k_i|d_i}h\left(\frac{d_1}{k_1},\frac{d_2}{k_2},\frac{d_3}{k_3}  \right)R(\mathbf{k}).
$$
\begin{lemme}
Soient $F_1$ et $F_2$ deux formes linéaires non proportionnelles et $F_3$ une forme quadratique irréductible sur $\mathbb{Q}$. Avec la notation (\ref{rrrr}), on a pour tout $A>0$,
$$
\sum_{\mathbf{k} \in \mathbb{N}^3} \frac{|h(\mathbf{k})|\log(k_1k_2k_3)^A}{k_1 k_2 k_3} \ll L_{\infty}^{\varepsilon}.
$$
En particulier,
\begin{eqnarray}
\prod_{p>2} \sigma_p=L(1,\chi_{\Delta}\chi)\sum_{\mathbf{k} \in \mathbb{N}^3} \frac{h(\mathbf{k})\chi(k_1k_2k_3)}{k_1k_2k_3} \ll L_{\infty}^{\varepsilon}.
\label{convol}
\end{eqnarray}
\label{lemme3}
\end{lemme}
\noindent
\textit{Démonstration.}-- Il suffit de remarquer que la preuve du lemme 4 de \cite{30} dans le cas $A=1$ s'adapte immédiatement pour fournir le résultat.
\hfill
$\square$\\
\newline
\indent
On définit également pour un polynôme $g \in \mathbb{Z}[X]$ la quantité
$$
 \rho_g(n)=\#\left\{ x \in \mathbb{Z}/n\mathbb{Z} \quad | \quad g(x) \equiv 0\Mod{n} \right\}.
$$
Il s'agit d'une fonction multiplicative qu'il suffit donc d'étudier sur les puissances de nombres premiers. On aura besoin du résultat suivant (voir \cite[lemme 1]{36} par exemple).
\begin{lemme}
Soient $g \in \mathbb{Z}[X]$ de degré $d \geqslant 2$ et $p$ un nombre premier qui ne divise pas le contenu de $g$ et tel que $p^{\mu} \parallel  {\rm{disc}}(g)$. Alors, pour tout $\nu \geqslant 1$, on a
$$
\rho_g\left( p^{\nu}\right) \leqslant d\min\left\{ p^{\frac{\mu}{2}},p^{\left( 1-\frac{1}{d} \right)\nu},p^{\nu-1} \right\}.
$$
\label{lemme4}
\end{lemme}

\section{Démonstration du Théorème \ref{theor2}}
\subsection{Extraction des valuations 2-adiques}
On effectue ici un raisonnement préliminaire similaire à celui effectué à la section 3 de \cite{36} afin de pouvoir utiliser des décompositions semblables à celles de la section 3 de \cite{HB03}. En utilisant les formules $r(2n)=r(n)$ et~$r(n)=0$ si~$n \equiv 3 \Mod{4}$, on obtient l'égalité
$$
S(X)=\sum_{k_0 \geqslant 0} \sum_{\mathbf{x} \in \mathbb{Z}^2 \cap X\mathcal{R} \atop 2^{k_0}\parallel \mathbf{x}} r\big( F_1(\mathbf{x})\big) r\big( F_2(\mathbf{x})\big) r\big( F_3(\mathbf{x})\big)=\sum_{k_0 \geqslant 0} S^{\ast}\left( 2^{-k_0}X \right),
$$
où 
$$
S^{\ast}(X)=\sum_{\mathbf{x} \in \mathbb{Z}^2 \cap X\mathcal{R} \atop 2 \nmid(x_1,x_2)}r\big( F_1(\mathbf{x})\big) r\big( F_2(\mathbf{x})\big) r\big( F_3(\mathbf{x})\big). 
$$
En regroupant les termes selon la valuation $2-$adique de $F_i(\mathbf{x})$, on a
$$
 S(X)=\sum_{k_0 \geqslant 0} \sum_{\mathbf{k} \in \mathbb{N}^3} S_{\mathbf{k}}\left( 2^{-k_0}X \right),
$$
avec $S_{\mathbf{k}}(X)$ la restriction de $S^{\ast}(X)$ aux $\mathbf{x}$ tels que 
$$
\nu_2(F_i(\mathbf{x}))=k_i, \quad 2^{-k_i}F_i(\mathbf{x}) \equiv 1\Mod{4}
$$
et $2 \nmid \mathbf{x}$. On a clairement que $2^{k_i}\leqslant F_i(\mathbf{x}) \leqslant (X')^{{\rm{deg}}(F_i)}$ donc
$
k_i \ll \log(X').
$
On remarque également que $\min\{k_i,k_j\}\leqslant\nu_2\left( \Delta_{ij} \right)$.

\begin{lemme}
Lorsque les conditions 
$$
\nu_2(F_i(\mathbf{x}))=k_i \quad \mbox{et} \quad 2^{-k_i}F_i(\mathbf{x}) \equiv 1\Mod{4}
$$
et $2 \nmid \mathbf{x}$ sont réalisées, il existe une matrice $\mathbf{M} \in \mathcal{M}_2(\mathbb{Z})$ inversible dans $\mathbb{Q}$ telle que les solutions s'écrivent $\mathbf{x}=\mathbf{Mx'}$ avec $x'_1 \equiv 1\Mod{4}$.
\label{lemme5}
\end{lemme}

\noindent
\textit{Démonstration.}-- 
En raisonnant comme dans la section 3 de \cite{36}, on peut supposer sans perte de généralité que $a_1$ est impair.  La première condition 
\begin{eqnarray}
2^{-k_1}F_1(\mathbf{x}) \equiv 1\Mod{4}
\label{cccc}
\end{eqnarray}
équivaut alors à l'existence de $x'_1 \equiv 1\Mod{4}$ tel que $x_1=cx_2+c'2^{k_1}x'_1$ où $c' \in\{\pm 1\}$ tel que $c' \equiv a_1\Mod{4}$ et $c \in \left[ 0, 2^{k_1+2} \right[\cap \mathbb{Z}$ tel que $a_1c \equiv-b_1 \Mod{2^{k_1+2}}$. Si $k_1=0$, on a alors automatiquement que $2 \nmid \mathbf{x}$ et sinon, cette condition devient équivalente au fait que $x_2$ soit impair.\par
Intéressons-nous à présent à la deuxième condition
$$
2^{-k_2}F_2(\mathbf{x}) \equiv 1 \equiv x'_1\Mod{4}.
$$
Si on suppose la première condition vérifiée, cette condition est équivalente à
$$
F_{2}\left( cx_2+c'2^{k_1}x'_1,x_2 \right) \equiv 2^{k_2} x'_1 \Mod{2^{k_2+2}}.
$$
On considère alors $F'_2(X,Y)=F_2\left( cY+c'2^{k_1}X,Y \right)=2^{k'_2}(aX+bY)$ pour $k'_2$, $a$ et $b$ trois entiers tels que $(2,a,b)=1$. Soit $k'_2 > k_2$ et alors les deux premières conditions n'ont pas de solution commune. Soit $k'_2 \leqslant k_2$ et on pose
\begin{eqnarray}
k''_2=k_2-k'_2 \geqslant 0
\label{kseconde}
\end{eqnarray} 
et $F''_2(X,Y)=2^{-k'_2}F'_2(X,Y)$ de sorte que
$$
F''_2(x'_1,x_2) \equiv 2^{k''_2}x'_1 \Mod{2^{k''_2+2}}.
$$
On peut écrire $x_2 \equiv \alpha_2 x'_1\Mod{2^{k''_2+2}}$ pour un unique $\alpha_2 \in \left[ 0, 2^{k''_2+2} \right[\cap \mathbb{Z}$ pour obtenir
\begin{eqnarray}
F''_2(1,\alpha_2) \equiv 2^{k''_2} \Mod{ 2^{k''_2+2} }.
\label{cond1}
\end{eqnarray}
Finalement, lorsque la première condition est vérifiée, la deuxième est vérifiée si, et seulement si, $k'_2 \leqslant k_2$ et $x_2=\alpha_2 x'_1+2^{k''_2+2}x'_2$ pour $\alpha_2$ solution de (\ref{cond1}). De plus, si $k_1\neq 0$, $\alpha_2$ doit être choisi impair.\par
Posons $F'_3(X,Y)=F_3\left( cY+c'2^{k_1}X,Y \right)=2^{k'_3}(cX^2+dXY+eY^2)$ pour $k'_3$, $c$, $d$ et $e$ quatre entiers tels que $(2,c,d,e)=1$ et $k''_3=k_3-k'_3$ lorsque $k'_3 \leqslant k_3$. L'égalité
$$
2^{-k_3}F_3(\mathbf{x}) \equiv 1 \equiv {x'_1}^2\Mod{4}
$$
conduit alors de même à montrer que si la première condition est remplie, la troisième l'est si, et seulement si, $k'_3 \leqslant k_3$ et s'il existe $\alpha_3 \in \left[ 0, 2^{k''_3+2} \right[\cap \mathbb{Z}$ solution de 
\begin{eqnarray}
 F_3''(1,\alpha_3) \equiv 2^{k''_3} \Mod{ 2^{k''_3+2} }
 \label{cond2}
\end{eqnarray}
tel que $x_2=\alpha_3 x'_1+2^{k''_3+2}x''_2,$ $\alpha_3$ devant être choisi impair lorsque $k_1 \neq 0$. Dans le cas où les trois conditions sont remplies simultanément, d'après ce qui précède, soit il n'existe pas de solution soit il existe $\alpha_2$ solution de (\ref{cond1}), $\alpha_3$ solution de (\ref{cond2}) tels que
$$
x_2\equiv \alpha_2 x'_1\Mod{2^{k''_2+2}} \quad \mbox{et} \quad x_{2}\equiv \alpha_3 x'_1\Mod{2^{k''_3+2}}.
$$
On a donc
$$
\alpha_2 \equiv \alpha_3  \Mod{ 2^{\min\{k''_2,k''_3\}+2} }
$$
et par conséquent $x_2$ est de la forme
$
x_2=\alpha x'_1+2^{\max(k''_2,k''_3)+2}x'_2$ avec $\alpha=\alpha_2$ si $k_2 \geqslant k_3$ et $\alpha=\alpha_3$ si $k_3 \geqslant k_2$ et pour un certain entier $x'_2$.
On vient donc de montrer que si les conditions $$
\nu_2(F_i(\mathbf{x}))=k_i, \quad 2^{-k_i}F_i(\mathbf{x}) \equiv 1\Mod{4}
$$
et $2 \nmid \mathbf{x}$ admettent des solutions, ces dernières peuvent s'écrire $\mathbf{x}=\mathbf{Mx'}$ avec $x'_1 \equiv 1\Mod{4}$ et
$$
 \mathbf{M}=\mathbf{M}_{\alpha}=
 \left(
 \begin{array}{ll}
  c'2^{k_1}& c\\
  0&1\\ 
 \end{array}
 \right)
  \left(
 \begin{array}{ll}
  1& 0\\
  \alpha&2^{\max\{k''_2,k''_3\}+2}\\ 
 \end{array}
 \right)=  \left(
 \begin{array}{ll}
  c'2^{k_1}+c\alpha& c2^{\max\{k''_2,k''_3\}+2}\\
  \alpha&2^{\max\{k''_2,k''_3\}+2}\\ 
 \end{array}
 \right),
$$
où $\alpha \in \left[0, 2^{\max\{k''_2,k''_3\}+2} \right[\cap \mathbb{Z}$ vérifie
\begin{eqnarray}
\left\{
\begin{array}{l}
F''_2(1,\alpha) \equiv 2^{k''_2} \Mod{ 2^{k''_2+2} },\\
  F_3''(1,\alpha) \equiv 2^{k''_3} \Mod{ 2^{k''_3+2} }, \\
\alpha \equiv 1 \Mod{2} \quad \mbox{lorsque } \quad k_1 \geqslant 1.\\
\end{array}
\right.
\label{conds}
\end{eqnarray}
On remarque pour conclure que
$$
 |\det(\mathbf{M})|=2^{k_1+\max\{k''_2,k''_3\}+2}.
$$
\hfill
$\square$\par
On majore alors 
\begin{eqnarray}
n(\mathbf{k})=n(k_1,k_2,k_3)
\label{nk}
\end{eqnarray}
le nombre d'entiers $\alpha$ dans $\left[0, 2^{\max\{k''_2,k''_3\}+2} \right[\cap \mathbb{Z}$ vérifiant les conditions (\ref{conds}). 
\begin{lemme}
 On a $n(\mathbf{k}) \ll 1$ avec $n(\mathbf{k})$ défini en (\ref{nk}) et où la borne ne dépend que de la valuation $2$-adique de $\Delta^{\ast}:=\Delta_{12}\Delta_{13}\Delta_{23}\Delta$.
\label{lemme6}
\end{lemme}
\noindent
\textit{Démonstration.}-- On commence par traiter le cas où $k_1 \geqslant 1$ dans lequel on sait que $\alpha$ doit être impair. On a clairement la majoration suivante
$$
n(\mathbf{k}) \leqslant \# \left\{ x \in \mathbb{Z}/{2^{k_1+\max\{k_2,k_3\}+2}}\mathbb{Z} \quad \Bigg| \quad 
\begin{array}{lll}
x \equiv -b_1 \overline{a_1} \Mod{2^{k_1}}\\[1mm]
 F_2(x,1) \equiv 0 \Mod{2^{k_2}}\\[1mm]
  F_3(x,1) \equiv 0 \Mod{2^{k_3}}\\
\end{array} \right\},
$$
où $\overline{a}_1$ désigne l'inverse multiplicatif de $a_1$ dans $\mathbb{Z}/{2^{k_1}}\mathbb{Z}$.
On a alors, dans le cas où $k_1\geqslant \max\{k_2,k_3\}$, l'inégalité
$$
n(\mathbf{k}) \ll \frac{2^{k_1+\max\{k_2,k_3\}}}{2^{k_1}}=2^{\max\{k_2,k_3\}}.
$$
En utilisant le fait que $\min\{k_i,k_1\} \leqslant \nu_2(\Delta_{i1})$ pour $i \in \{2,3\}$, on en déduit la majoration
$$
n(\mathbf{k}) \ll 2^{\nu_2(\Delta^{\ast})} \ll 1.
$$
Passons aux cas où $\max\{k_2,k_3\}>k_1$. On traite tout d'abord le cas $k_3\geqslant k_2$. Par le Lemme~\ref{lemme4}, on a alors la majoration
$$
n(\mathbf{k}) \ll \frac{2^{k_1+k_3}}{2^{k_3}}\rho_{F_3(x,1)}\left( 2^{k_3} \right) \ll 2^{k_1}\rho_{F_3(x,1)}\left( 2^{k_3} \right) \ll 2^{k_1+\frac{\nu_2(\Delta)}{2}} \ll 1.
$$
Enfin, de la même manière, dans les cas où $k_2\geqslant k_3$, on obtient
$$
n(\mathbf{k}) \ll 2^{k_1}\rho_{F_2(x,1)}\left( 2^{k_2} \right) \ll 2^{k_1} \ll 1.
$$
\par
Il reste à majorer $n(\mathbf{k})$ lorsque $k_1=0$, dans ce cas, on procède de même que ci-dessus pour les $\alpha$ impairs et on doit ajouter tous les $\alpha$ pairs de $\left[0, 2^{\max\{k_2,k_3\}} \right[\cap \mathbb{Z}$ qui vérifient (\ref{conds}). Avec la notation (\ref{cond2}), le nombre de ces $\alpha$ est majoré par celui des $\beta$ tels que
$$
F_3''(1,\beta) \equiv 0 \Mod{ 2^{k''_3} }.
$$
Par le Lemme \ref{lemme4}, cette quantité est majorée par $2^{\nu_2({\rm{disc}}(F_3''))} \ll 1$, ce qui achève la preuve.
\hfill
$\square$\\
\par
On relie maintenant la quantité $n(\mathbf{k})$ et la constante $\sigma_2$ qui apparaît dans le Théorème~\ref{theor2}.
\begin{lemme}
On a 
$$
 \sigma_2=2\sum_{k_0 \geqslant 0} \frac{1}{2^{2k_0}} \sum_{k_1,k_2,k_3 \in \mathbb{N}^3} \frac{n(\mathbf{k})}{2^{k_1+\max\{k''_2,k''_3\}+2}},
$$
où les entiers $k''_2$ et $k''_3$ se déduisent de $k_2$ et $k_3$ comme en (\ref{kseconde}) et (\ref{cond2}).
\label{lemme7}
\end{lemme}
\noindent
\textit{Démonstration.}-- En partitionnant $\left(\mathbb{Z}/2^{n}\mathbb{Z}\right)^2$ selon la valuation $2-$adique de $(x_1,x_2)$, 
on obtient l'égalité
$$
\sigma_2=\lim_{n \rightarrow +\infty}\frac{8}{2^{2n}} \sum_{0 \leqslant k_0\leqslant n\atop 0 \leqslant k_1,k_2,k_3 \leqslant n}\#\left\{ 
 \mathbf{x} \in \left(\mathbb{Z}/2^{n-k_0}\mathbb{Z}\right)^2 \hspace{0.5mm} \bigg| \hspace{0.5mm}  2\nmid \mathbf{x}, \hspace{0.5mm} 
  F_i(\mathbf{x}) \equiv 2^{k_i} \Mod{ 2^{\min\{k_i+2,n-k_0\}}}
 \right\}.
$$
On vérifie aisément grâce au Lemme \ref{lemme4} que la contribution des termes tels que pour au moins un des $k_i$ ($1 \leqslant i \leqslant 3$) on ait $k_i+2>n-k_0$ à la triple somme intérieure est négligeable. On peut donc s'intéresser uniquement aux triplets tels que $\max\{k_1,k_2,k_3\} \leqslant n-k_0-2$, ce qui revient à estimer
$$
\sigma_2= 8\lim_{n \rightarrow +\infty} 2^{-2n} \sum_{0 \leqslant k_0\leqslant n\atop0 \leqslant k_1,k_2,k_3 \leqslant n-k_0-2}  \#\left\{ 
 \mathbf{x} \in \left(\mathbb{Z}/2^{n-k_0}\mathbb{Z}\right)^2 \hspace{0.5mm} \bigg| \hspace{0.5mm}  2\nmid \mathbf{x}, \hspace{0.5mm} 
  F_i(\mathbf{x}) \equiv 2^{k_i} \Mod{2^{k_i+2} }
 \right\}.
$$
On a vu en (\ref{cccc}) lors de la preuve du Lemme \ref{lemme5} qu'on pouvait supposer $a_1$ impair et que les conditions
$$
\nu_2(F_i(\mathbf{x}))=k_i, \quad 2^{-k_i}F_i(\mathbf{x}) \equiv 1\Mod{4}, \quad 2 \nmid \mathbf{x}
$$
impliquent l'existence de $x'_1\equiv 1 \Mod{4}$ tel que
$$
x_1 \equiv cx_2+c'2^{k_1}x'_1\Mod{2^{k_1+2} }
$$
avec $c' \in \{\pm1\}$ et $c \equiv a_1 \Mod{4}$. Par conséquent, $x_1$ est entièrement déterminé par la valeur de $x_2$ modulo $2^{k_1+2}$. On sait ensuite que si on pose $x_2 \equiv \alpha x'_1 \Mod{2^{\max\{k''_2,k''_3\}+2} }$, alors $\alpha$ est une des $n(\mathbf{k})$ solutions de (\ref{conds}) dans $\left[  0,2^{\max\{k''_2,k''_3\}+2}\right[\cap\mathbb{Z}$. Posant $k=\max\{k_1,k''_2,k''_3\}$, on constate alors que les trois conditions de congruence ne dépendent que de la classe de $\mathbf{x}$ modulo $2^{k+2}$. Autrement dit, on obtient 
$$
\begin{aligned}
\sigma_2= &8\lim_{n \rightarrow +\infty} 2^{-2n} \sum_{0 \leqslant k_0\leqslant n \atop 0 \leqslant k_1,k_2, k_3 \leqslant n-k_0-2} 2^{2(n-k_0-k-2)} \\
&\qquad\qquad \qquad \qquad \qquad\times\#\left\{ 
 \mathbf{x} \in \left(\mathbb{Z}/2^{k+2}\mathbb{Z}\right)^2 \hspace{0.5mm} \bigg| \hspace{0.5mm}  2\nmid \mathbf{x}, \hspace{0.5mm} 
  F_i(\mathbf{x}) \equiv 2^{k_i} \Mod{2^{k_i+2} }
 \right\}\\
\end{aligned}
$$
soit
$$
\sigma_2=2\sum_{k_0 \geqslant 0} \frac{1}{2^{2k_0}} \sum_{0 \leqslant k_1,k_2,k_3} \frac{1}{2^{2k+2}} \#\left\{ 
 \mathbf{x} \in \left(\mathbb{Z}/2^{k+2}\mathbb{Z}\right)^2 \hspace{0.5mm} \bigg| \hspace{0.5mm}  2\nmid \mathbf{x}, \hspace{0.5mm} 
  F_i(\mathbf{x}) \equiv 2^{k_i} \Mod{ 2^{k_i+2} }
 \right\}.
$$
Dans le cas par exemple où $k=k_1$, cela fournit une contribution
$$
2\sum_{k_0\geqslant 0} \frac{1}{2^{2k_0}} \sum_{0 \leqslant k_2,k_3 \leqslant k_1} \frac{1}{2^{2k_1+2}}2^{k_1-\max\{k''_2,k''_3\}}n(\mathbf{k})
$$
soit
$$
2\sum_{k_0\geqslant 0} \frac{1}{2^{2k_0}} \sum_{0 \leqslant k_2,k_3 \leqslant k_1} \frac{n(\mathbf{k})}{2^{k_1+\max\{k''_2,k''_3\}+2}}.
$$
On traite les deux autres cas de la même façon et en regroupant les trois contributions, on obtient finalement la formule
\begin{samepage}
$$
\sigma_2=2\sum_{k_0 \geqslant 0} \frac{1}{2^{2k_0}} \sum_{k_1,k_2,k_3 \geqslant 0} \frac{n(\mathbf{k})}{2^{k_1+\max\{k''_2,k''_3\}+2}}=\frac{2}{3}\sum_{k_1,k_2,k_3 \geqslant 0} \frac{n(\mathbf{k})}{2^{k_1+\max\{k''_2,k''_3\}}}.
$$
\hfill $\square$
\end{samepage}
\\
\newline
\indent
On écrit ensuite
\begin{eqnarray}
S(X)=\sum_{0 \leqslant k_0} \sum_{0 \leqslant \max(k_i) \leqslant \log \log X} S_{\mathbf{k}}\left( 2^{-k_0}X \right)+\sum_{k_0 \geqslant 0} \sum_{\max(k_i)>\log\log X} S_{\mathbf{k}}\left( 2^{-k_0}X \right).
\label{summ}
\end{eqnarray}
D'après ce qui précède, on a
\begin{eqnarray}
S_{\mathbf{k}}\left( X \right)=\sum_{\alpha} S_{\mathbf{k},\alpha}\left( X \right)
\label{skalpha}
\end{eqnarray}
où $\alpha$ parcourt les $n(\mathbf{k})$ solutions de (\ref{conds}) et où
\begin{eqnarray}
S_{\mathbf{k},\alpha}\left( X \right)=\sum_{\mathbf{x}' \in \mathbb{Z}^2 \cap X \mathcal{R}_{\mathbf{M}} \atop x'_1 \equiv 1[4]} r\big( F_{1}\left(\mathbf{M}\mathbf{x}'  \right) \big)r\big( F_{2}\left(\mathbf{M}\mathbf{x}'  \right) \big)r\big( F_3\left(\mathbf{M}\mathbf{x}'  \right) \big),
\label{sommee}
\end{eqnarray}
avec $\mathcal{R}_{\mathbf{M}}=\{ \mathbf{x}' \in \mathbb{R}^2 \quad | \quad \mathbf{Mx}' \in \mathcal{R} \}.$ L'ensemble $\mathbf{M}\mathbb{Z}^2=\{ \mathbf{M}\mathbf{x'}\quad | \quad  \mathbf{x'} \in \mathbb{Z}^2 \}$ est un réseau de covolume
$$
\det(\mathbf{M}\mathbb{Z}^2)=\det(\mathbf{M}).
$$
On considère alors une base réduite $(\mathbf{e}_1,\mathbf{e}_2)$ de $\mathcal{R}_{\mathbf{M}}$, c'est-à-dire une base telle que $\mathbf{e}_1$ soit un vecteur non nul de norme minimale et $\mathbf{e}_2$ soit un vecteur de $\mathcal{R}_{\mathbf{M}} \smallsetminus \mathbf{e}_1\mathbb{Z}$ de norme minimale. Notant $\parallel \mathbf{x}\parallel =\max\{|x_1|,|x_2|\}$, on écrit alors tout $\mathbf{x} \in \mathbb{Z}^2 \cap X \mathcal{R}$ sous la forme
\begin{eqnarray}
\mathbf{x}=\mathbf{Mx}'=\lambda \mathbf{e}_1+\mu \mathbf{e}_2
\label{basemini}
\end{eqnarray}
pour deux entiers $\lambda$ et $\mu$. On pose alors 
\begin{eqnarray}
F'_i(\lambda,\mu)=F_i(\lambda \mathbf{e}_1+\mu \mathbf{e}_2) \quad \mbox{et} \quad F_{\mathbf{M}}=F'_1F'_2F'_3.
\label{formesprime}
\end{eqnarray}
\indent
On montre à présent que le second terme de (\ref{summ}) est un terme d'erreur acceptable dans l'optique du Théorème \ref{theor2}. Pour ce faire, on a besoin des deux lemmes suivants.
\begin{lemme}
Pour tous $\varepsilon>0$ et $X \geqslant 1$, on a $E \ll (\det(\mathbf{M})L_{\infty})^{\varepsilon}$, où 
$$
E=\prod_{4<p\leqslant X}\left( 1+\frac{\rho^{\ast}_{F_{\mathbf{M}}}(p)\chi(p)}{p} \right)\prod_{i=1,2} \prod_{p \leqslant X}\left( 1+\frac{d_i\chi(p)}{p} \right),
$$
avec
$$
\rho^{\ast}_{F_{\mathbf{M}}}(n)=\frac{1}{\varphi(n)}\#\left\{ (x_1,x_2) \in \left(\mathbb{Z}/n\mathbb{Z}\right)^2 \quad | \quad F_{\mathbf{M}}(x_1,x_2) \equiv 0 \Mod{n} ,\quad (x_1,x_2,n)=1 \right\},
$$
$d_i \in \{0,1\}$ et $F_{\mathbf{M}}$ défini en (\ref{formesprime}).
\label{lemme9}
\end{lemme}
\noindent
\textit{Démonstration.}-- On a
$$
\prod_{i=1,2} \prod_{p \leqslant X}\left( 1+\frac{d_i\chi(p)}{p} \right) \ll 1.
$$
Notons $c(f)$ le contenu d'un polynôme $f$ à coefficients entier. Lorsque $p$ divise $\delta:=c(F'_1)c(F'_2)c(F'_3)$, on a une contribution majorée par
$$
\begin{aligned}
\prod_{4<p\atop p \mid \delta}\left( 1+\frac{\rho^{\ast}_{F_{\mathbf{M}}}(p)\chi(p)}{p} \right)&=\prod_{4<p \atop p|\delta} \left( 1+\frac{(p+1)\chi(p)}{p} \right)\\
& \leqslant \prod_{p|\delta} \left( 1+\frac{p+1}{p}\right) \leqslant 3^{\omega(\delta)} \ll \delta^{\varepsilon} \ll (\det(\mathbf{M})L_{\infty})^{\varepsilon}.
\end{aligned}
$$
L'inégalité élémentaire 
$
\rho^{\ast}_{F'_1F'_2F_3'}(p) \leqslant \rho^{\ast}_{F'_1}(p)+\rho^{\ast}_{F'_2}(p)+\rho^{\ast}_{F_3'}(p),
$
valable pour tout nombre premier $p$, implique que les premiers $p$ qui ne divisent pas $\delta$ contribuent pour
$$
\prod_{4<p\leqslant X \atop p\nmid \delta}\left( 1+\frac{\rho^{\ast}_{F_{\mathbf{M}}}(p)\chi(p)}{p} \right)\ll \prod_{4<p\leqslant X \atop p\nmid \delta}\left( 1+\frac{(\rho^{\ast}_{F'_1}(p)+\rho^{\ast}_{F'_2}(p)+\rho^{\ast}_{F_3'}(p))\chi(p)}{p} \right).
$$
On peut alors aisément majorer cette contribution par
$$
\exp\left(\sum_{4<p\leqslant X \atop p\nmid \delta} \frac{(\rho^{\ast}_{F'_1}(p)+\rho^{\ast}_{F'_2}(p)+\rho^{\ast}_{F_3'}(p))\chi(p)}{p}\right).
$$
Or, on a d'une part pour $i \in \{1,2\}$,
$$
\sum_{4<p\leqslant X \atop p\nmid \delta} \frac{\rho^{\ast}_{F'_i}(p)\chi(p)}{p}=\sum_{4<p\leqslant X \atop p\nmid \delta} \frac{\chi(p)}{p}=\sum_{p\nmid \delta} \frac{\chi(p)}{p}+O(1) \leqslant  \log\log(\det(\mathbf{M})L_{\infty})+O(1).
$$
D'autre part,
$$
\sum_{4<p\leqslant X \atop p\nmid \delta}  \frac{\rho^{\ast}_{F_3'}(p)\chi(p)}{p}\ll \sum_{4<p\leqslant X \atop p\nmid \delta}  \frac{\rho_{F_3'(x,1)}(p)\chi(p)}{p}+\log\log(\det(\mathbf{M})L_{\infty})+O(1).
$$
Pour pouvoir conclure, il nous reste donc à montrer que 
$$
\sum_{4<p\leqslant X \atop p\nmid \delta} \frac{\rho_{F_3'(x,1)}(p)\chi(p)}{p} \ll \log\log(\det(\mathbf{M})L_{\infty})+O(1)
$$
sous l'hypothèse que $F_3$ est irréductible sur $\mathbb{Q}[i]$. On voit tout de suite que
$$
\rho_{F_3'(x,1)}(p)=1+\left( \frac{\mbox{disc}(F_3'(x,1))}{p} \right)
$$
où $\left( \frac{.}{p} \right)$ désigne le symbole de Legendre modulo $p$. On décompose $\mbox{disc}(F_3'(x,1))=\varepsilon_{3} u_{3} v_{3}$ où~$\varepsilon_{3} \in \{-1,1\}$, $v_{3}$ est un carré et $u_{3}$ est positif sans facteur carré. On a alors par multiplicativité
$$
\rho_{F_3'(x,1)}(p)=1+\chi(p)^{\frac{1-\varepsilon_{3}}{2}}\left( \frac{u_{3}}{p} \right).
$$
Comme la forme est irréductible sur $\mathbb{Q}[i]$, on a $u_{3} \neq 1$ et on obtient donc un caractère de Dirichlet modulo $u_{3}$ non principal et distinct de $\chi$, ce qui permet de conclure la preuve.
\hfill
$\square$\\
\newline
\begin{lemme}
Soit $\varepsilon>0$. Lorsque $X \geqslant 1$ et $\mathbf{k} \in \mathbb{N}^3$, on a
\begin{eqnarray}
S_{\mathbf{k}}\left( X \right) \ll 2^{\varepsilon(k_1+\max\{k_2,k_3\})}L_{\infty}^{\varepsilon} \left(\frac{X^2}{2^{\max\{k_1,k_2,k_3\}}} +X^{1+\varepsilon} \right),
\label{pff}
\end{eqnarray}
où $S_{\mathbf{k}}$ est définie en (\ref{skalpha}).
\label{leme}
\end{lemme}
\noindent
\textit{Démonstration.}-- On reprend ici les notations (\ref{basemini}). Un cas particulier du résultat de Davenport \cite[lemma 5]{D} permet d'obtenir les estimations
$$
\lambda \ll \frac{\parallel \mathbf{x}\parallel }{\parallel \mathbf{e}_1\parallel } \quad \mbox{et} \quad \mu \ll \frac{\parallel \mathbf{x}\parallel }{\parallel \mathbf{e}_2\parallel }.
$$
D'où, avec les notations (\ref{formesprime}),
$$
S_{\mathbf{k},\alpha}\left(X \right) \ll \sum_{\lambda \ll X/\parallel \mathbf{e}_1\parallel  \atop \mu\ll X/\parallel \mathbf{e}_2\parallel } r\left( F'_{1}\left(\lambda,\mu  \right) \right)r\left( F'_{2}\left(\lambda,\mu   \right) \right)r\left( F_3'\left(\lambda,\mu  \right) \right).
$$
On introduit alors, en notant $r_0=\frac{1}{4}r$, une fonction multiplicative $r_1$ définie de la manière suivante
$$
\forall p \hspace{2mm} \mbox{premier}, \quad \forall \nu \geqslant 1, \quad r_1(p^{\nu})=
\left\{
\begin{array}{ll}
r_0(p)=1+\chi(p) & \mbox{si } \nu=1\\
(1+\nu)^3 &  \mbox{sinon.}\\
\end{array}
\right.
$$ 
Pour tous $k$, $m$ et $n$ entiers, on a ainsi
$$
r_0(k)r_0(m)r_0(n) \leqslant r_1(kmn).
$$
On en déduit la majoration
$$
S_{\mathbf{k},\alpha}\left( X \right) \ll \sum_{\lambda \ll X/\parallel \mathbf{e}_1\parallel  \atop \mu\ll X/\parallel \mathbf{e}_2\parallel } r_1\left(  F_{\mathbf{M}}(\lambda,\mu)   \right)
$$
où $F_{\mathbf{M}}$ est un polynôme de $\mathbb{Z}[x_1,x_2]$ de degré 4 défini en (\ref{formesprime}). On a alors
$$
\parallel F'_i\parallel\ \leqslant (\mbox{deg}(F'_i)+1)\parallel \mathbf{M}\parallel ^{{\rm{deg}}(F'_i)} L_{\infty} 
$$
où $\parallel \mathbf{M}\parallel $ désigne le plus grand coefficient de la matrice en valeur absolue. Ici, à partir de l'expression de la matrice $\mathbf{M}$, on voit immédiatement que $\parallel \mathbf{M}\parallel  \ll 2^{k_1+\max\{k_2,k_3\}}$. On a donc les inégalités
$$
\parallel F_{\mathbf{M}}\parallel  \ll 2^{4(k_1+\max\{k_2,k_3\})}L_{\infty}^3 \ll (2^{k_1+\max\{k_2,k_3\}}L_{\infty})^4.
$$
On déduit finalement du corollaire 1 de \cite{27} et du Lemme \ref{lemme9} la majoration
$$
S_{\mathbf{k},\alpha}\left( X \right) \ll 2^{\varepsilon(k_1+\max\{k_2,k_3\})}L_{\infty}^{\varepsilon} \left(\frac{X^2}{\parallel \mathbf{e}_1\parallel .\parallel \mathbf{e}_2\parallel } +X^{1+\varepsilon} \right)
$$
pour tout $\varepsilon>0$. Or, on sait que pour tout réseau $\Gamma$, de base réduite $(\mathbf{b}_1,\mathbf{b}_2)$, on a
$$
\det(\Gamma) \leqslant \parallel \mathbf{b}_1\parallel .\parallel \mathbf{b}_2\parallel  \ll \det(\Gamma).
$$
On a donc, avec les notations (\ref{basemini}),
$$
\frac{1}{\parallel \mathbf{e}_1\parallel \hspace{0.5mm}.\hspace{0.5mm}\parallel \mathbf{e}_2\parallel } \leqslant \frac{1}{\det(\mathbf{M})}\leqslant\frac{1}{2^{\max\{k_1,k_2,k_3\}}},
$$
ce qui fournit finalement
$$
S_{\mathbf{k},\alpha}\left( X \right) \ll 2^{\varepsilon(k_1+\max\{k_2,k_3\})}L_{\infty}^{\varepsilon} \left(\frac{X^2}{2^{\max\{k_1,k_2,k_3\}}} +X^{1+\varepsilon} \right).
$$
On conclut la preuve grâce au Lemme \ref{lemme6}.
\hfill
$\square$\\
\newline
On déduit de (\ref{skalpha}) et du Lemme \ref{leme} l'estimation
$$
S_{\mathbf{k}}\left( 2^{-k_0}X \right) \ll 2^{\varepsilon(k_1+\max\{k_2,k_3\})}L_{\infty}^{\varepsilon} \left(\frac{X^2}{2^{2k_0+\max\{k_1,k_2,k_3\}}} +\frac{X^{1+\varepsilon}}{2^{(1+\varepsilon)k_0}} \right).
$$
On utilise alors cette majoration pour estimer le deuxième terme de (\ref{summ}). Traitons par exemple le cas de la somme
$$
S_1:=\sum_{0 \leqslant k_0} \sum_{k_2,k_3 \leqslant \log\log(X)} \sum_{k_1>\log\log X} S_{\mathbf{k}}\left( 2^{-k_0}X \right).
$$
Quitte à prendre $\varepsilon$ assez petit, on a par des majorations élémentaires
$$
S_1
\begin{aligned}[t]
&
\ll \sum_{k_2,k_3 \leqslant \log\log(X)} \sum_{k_1>\log\log X} 2^{\varepsilon(k_1+\max\{k_2,k_3\})}L_{\infty}^{\varepsilon} \left(\frac{X^2}{2^{\max\{k_1,k_2,k_3\}}}+X^{1+\varepsilon}\right)\\
& \ll L_{\infty}^{\varepsilon}\log(X)^{\varepsilon}\sum_{k_1>\log\log X} 2^{2\varepsilon k_1} \left(\frac{X^2}{2^{k_1}}+X^{1+\varepsilon}\right)\\
& \ll  L_{\infty}^{\varepsilon}X^2(\log(X))^{3 \varepsilon \log(2)-\log(2)} \ll  L_{\infty}^{\varepsilon}X^2(\log(X))^{3 \varepsilon \log(2)-\eta}.
\end{aligned}
$$
On traite de la même façon le terme d'erreur et les autres cas de sorte qu'avec les notations~(\ref{sommee}),
\begin{eqnarray}
S(X)=\sum_{k_0 \geqslant 0} \sum_{0 \leqslant k_1,k_2,k_3 \leqslant \log\log(X)} S_{\mathbf{k}}\left(2^{-k_0}X \right)+O\left( \frac{L_{\infty}^{\varepsilon}X^2}{(\log(X))^{\eta-\varepsilon}} \right).
\label{sx}
\end{eqnarray}
\newline
\noindent
On pose alors $X'=r'X$, $Y=(r'X)^{1/2}/(\log(X))^C$ pour une constante $C>0$ que l'on fixera plus tard et $F_{i,\mathbf{M}}(\mathbf{x}')=F_i(\mathbf{Mx}')$. Valables lorsque $2^{-\nu_2(m)}m\equiv 1\Mod{4}$ et $0 \leqslant m \leqslant (X')^2$, les décompositions suivantes permettent de restreindre les intervalles dans lesquels varient les variables de façon acceptable:
\begin{enumerate}[i)]
\item
la décomposition
$$
r(m)=4A_+(m)+4A_-(m)
\quad
\mbox{avec}
\quad
A_+(m)=\sum_{d|m \atop d \leqslant \sqrt{X'}} \chi(d) \quad \mbox{et} \quad A_-(m)=\sum_{e|m \atop m> e\sqrt{X'} } \chi(e),
$$
appliquée à $F_{2,\mathbf{M}}(\mathbf{x}')$;
\item
la décomposition
$$
r(m)=4D_+(m)+4D_-(m)
\quad
\mbox{avec}
\quad
D_+(m)=\sum_{d|m \atop d \leqslant X'} \chi(d) \quad \mbox{et} \quad D_-(m)=\sum_{e|m \atop m > eX'} \chi(e),
$$
appliquée à $F_{3,\mathbf{M}}(\mathbf{x}')$;
\item
la décomposition
$$
r(m)=4B_+(m)+4C(m)+4B_-(m),
$$
avec
$$
B_+(m)=\sum_{d|m \atop d \leqslant Y} \chi(d), \quad C(m)=\sum_{d|m \atop Y< d \leqslant X'/Y} \chi(d) \quad \mbox{et} \quad B_-(m)=\sum_{e|m \atop m > eX'/Y} \chi(e),
$$
appliquée à $F_{1,\mathbf{M}}(\mathbf{x}')$.
\end{enumerate}
On notera que, dans $A_-$, $D_-$ et $B_-$, on a respectivement $e \leqslant \sqrt{X'}$, $e \leqslant X'$ et $e \leqslant Y$. On introduit alors les quantités
$$
S_{\pm,\pm,\pm}(X;\mathbf{k},\alpha)=\sum_{\mathbf{x}'\in \mathbb{Z}^2 \cap X\mathcal{R}_{\mathbf{M}} \atop x'_1 \equiv 1[4]}B_{\pm}\left(F_{1,\mathbf{M}}(\mathbf{x}')\right)A_{\pm}\left(F_{2,\mathbf{M}}(\mathbf{x}')\right)D_{\pm}\left(F_{3,\mathbf{M}}(\mathbf{x}')\right)
$$
et
$$
S_0(X;\mathbf{k},\alpha)=\sum_{\mathbf{x}'\in \mathbb{Z}^2 \cap X\mathcal{R}_{\mathbf{M}} \atop x'_1 \equiv 1[4]}C\left(L_{1,\mathbf{M}}(\mathbf{x}')\right)r\left(F_{2,\mathbf{M}}(\mathbf{x}')\right)r\left(F_{3,\mathbf{M}}(\mathbf{x}')\right)
$$
de sorte que
\begin{eqnarray}
S_{\mathbf{k}}\left(X \right)=4^3\sum_{\alpha}\sum_{\pm,\pm,\pm}S_{\pm,\pm,\pm}(X;\mathbf{k},\alpha)+4\sum_{\alpha}S_0(X;\mathbf{k},\alpha).
\label{skx}
\end{eqnarray}

\subsection{Traitement de $S_0$}

 La contribution des sommes $S_0(X;\mathbf{k},\alpha)$ est
$$
S_0(X)=\sum_{k_0 \geqslant 0} \sum_{0 \leqslant k_1,k_2,k_3 \leqslant \log\log(X)} \sum_{\alpha} S_0\left( 2^{-k_0}X;\mathbf{k},\alpha\right).
$$
On montre dans cette section que la somme $S_0$ constitue un terme d'erreur convenable dans l'optique du Théorème \ref{theor2}.
\begin{lemme}
Soient $\varepsilon>0$ et $X \geqslant 1$. Avec la notation (\ref{eta}), on a
$$
S_0(X) \ll \frac{L_{\infty}^{\varepsilon} r_{\infty}r' X^2}{(\log(X))^{\eta-\varepsilon}}.
$$
\label{oups}
\end{lemme}
Le reste de la section est consacrée à la démonstration de ce lemme. On s'inspire ici de la méthode utilisée dans la section 4 de \cite{27}. On a clairement que
$$
S_0\left( X;\mathbf{k},\alpha\right) \ll \sum_{\mathbf{x}'\in \mathbb{Z}^2 \cap X\mathcal{R}_{\mathbf{M}} \atop x'_1 \equiv 1[4]}\left|C\left(F_{1,\mathbf{M}}(\mathbf{x}')\right)\right|r\left(F_{2,\mathbf{M}}(\mathbf{x}')\right)r\left(F_{3,\mathbf{M}}(\mathbf{x}')\right).
$$
On pose alors
$$
E=\left\{m \in \mathbb{Z} \quad | \quad \exists d|m, \quad \mbox{tel que} \quad Y<d \leqslant \frac{X}{Y}\right\},
$$
$$
E_{k_0}=\{ m \in \mathbb{Z} \quad | \quad \exists \mathbf{x} \in 2^{-k_0}X\mathcal{R} \quad \mbox{tel que} \quad F_1(\mathbf{x})=m  \}
$$
et
$$
\mathcal{B}_{k_0}=E \cap E_{k_0}
$$
de sorte que
$$
S_0\left( 2^{-k_0}X;\mathbf{k},\alpha\right) \ll \sum_{m \in \mathcal{B}_{k_0}} S_{0,m}\left( 2^{-k_0}X \right)\left| C(m) \right|
$$
où
$$
S_{0,m}\left( X \right)=\sum_{\mathbf{x} \in \mathbb{Z}^2 \cap X\mathcal{R} \atop F_1(\mathbf{x})=m}r\left(F_{2}(\mathbf{x})\right)r\left(F_3(\mathbf{x})\right). 
$$
On en déduit donc que
$$
S_0(X) \ll (\log\log(X))^3 \sum_{k_0 \geqslant 0}  \sum_{m \in \mathcal{B}_{k_0}} S_{0,m}\left( 2^{-k_0}X \right)\left| C(m) \right|.
$$
D'après \cite[lemma 6]{27}, on a, lorsque $2^{k_0} \leqslant \sqrt{X}$,
$$
\sum_{m \in \mathcal{B}_{k_0}} \left| C(m) \right| \ll \frac{r' 2^{-k_0} X \log \log(X')^{9/4}}{(\log(X'))^{\eta}}.
$$
On en déduit, pour tout $X \geqslant 1$, la majoration
\begin{eqnarray}
S_0(X) \ll \frac{r' X \log \log(X')^{21/4}}{(\log(X'))^{\eta}}\sum_{k_0 \geqslant 0} 2^{-k_0} \max_{m \in \mathbb{N}}\left| S_{0,m}\left( 2^{-k_0}X \right) \right|+X.
\label{lala}
\end{eqnarray}
On montre alors le lemme suivant qui permet de conclure que la contribution de $S_0$ donne un terme d'erreur convenable.
 \begin{lemme}
Il existe une constante absolue $c>0$ telle que
$$
S_{0,m}(X) \ll L_{\infty}^{\varepsilon} r_{\infty} X (\log\log(X'))^{c}.
$$
\label{lemme11}
\end{lemme}
\noindent
\textit{Démonstration.}-- On adapte ici la démonstration du lemme 5 de \cite{27}. Supposons que $a_1 \neq 0$, alors on a
$$
x_1=\frac{m-b_1x_2}{a_1}
$$
et par conséquent
$$
F_2(\mathbf{x})=\frac{A_2m+B_2n}{a_1}=F'_2(m,n)
$$
avec $A_2=a_2$, $B_2=a_1b_2-a_2b_1$ et $n=x_2$. On a alors
$$
F_3(\mathbf{x})=\frac{A_3m^2+B_3 n^2+C_3mn}{a_1^2}=F_3'(m,n)
$$
avec $A_3=a_3$, $B_3=a_3b_1^2+b_3a_1^2-c_3b_1a_1$ et $C_3=c_3a_1-2a_3b_1$. On peut remarquer que $B_2B_3 \neq 0$. On pose alors 
$$
B'_2=\frac{B_2}{\gcd(A_2m,B_2)}, \quad B'_3=\frac{B_3}{\gcd(A_3m^2,B_3,C_3m)}
$$
et
$$
A'_2(m)=\frac{A_2m}{\gcd(A_2m,B_2)}, \hspace{1mm} A'_3(m)=\frac{A_3m^2}{\gcd(A_3m^2,B_3,C_3m)}, \hspace{1mm} C'_3(m)=\frac{C_3m}{\gcd(A_3m^2,B_3,C_3m)}.
$$
On introduit enfin $h:=2 \times 3 \times a_1B_2B_3\Delta^{\ast}$ avec la notation du Lemme \ref{lemme6} et la fonction multiplicative~$r_2$ définie pour $p$ premier et $\nu$ un entier naturel par
$$
r_2\left( p^{\nu} \right)=\left\{
\begin{array}{ll}
 (\nu+1)^2 & \mbox{si } p|h \quad \mbox{ou} \quad \nu \geqslant 2\\
 r_0\left( p^{\nu}\right)=1+\chi(p) & \mbox{sinon}\\
\end{array}
\right.
$$
de sorte que
$$
S_{0,m}(X) \leqslant 16 \tau(h)^2\sum_{n \leqslant R_{\infty}X}r_2(G_m(n)) \ll L_{\infty}^{\varepsilon} \sum_{n \leqslant R_{\infty}X}r_2(G_m(n)),
$$
où 
$$
G_m(n)=\big(A'_2(m)+B'_2n\big)\left(A'_3(m)+B'_3n^2+C'_3(m)n\right).
$$
On rappelle qu'un nombre premier est dit fixe par un polynôme $P$ à coefficients entiers si pour tout entier $n$, on a $P(n) \equiv 0 \Mod{p}$. Le polynôme $G_m$ étant irréductible et primitif de degré 3, seuls $p=2$ et $p=3$ peuvent  être des premiers fixes. Deux applications successives éventuelles du lemme 5 couplées à la remarque 8 de \cite{B} montrent qu'il existe $a \in \{0,1\}$, $b \in \{0,1,2\}$ et 
$
0 \leqslant \mu \leqslant 4
$
tels que
$$
G_{m,2}(x)=\frac{G_m(3^{a}\times 2^{b}x+k)}{3^{a} \times 2^{\mu}}
$$
soit sans premier fixe pour un certain entier $k \leqslant 10$.  On déduit finalement la majoration
$$
S_{0,m}(X) \ll L_{\infty}^{\varepsilon} \sum_{n \ll r_{\infty}X} r_2(G_{m,2}(n)).
$$
Pour $r_{\infty}X \gg  L_{\infty}^{\varepsilon}(X')^{\varepsilon}$ et $\varepsilon \in ]0,1[$, une application du théorème 2 de \cite{B} permet d'obtenir l'estimation
$$
S_{0,m}(X) \ll L_{\infty}^{\varepsilon}r_{\infty} X \prod_{p \ll r_{\infty}X}\left(\left( 1-\frac{\rho_{G_{m,2}}(p)}{p}\right)\sum_{\nu \geqslant 0} \frac{r_2(p^{\nu})\rho_{G_{m,2}}(p^{\nu})}{p^{\nu}} \right).
$$
En utilisant le Lemme \ref{lemme4}, on obtient la majoration
\begin{eqnarray}
\sum_{\nu \geqslant 1} \frac{r_2(p^{\nu})\rho_{G_{m,2}}(p^{\nu})}{p^{\nu}} \leqslant \frac{6}{p}+\frac{27}{p}+3\sum_{\nu \geqslant 3} \frac{(\nu+1)^2}{p^{\frac{\nu}{3}}} \leqslant \frac{6}{p}+\frac{27}{p}+\frac{1167}{p}=\frac{1200}{p}.
\label{kk}
\end{eqnarray}
On déduit de tout cela que les nombres premiers $p$ qui divisent le discriminant de $G_{m,2}$ ont une contribution à $S_{0,m}(X)$
$$
\ll L_{\infty}^{\varepsilon}r_{\infty}X \prod_{p | \text{disc}(G_{m,2})} \left(1+\frac{1200}{p}  \right).
$$
Une étude attentive du passage de $F_m$ à $G_{m,2}$ ainsi que le lemme 1 de \cite{B} montrent que $\mbox{disc}(G_{m,2})\ll L_{\infty}^{19}m^{6}$. On en déduit donc
$$
L_{\infty}^{\varepsilon}r_{\infty}X  \prod_{p | \text{disc}(G_{m,2})} \left(1+\frac{1200}{p}  \right) \ll L_{\infty}^{\varepsilon}r_{\infty}X(\log\log(X'))^{1200},
$$
puisqu'il est nécessaire que $m \leqslant X'$ pour que $S_{0,m}(X)$ soit non nulle. On traite maintenant la contribution des $p$ qui ne divisent pas le discriminant. Pour ces nombres premiers, on a une meilleure majoration de $\rho_{G_{m,2}}(p^{\nu}) \ll 1$, qui donne de la même façon qu'en (\ref{kk}) les majorations
$$
\sum_{\nu \geqslant 2} \frac{r_2(p^{\nu})\rho_{G_{m,2}}(p^{\nu})}{p^{\nu}} \ll \sum_{\nu \geqslant 2} \frac{(\nu+1)^2}{p^{\nu}} \ll \frac{1}{p^2}.
$$
Cela permet de négliger tous les exposants $\nu \geqslant 2$. On en déduit donc finalement
$$
S_{0,m}(X) \ll  L_{\infty}^{\varepsilon}r_{\infty}X(\log\log(X'))^{1200} \prod_{p \ll r_{\infty}X \atop p\nmid \text{disc}(G_{m,2})}\left(1-\frac{\rho_{G_{m,2}}(p)}{p}\right)\left( 1+\frac{r_2(p)\rho_{G_{m,2}}(p)}{p}\right)
$$
soit
$$
\begin{aligned}
S_{0,m}(X)& \ll  L_{\infty}^{\varepsilon}r_{\infty}X(\log\log(X'))^{1200} \underset{p \ll r_{\infty}X \atop p\nmid \text{disc}(G_{m,2})}{\prod}\left( 1+\frac{(r_2(p)-1)\rho_{G_{m,2}}(p)}{p}\right)\\[0mm]
& \ll L_{\infty}^{\varepsilon}r_{\infty}X(\log\log(X'))^{1200} \underset{p \ll r_{\infty}X}{\prod}\left( 1+\frac{\chi(p)\rho_{G_{m,2}}(p)}{p}\right)\\[0mm]
&\ll L_{\infty}^{\varepsilon}r_{\infty}X(\log\log(X'))^{1200+O(1)}.\\
\end{aligned}
$$
On raisonne ici comme lors de la preuve du Lemme \ref{lemme9} et on utilise le fait que pour tout~$p \geqslant 5$, $\rho_{G_{m,2}}(p)=\rho_{G_m}(p)$ et le fait que $G_m$ est irréductible sur $\mathbb{Q}[i]$. 
Dans le cas contraire, lorsque $n \leqslant  L_{\infty}^{\varepsilon}(X')^{\varepsilon}$, on a
$$
r_2(G_{m,2}(n)) \ll G_{m,2}(n)^{\varepsilon} \ll L_{\infty}^{\varepsilon} (X')^{\varepsilon} \ll L_{\infty}^{\varepsilon} {r_{\infty}}^{\varepsilon} X^{\varepsilon}.
$$
La dernière majoration découle de l'inégalité $r' \leqslant 2L_{\infty}r_{\infty}$. D'où,
$$
S_{0,m}(X) \ll  L_{\infty}^{\varepsilon} r_{\infty}X^{\varepsilon} \ll  L_{\infty}^{\varepsilon} r_{\infty}X(\log\log(X'))^{O(1)}
$$
pour $\varepsilon$ assez petit.
\hfill
$\square$\\
\newline
Pour conclure, on tire finalement de la majoration (\ref{lala}) et du Lemme \ref{lemme11} la majoration
$$
S_0(X) \ll \frac{L_{\infty}^{\varepsilon} r_{\infty}r' X^2 \log \log(X')^{21/4+c_0}}{(\log(X'))^{\eta}}\sum_{k_0 \geqslant 0} 2^{-2k_0} \ll \frac{L_{\infty}^{\varepsilon} r_{\infty}r' X^2 \log \log(X')^{21/4+c_0}}{(\log(X'))^{\eta}},
$$
ce qui fournit
$$
S_0(X) \ll \frac{L_{\infty}^{\varepsilon} r_{\infty}r' X^2}{(\log(X'))^{\eta-\varepsilon}}.
$$
Puisque l'hypothèse $r'X^{1-\varepsilon} \geqslant 1$ garantit que
$
\log(X') \gg \log(X),
$
cette estimation est convenable pour obtenir le Lemme \ref{oups}. Ceci achève le traitement de $S_0$.

\subsection{Traitement des $S_{\pm,\pm,\pm}$ et fin de la preuve du Théorème \ref{theor2}}

On traite ici le cas de $S_{-,-,-}(X; \mathbf{k}, \alpha)$ qui est le plus délicat, les autres cas se traitant de manière similaire. On utilise dans la suite les notations $V_1=Y$, $V_2=\sqrt{X'}$, $V_3=X'$ et on pose
$$
\mathcal{R}_{\mathbf{M},4,\mathbf{d}}^{-,-,-,}=\left\{ \mathbf{x}' \in \mathbb{Z}^2 \cap\mathcal{R}_{\mathbf{M}} \quad \Bigg| \quad 
\begin{array}{l}
F_1(\mathbf{Mx'})>r'd_1Y^{-1}\\
F_2(\mathbf{Mx'})>(r')^{1/2}d_2X^{-1/2}\\
F_3(\mathbf{Mx'})> \frac{r'd_3}{X}\\
x'_1 \equiv 1[4]
\end{array}
\right\}.
$$
On a
$$
 S_{-,-,-}(X; \mathbf{k}, \alpha)=\sum_{\mathbf{d} \in \mathbb{N}^3 \atop d_i \leqslant V_i} \chi(d_1d_2d_3) \#\left( \Lambda_{\mathbf{M}}(\mathbf{d}) \cap X\mathcal{R}_{\mathbf{M},4,\mathbf{d}}^{-,-,-,}  \right)
$$
où $\Lambda_{\mathbf{M}}(\mathbf{d})=\Lambda(\mathbf{d},F_{1,\mathbf{M}},F_{2,\mathbf{M}},F_{3,\mathbf{M}})$. On applique alors le lemme de géométrie des nombres suivant qui est dû à Daniel \cite[lemme 3.2]{D99} et La Bretèche et Browning \cite[lemme 5]{30}.
\begin{lemme}
 Soient $\varepsilon>0$, $F_1$, $F_2$ et $F_3$ des formes vérifiant les hypothèses \textbf{NH} (2.0), $X \geqslant 1$, $V_1$, $V_2$, $V_3 \geqslant 2$ et $V=V_1V_2V_3$. Alors il existe une constante absolue $A>0$ telle que
 $$
\sum_{\mathbf{d} \in \mathbb{N}^3 \atop d_i \leqslant V_i}\left| \#\left( \Lambda(\mathbf{d}) \cap X\mathcal{R}_{4,\mathbf{d}} \right)-{\rm{vol}}(\mathcal{R}_{\mathbf{d}})X^2 \frac{\rho(\mathbf{d})}{4(d_1d_2d_3)^2} \right| \ll L_{\infty}^{\varepsilon}(r_{\infty}X\sqrt{V}+V)\log(V)^A,
 $$
\nopagebreak
 avec $\mathcal{R}_{\mathbf{d}} \subset \mathcal{R}$ une région de frontière continûment différentiable dépendant de $\mathbf{d}$ et où $X\mathcal{R}_{4,\mathbf{d}}=\{ \mathbf{x} \in \mathbb{Z}^2 \cap X\mathcal{R}_{\mathbf{d}} \mid x_1 \equiv 1 \Mod{4}\}$.
\label{lemme12}
\end{lemme}
\noindent
\textit{Démonstration.}-- Les éléments de la preuve sont dans l'article \cite{D99} de Daniel ainsi que dans~\cite{M1}, \cite{M2} et \cite{B}. Il faut simplement préciser en plus la dépendance par rapport aux formes et noter que la dépendance de la région par rapport à l'indice de sommation ne change absolument rien dans la preuve. Dernier point, le lemme se démontre pour des formes primitives et on passe aux formes vérifiant \textbf{NH} comme dans \cite{30}.
\hfill
$\square$\\
\newline
\indent
On obtient ici l'égalité
\begin{eqnarray}
 S_{-,-,-}(X; \mathbf{k}, \alpha)=\sum_{\mathbf{d} \in \mathbb{N}^3 \atop d_i \leqslant V_i} \chi(d_1d_2d_3)\frac{X^2}{2^{k_1+\max\{k''_2,k'_3\}+2}}\frac{\rm{vol}(\mathcal{R}_{\mathbf{d}}^{-,-,-})}{4} \frac{\rho(\mathbf{d})}{(d_1d_2d_3)^2}
 \label{smmm}
\end{eqnarray}
 $$
 +O\left( 2^{\varepsilon(k_1+k_2+k_3)}L_{\infty}^{\varepsilon}(r_{\infty}X\sqrt{V}+V)\log(V)^A \right)
$$
où
$$
\mathcal{R}_{\mathbf{d}}^{-,-,-}=\left\{ \mathbf{x} \in \mathcal{R} \quad \bigg| \quad F_1(\mathbf{x})>\frac{r'd_1}{Y}, \quad F_2(\mathbf{x})>\frac{(r')^{1/2}d_2}{X^{1/2}}, \quad F_3(\mathbf{x})> \frac{r'd_3}{X} \right\}.
$$
Posant
$$
T= 2^{\varepsilon(k_1+k_2+k_3)}L_{\infty}^{\varepsilon}(r_{\infty}X\sqrt{V}+V)\log(V)^A,
$$
un calcul élémentaire fournit alors immédiatement
$$
T\ll 2^{\varepsilon(k_1+k_2+k_3)}L_{\infty}^{\varepsilon}X^2\left( \frac{r_{\infty}r'}{(\log(X))^{C/2-A}}+\frac{(r')^2}{(\log(X))^{C-A}} \right).
$$
Si $r' \leqslant r_{\infty}(\log(X))^{A+1}$, l'estimation précédente donne avec $C=2A+8$:
$$
T \ll 2^{\varepsilon(k_1+k_2+k_3)}L_{\infty}^{\varepsilon}X^2\left( \frac{r_{\infty}r'}{(\log(X))^{4}}+\frac{r'r_{\infty}}{(\log(X))^{7}} \right) \ll 2^{\varepsilon(k_1+k_2+k_3)}L_{\infty}^{\varepsilon}X^2\frac{r_{\infty}r'}{(\log(X))^{4}}.
$$
On remarque alors en utilisant le théorème 1 de \cite{B} que
$$
\begin{array}{ll}
T & \ll S_{-,-,-}(X; \mathbf{k}, \alpha)\ll \sum_{x_i \leqslant r_{\infty}X} \tau(F_{1,\mathbf{M}}(\mathbf{x}))\tau(F_{2,\mathbf{M}}(\mathbf{x}))\tau(F_{3,\mathbf{M}}(\mathbf{x}))\\[4mm]
& \ll \sum_{x_i \leqslant r_{\infty}X} \tau_0(F_{1,,\mathbf{M}}(\mathbf{x})F_{2,\mathbf{M}}(\mathbf{x})F_{3,\mathbf{M}}(\mathbf{x})) \ll 2^{\varepsilon(k_1+k_2+k_3)}L_{\infty}^{\varepsilon}r_{\infty}^2X^2(\log(X))^3.\\
\end{array}
$$
où l'on a introduit la fonction multiplicative
$$
\tau_0(p^{\nu})=
\left\{
\begin{array}{ll}
 2=\tau(p^{\nu}) & \mbox{si } \nu=1\\
 (\nu+1)^3 & \mbox{sinon}.
\end{array}
\right.
$$
Cela permet d'obtenir, si $r' \geqslant r_{\infty}(\log(X))^{A+1}$, la majoration
$$
T \ll 2^{\varepsilon(k_1+k_2+k_3)}L_{\infty}^{\varepsilon}r_{\infty}r'X^2(\log(X))^{3-A-1} \ll 2^{\varepsilon(k_1+k_2+k_3)}L_{\infty}^{\varepsilon}r_{\infty}r'X^2(\log(X))^{2-A},
$$
 ce qui implique que dans ce cas également
$$
T \ll 2^{\varepsilon(k_1+k_2+k_3)}L_{\infty}^{\varepsilon}X^2\frac{r_{\infty}r'}{(\log(X))^{4}}.
$$
Pour conclure le traitement du terme d'erreur, il ne reste plus qu'à remplacer $X$ par $2^{-k_0}X$ puis à sommer pour obtenir une contribution
$$
\ll L_{\infty}^{\varepsilon}r_{\infty}r'\sum_{k_0 \geqslant 0} \frac{2^{-2k_0}X^2}{\left(\log\left(2^{-k_0}X\right)\right)^{4}} \sum_{k_1,k_2,k_3  \leqslant \log\log(X)} 2^{\varepsilon(k_1+k_2+k_3)}n(\mathbf{k}).
$$
En utilisant le Lemme \ref{lemme6}, on obtient
$$
\ll L_{\infty}^{\varepsilon}r_{\infty}r'X^2 \frac{(\log(X))^{3\varepsilon\log(2)}}{(\log(X))^4} \ll  L_{\infty}^{\varepsilon}r_{\infty}r' \frac{X^2}{\log(X)}
$$
pour $\varepsilon \leqslant 1/ \log(2)$. Ceci est satisfaisant dans l'optique du Théorème \ref{theor2}.\\
\newline
On passe maintenant à l'étude du terme principal de (\ref{smmm}). On exploite pour ce faire les deux lemmes élémentaires suivants sur les séries de Dirichlet associées à des convolutions pour étudier
\begin{eqnarray}
\mathcal{S}(V_1,V_2,V_3)=\sum_{d_i \leqslant V_i} \chi(d_1d_2d_3) \frac{\rho(d_1,d_2,d_3)}{(d_1d_2d_3)^2}.
\label{S}
\end{eqnarray}
\begin{lemme}
Soient $A>0$, $g$, $w$ deux fonctions arithmétiques et $C$, $C'$, $C''$ trois constantes telles que
$$
\sum_{d=1}^{+\infty} \frac{|w(d)|\log(2d)^A}{d} \leqslant C'' \quad \mbox{et} \quad \sum_{d \leqslant x} \frac{g(d)}{d}=C+O\left( \frac{C'}{(\log(2x))^A} \right).
$$
On a alors que
$$
\sum_{n \leqslant x} \frac{(g \ast w)(n)}{n}=C\sum_{d=1}^{+\infty} \frac{w(d)}{d}+O\left( \frac{C''(C+C')}{(\log(2x))^A} \right).
$$
\label{lemme13}
\end{lemme}
\noindent
\textit{Démonstration}-- La preuve se trouve dans \cite[lemma 5]{36}.
\hfill
$\square$\\
\newline
\begin{lemme}
Soient $g$, $w$ deux fonctions arithmétiques et $C$, $C'$, $C''$ trois constantes telles que
$$
\sum_{d=1}^{+\infty} \frac{|w(d)|\log(2d)}{d} \leqslant C'' \quad \mbox{et} \quad \sum_{d \leqslant x} \frac{g(d)}{d}=C\log(x)+O\left( C' \right).
$$
On a alors que
$$
\sum_{n \leqslant x} \frac{(g \ast w)(n)}{n}=C\log(x)\sum_{d=1}^{+\infty} \frac{w(d)}{d}+O\left( C''(C+C') \right).
$$
\label{lemme14}
\end{lemme}
\noindent
\textit{Démonstration}-- La preuve est très similaire à celle du lemme précédent et nous ne la rédigeons pas ici.
\hfill
$\square$\\
\newline
On établit alors le lemme suivant.
\begin{lemme}
Si on note $V_m=\min\{V_i\}$, pour tout $A>0$, on~a
$$
\mathcal{S}(V_1,V_2,V_3)=\left( \frac{\pi}{4} \right)^3 \prod_{p>2} \sigma_p+O\left(\frac{L_{\infty}^{\varepsilon}}{\log(V_m)^A}  \right).
$$
\label{lemme15}
\end{lemme}
\noindent
\textit{Démonstration}-- On utilise les notations du Lemme \ref{lemme3} pour écrire
$$
\mathcal{S}(V_1,V_2,V_3)=\sum_{d_i \leqslant V_i} \chi(d_1d_2d_3) \frac{(h\ast R)(d_1,d_2,d_3)}{d_1d_2d_3}
$$
où
$$
\chi(d_1d_2d_3)(h\ast R)(d_1,d_2,d_3)=\sum_{e_i|d_i}\chi\left(\frac{d_1}{e_1}\frac{d_2}{e_2}\frac{d_3}{e_3}\right) h\left( \frac{d_1}{e_1},\frac{d_2}{e_2},\frac{d_3}{e_3} \right)\chi(e_1e_2e_3)r_{\Delta}(e_3).
$$
Supposons que $V_3  \geqslant V_2 \geqslant V_1$. On commence alors par sommer sur $d_3$ et $e_3$ si bien qu'on doit estimer la somme suivante
$$
\sum_{d_3 \leqslant V_3} \sum_{e_3|d_3}\frac{\chi\left(\frac{d_3}{e_3}\right) h\left( \frac{d_1}{e_1},\frac{d_2}{e_2},\frac{d_3}{e_3} \right)\chi(e_3)r_{\Delta}(e_3)}{d_3}.
$$
Le Lemme \ref{lemme13}, avec $g=\chi r_{\Delta}=\chi \ast (\chi \chi_{\Delta})$ et $w=\chi h$ où $h$ est vue simplement comme fonction de sa troisième variable fournit alors que cette somme est égale à
$$
\frac{\pi}{4}L(1,\chi\chi_{\Delta})\sum_{k_3 \geqslant 1} \frac{\chi(k_3)h\left( \frac{d_1}{e_1},\frac{d_2}{e_2},k_3 \right)}{k_3}+O\left( \frac{1}{\log(V_3)^A}\sum_{k_3 \geqslant 1}\frac{\left|h\left( \frac{d_1}{e_1},\frac{d_2}{e_2},k_3\right)\right| \log(2k_3)^A }{k_3} \right),
$$
où on a utilisé le fait que $L(1,\chi)=\frac{\pi}{4}.$ On somme alors désormais sur $d_2$ et $e_2$. Le terme principal ci-dessus, devient par une nouvelle application du Lemme \ref{lemme13} avec toujours $w=\chi h$ où $h$ est vue comme fonction de la deuxième variable uniquement cette fois et avec $g=\chi$
\begin{small}
$$
\left(\frac{\pi}{4}\right)^2L(1,\chi\chi_{\Delta})\sum_{k_2,k_3 \geqslant 1} \frac{\chi(k_2k_3)h\left( \frac{d_1}{e_1},k_2,k_3 \right)}{k_2k_3}+O\left( \frac{1}{\log(V_2)^A}\sum_{k_2,k_3  \geqslant 1}\frac{\left|h\left( \frac{d_1}{e_1},k_2,k_3 \right)\right| \log(2k_2)^A}{k_2k_3} \right).
$$
\end{small}
Pour traiter le terme d'erreur de la sommation sur $d_3$ et $e_3$, on va utiliser le Lemme \ref{lemme14} avec~$g=1$ et $w=|h|$ vue comme fonction de sa deuxième variable pour obtenir une contribution
$$
O\left( \frac{\log(V_2)}{\log(V_3)^A}\sum_{k_2,k_3 \geqslant 1}\frac{\left|h\left( \frac{d_1}{e_1},k_2,k_3\right)\right| \log(2k_3)^A\log(2k_2)^A }{k_2k_3} \right).
$$
On effectue alors la même manipulation sur la somme sur $d_1$ et $e_1$ et le Lemme \ref{lemme3} permet d'obtenir que les quantités de la forme
$$
\sum_{k_1,k_2,k_3 \geqslant 1}\frac{\left|h\left( k_1,k_2,k_3\right)\right| \log(2k_3)^A\log(2k_2)^A\log(2k_1)^A }{k_1k_2k_3}
$$
ont une contribution
$$
\ll\sum_{k_1,k_2,k_3 \geqslant 1}\frac{\left|h\left( k_1,k_2,k_3\right)\right| \log(2k_1k_2k_3)^{3A}}{k_1k_2k_3} \ll L_{\infty}^{\varepsilon},
$$
ce qui permet de conclure la preuve du lemme en utilisant l'expression (\ref{convol}) obtenue dans le Lemme \ref{lemme3}.
\hfill
$\square$\\
\newline
On est alors en mesure de terminer la démonstration du Théorème \ref{theor2}. Il ne reste plus qu'à introduire le terme $\rm{vol}(\mathcal{R}_{\mathbf{d}}^{-,-,-})$ dans (\ref{smmm}). Pour ce faire, on pose
\begin{eqnarray}
\mathcal{S}'(X;\mathcal{R})=\sum_{d_i \leqslant V_i} \frac{\chi(d_1d_2d_2)\rho(\mathbf{d}){\rm{vol}}(\mathcal{R}_{\mathbf{d}}^{-,-,-})}{(d_1d_2d_3)^2}.
\label{sprimex}
\end{eqnarray}
On s'inspire alors de \cite[section 7.3]{39} et on écrit
$$
\rm{vol}(\mathcal{R}_{\mathbf{d}}^{-,-,-})=\iint_{\mathbf{x} \in \mathcal{R}} \mathbf{1}_{\mathcal{R}_{\mathbf{d}}^{-,-,-}}(\mathbf{x}) \mbox{d}\mathbf{x}.
$$
En r\'einjectant dans (\ref{sprimex}) et en intervertissant les sommations, on aboutit \`a
$$
\mathcal{S}'(X;\mathcal{R})=\iint_{\mathbf{x} \in \mathcal{R}}  \sum_{\substack{d_1 \leqslant \min\{Y,F_1(\mathbf{x})r'^{-1}Y\},\\ d_2\leqslant \min\{\sqrt{X'},F_2(\mathbf{x})r'^{-1/2}\sqrt{X}\},\\ d_3\leqslant \min\{X',F_3(\mathbf{x})r'^{-1}X\}}} \frac{\chi(d_1d_2d_2)\rho(\mathbf{d})}{(d_1d_2d_3)^2}  \mbox{d}\mathbf{x}.
$$
Le Lemme \ref{lemme15} permet d'obtenir l'estimation
$$
\mathcal{S}'(X;\mathcal{R})=\left(\frac{\pi}{4}\right)^3{\rm{vol}}(\mathcal{R})\prod_{p>2}\sigma_p+
O\left(L_{\infty}^{\varepsilon}I(X) \right),
$$
avec
$$
I(X)=I(X;F_1,F_2,Q)=\iint_{\mathbf{x} \in \mathcal{R}}\frac{1}{\log(\min\{F_1(\mathbf{x})r'^{-1}Y,F_2(\mathbf{x})r'^{-1/2}\sqrt{X},F_3(\mathbf{x})r'^{-1}X\})} \mbox{d}\mathbf{x}.
$$
On établit alors l'estimation 
\begin{eqnarray}
I(X) \ll_{\varepsilon} \frac{r_{\infty}^2L_{\infty}^{\varepsilon}}{\log(X)^{\eta}}.
\label{enfin}
\end{eqnarray}
Lorsque $L_{\infty}> (\log(X))^{1/\varepsilon}$, le Lemme \ref{lemme3} puis le changement de variables $\mathbf{x}=r_{\infty}\mathbf{z}$ permettent d'obtenir les majorations
$$
\mathcal{S}'(X;\mathcal{R}) \ll L_{\infty}^{\varepsilon} \iint_{\mathbf{x} \in \mathcal{R}} \mbox{d}\mathbf{x} \ll r_{\infty}^2L_{\infty}^{\varepsilon}\iint_{\|\mathbf{z}\| \leqslant 1} \mbox{d}\mathbf{z} \ll r_{\infty}^2L_{\infty}^{\varepsilon}.
$$
Ces dernières sont alors suffisantes pour établir (\ref{enfin}). Dans le cas où $L_{\infty}\leqslant (\log(X))^{1/\varepsilon}$, on peut majorer $I(X)$ par la somme 
\begin{small}
$$
\iint_{\mathbf{x} \in \mathcal{R}}\frac{1}{\log(F_3(\mathbf{x})r'^{-1}X+2)} \mbox{d}\mathbf{x} +\iint_{\mathbf{x} \in \mathcal{R}}\frac{1}{\log(F_2(\mathbf{x})r'^{-1/2}X+2)} \mbox{d}\mathbf{x}+\iint_{\mathbf{x} \in \mathcal{R}}\frac{1}{\log(F_1(\mathbf{x})r'^{-1}Y+2)} \mbox{d}\mathbf{x}.
$$
\end{small}
Traitons par exemple le cas du premier terme ci-dessus. En effectuant un changement de variables $\mathbf{x}=r_{\infty}\mathbf{z}$, on a
$$
\iint_{\mathbf{x} \in \mathcal{R}}\frac{1}{\log(F_3(\mathbf{x})r'^{-1}X+2)} \mbox{d}\mathbf{x} \leqslant r_{\infty}^2\iint_{\parallel \mathbf{z}\parallel  \leqslant 1}\frac{1}{\log(F_3(\mathbf{z})r_{\infty}^{2}r'^{-1}X+2)} \mbox{d}\mathbf{z}.
$$
On peut alors constater qu'un changement de variables $\mathbf{u}=E \mathbf{z}$ (voir \cite[section 6]{30}) fait apparaître l'inverse de $\det(E)$ et ainsi, on peut intégrer plutôt sur la forme $J$ définie par $J(\mathbf{z})=F_3(E\mathbf{z})$ (voir \cite{39}) et obtenir par exemple
la majoration
$$
\iint_{\mathbf{x} \in \mathcal{R}}\frac{1}{\log(F_3(\mathbf{x})r'^{-1}X+2)} \mbox{d}\mathbf{x}
\begin{aligned}[t]
 &\leqslant \frac{r_{\infty}^2}{\sqrt{\log(r_{\infty}^{2}r'^{-1}X)}}\iint_{\parallel \mathbf{z}\parallel  \leqslant 1}\frac{1}{\sqrt{|\log(J(\mathbf{z})+2)|}} \mbox{d}\mathbf{z}\\
 &\ll \frac{r_{\infty}^2}{\sqrt{\log(r_{\infty}^{2}r'^{-1}X)}}.
\end{aligned}
$$
Mais les inégalités $r'/(2L_{\infty}) \leqslant r_{\infty} \leqslant 2L_{\infty}r'$ démontrées en section 5 de \cite{30} ainsi que l'hypothèse $r' X^{1-\varepsilon}\geqslant 1$ fournissent la majoration
$$
\frac{1}{\log(r_{\infty}^{2}r'^{-1}X)} \leqslant \frac{1}{\log\left( \frac{X^{\varepsilon}}{4\log(X)^{2/\varepsilon}} \right)},
$$
qui permet de conclure également dans ce cas à la majoration (\ref{enfin}). On traite de manière tout à fait analogue les autres cas, ce qui permet d'obtenir
\begin{eqnarray}
 S_{-,-,-}(X; \mathbf{k}, \alpha)=\frac{X^2}{2^{k_1+\max\{k''_2,k''_3\}+2}} \left(\frac{\pi}{4}\right)^3\frac{\rm{vol}(\mathcal{R})}{4}\prod_{p>2} \sigma_p+O_{\varepsilon}\left( L_{\infty}^{\varepsilon}r_{\infty}^2 \frac{X^2}{(\log(X))^{\eta}}\right).
\label{smmm2}
\end{eqnarray}
En combinant (\ref{sx}), (\ref{skx}), le Lemme \ref{oups}, (\ref{smmm2}) et le Lemme \ref{lemme7}, on conclut à l'estimation
\begin{small}
$$
S(X)=\sum_{k_0 \geqslant 0} \sum_{\mathbf{k} \in \mathbb{N}^3}  \left(\frac{2^{-2k_0+9}X^2 n(\mathbf{k})}{2^{k_1+\max\{k''_2,k''_3\}+2}} \left(\frac{\pi}{4}\right)^3\frac{\rm{vol}(\mathcal{R})}{4}\prod_{p>2} \sigma_p+O_{\varepsilon}\left( L_{\infty}^{\varepsilon}(r_{\infty}r'+r_{\infty}^2) \frac{2^{-2k_0}X^2}{\log(X)^{\eta-\varepsilon}}\right)\right),
\label{th1}
$$
\end{small}
qui permet d'achever la preuve du Théorème \ref{theor2}.

\section{Démonstration du Théorème \ref{theor4}: interprétation de la constante}

Pour $A \in \mathbb{Z}$, $\alpha \in \mathbb{N}$ et $p^n$ une puissance d'un nombre premier, on pose
$$
S_{\alpha}(A;p^n)=\#\left\{ (x,y) \in \left(\mathbb{Z}/p^n\mathbb{Z}\right)^2 \quad | \quad p^{\alpha}(x^2+y^2) \equiv A \Mod{p^n} \right\}.
$$
Si $\alpha \leqslant n$, on a clairement
$$
S_{\alpha}(A;p^n)=p^{2\alpha}S_0(A/p^{\alpha};p^{n-\alpha})
$$
lorsque $\alpha \leqslant \nu_p(A)$ et $S_{\alpha}(A;p^n)=0$ sinon. Il suffit donc de traiter le cas $\alpha=0$, ce que fournit le lemme suivant issu de la section 2 de \cite{27}.
\begin{lemme}
Lorsque $p \equiv 1\Mod{4}$, on a
$$
S_0(A;p^n)=
\left\{
\begin{array}{ll}
 p^n+np^n(1-1/p) & \mbox{si } \nu_p(A) \geqslant n\\
 (1+\nu_p(A))p^n(1-1/p) & \mbox{sinon}.\\
\end{array}
\right.
$$
Lorsque $p \equiv 3\Mod{4}$, on a 
$$
S_0(A;p^n)=
\left\{
\begin{array}{ll}
 p^{2[n/2]} & \mbox{si } \nu_p(A) \geqslant n\\
 p^n(1+1/p) & \mbox{si } \nu_p(A)<n \quad \mbox{et } 2|\nu_p(A)\\
 0 & \mbox{sinon}
\end{array}
\right.
$$
et pour $p=2$, on a
 $$
S_0(A;2^n)=
\left\{
\begin{array}{ll}
 2^{n} & \mbox{si } \nu_2(A) \geqslant n-1\\
 2^{n+1} & \mbox{si } \nu_p(A)<n-1 \quad \mbox{et } 2^{-\nu_2(A)}A \equiv 1 \Mod{4}\\
 0 & \mbox{sinon}.
\end{array}
\right.
$$
\label{lemme16}
\end{lemme}

\subsection{Le cas $p \equiv 1 \Mod{4}$}

On raisonne de manière similaire à la section 4 de \cite{30} à la différence qu'on utilise le Lemme \ref{lemme16}. On pose 
$$
M_{\mathbf{\nu}}(p^n)=\#\left\{\mathbf{x} \in \left(\mathbb{Z}/p^n\mathbb{Z}\right)^2 \quad | \quad \nu_p(F_i(\mathbf{x}))=\nu_i \right\}
$$
et
$$
M'_{\mathbf{\nu}}(p^n)=\#\left\{\mathbf{x} \in \left(\mathbb{Z}/p^n\mathbb{Z}\right)^2 \quad | \quad \nu_p(F_i(\mathbf{x}))\geqslant \nu_i \right\}.
$$
Lorsque $n \geq \nu_1+\nu_2+\nu_3$, on a clairement la formule
\begin{eqnarray}
M'_{\mathbf{\nu}}(p^n)=p^{2n-2\nu_1-2\nu_2-2\nu_3} \rho(p^{\nu_1},p^{\nu_2},p^{\nu_3})
\label{aux}
\end{eqnarray}
et
$$
M_{\mathbf{\nu}}(p^n)=\sum_{\mathbf{e}\in\{0,1\}^3}(-1)^{e_1+e_2+e_3}M'_{\mathbf{\nu}+\mathbf{e}}(p^n).
$$
Avec la notation (\ref{N2}) et posant $m_j=\max\{\lambda_j,\mu_j\}$, on a alors
$$
N_{\mathbf{\lambda},\mathbf{\mu}}(p^n)=p^{3n+\lambda_1+\lambda_2+\lambda_3}\left(1-\frac{1}{p}\right)^3 \sum_{m_j \leqslant \nu_j<n}M_{\mathbf{\nu}}(p^n)\prod_{1\leqslant j\leqslant 3}(1+\nu_j-\lambda_j)+O(n^4p^{4n}).
$$
En utilisant l'identité (\ref{aux}), en divisant par $p^{5n+\lambda_1+\lambda_2+\lambda_3}$ et en faisant tendre $n$ vers l'infini, on obtient
$$
\omega_{\mathbf{\lambda},\mathbf{\mu}}(p)=\left(1-\frac{1}{p}\right)^3 \sum_{m_j \leqslant \nu_j}\sum_{\mathbf{e}\in \{0,1\}^3}(-1)^{e_1+e_2+e_3}\frac{\rho(p^{\nu_1+e_1},p^{\nu_2+e_2},p^{\nu_3+e_3})}{p^{2(\nu_1+e_1+\nu_2+e_2+\nu_3+e_3)}}\prod_{1\leqslant j\leqslant 3}(1+\nu_j-\lambda_j).
$$
On effectue alors les changements de variables $n_j=\nu_j+e_j-\lambda_j$ pour obtenir
$$
\begin{aligned}
\omega_{\mathbf{\lambda},\mathbf{\mu}}(p)=&\left(1-\frac{1}{p}\right)^3 \sum_{m_j-\lambda_j \leqslant n_j}\frac{\rho(p^{n_1+\lambda_1},p^{n_2+\lambda_2},p^{n_3+\lambda_3})}{p^{2(n_1+\lambda_1+n_2+\lambda_2+n_3+\lambda_3)}}\\
&\times \sum_{0 \leqslant e_j \leqslant \min\{1,\lambda_j+n_j-m_j\}}(-1)^{e_1+e_2+e_3}\prod_{1\leq j\leq 3}(1+n_j-e_j).\\
\end{aligned}
$$
Or, on a
$$
\sum_{0\leqslant e \leqslant \min\{1,\lambda+n-m\}}(-1)^e(1+n-e)=
\left\{
\begin{array}{ll}
1 & \mbox{si } \lambda+n-m \geq 1\\
1+m-\lambda & \mbox{si } \lambda+n-m=0\\
\end{array}
\right.
$$
et $(1+m-\lambda)=\#\mathbb{Z} \cap [0,m-\lambda]$. On a donc l'égalité
\begin{small}
$$
\prod_{i=1}^3 \#\mathbb{Z} \cap [0,m_i-\lambda_i] \times \frac{\rho(p^{m_1},p^{m_2},p^{m_3})}{p^{2(m_1+m_2+m_3)}}=\sum_{n_i \leqslant m_i-\lambda_i} \frac{\rho(p^{\max\{m_1,\lambda_1+n_1\}},p^{\max\{m_2,\lambda_2+n_2\}},p^{\max\{m_3,\lambda_3+n_3\}})}{p^{2(\max\{m_1,\lambda_1+n_1\}+\max\{m_2,\lambda_2+n_2\}+\max\{m_3,\lambda_3+n_3\})}}
$$
\end{small}
et finalement
$$
\omega_{\mathbf{\lambda},\mathbf{\mu}}(p)=\left(1-\frac{1}{p}\right)^3 \sum_{n_i \geq 0} \frac{\rho(p^{\max\{m_1,\lambda_1+n_1\}},p^{\max\{m_2,\lambda_2+n_2\}},p^{\max\{m_3,\lambda_3+n_3\}})}{p^{2(\max\{m_1,\lambda_1+n_1\}+\max\{m_2,\lambda_2+n_2\}+\max\{m_3,\lambda_3+n_3\})}}.
$$
On a donc bien obtenu l'identité $\omega_{\boldsymbol{\lambda},\boldsymbol{\mu}}(p)=\sigma_p(\mathbf{d},\mathbf{D})$ pour tout nombre premier $p \equiv 1 \Mod{4}$.

\subsection{Le cas $p \equiv 3 \Mod{4}$}

Avec les mêmes notations que lors de la section précédente, lorsque $p \equiv 3 \Mod{4}$, on aboutit de la même manière à l'expression
\begin{small}
$$
N_{\mathbf{\lambda},\mathbf{\mu}}(p^n)=p^{5n+\lambda_1+\lambda_2+\lambda_3}\left(1+\frac{1}{p}\right)^3 \sum_{m_i \leqslant \nu_i \atop 2|\nu_i-\lambda_i} \sum_{\mathbf{e} \in \{0,1\}^3}(-1)^{e_1+e_2+e_3}\frac{\rho(p^{\nu_1+e_1},p^{\nu_2+e_2},p^{\nu_3+e_3})}{p^{2(\nu_1+e_1+\nu_2+e_2+\nu_3+e_3)}}+O(n^4p^{4n}). 
$$
\end{small}
Passant à la limite et effectuant le même changement de variables que dans le cas précédent, on obtient
$$
\omega_{\mathbf{\lambda},\mathbf{\mu}}(p)=\left(1+\frac{1}{p}\right)^3 \sum_{n_i \geqslant m_i-\lambda_i} \frac{\rho(p^{n_1+\lambda_1},p^{n_2+\lambda_2},p^{n_3+\lambda_3})}{p^{2(n_1+\lambda_1+n_2+\lambda_2+n_3+\lambda_3)}} \sum_{0 \leqslant e_i \leqslant \min\{1,\lambda_i+n_i-m_i\} \atop e_i \equiv n_i [2]} (-1)^{e_1+e_2+e_3}.
$$
Or,
$$
\sum_{0 \leqslant e \leqslant \min\{1,\lambda+n-m\}\atop e \equiv n[2]} (-1)^e=
\left\{
\begin{array}{ll}
(-1)^n & \mbox{si } \lambda+n-m \geq 1\\
1 & \mbox{si } \lambda+n-m=0 \quad \mbox{et } 2|m-\lambda\\
0 & \mbox{si } \lambda+n-m=0 \quad \mbox{et } 2\not|m-\lambda.\\
\end{array}
\right.
$$
On conclut alors comme lors de la section précédente à l'égalité $\omega_{\boldsymbol{\lambda},\boldsymbol{\mu}}(p)=\sigma_p(\mathbf{d},\mathbf{D})$ pour tout nombre premier $p \equiv 3 \Mod{4}$.

\subsection{Le cas $p=2$}

Le Lemme \ref{lemme16} fournit immédiatement la relation
$$
\omega_{\mathbf{d}}(2)=\lim_{n \rightarrow +\infty} 2^{-2n+3}\#\left\{ 
 \mathbf{x} \in \left(\mathbb{Z}/2^n\mathbb{Z}\right)^2 \quad \bigg| \quad 
  F_i(\mathbf{x}) \in d_i\mathcal{E}_{2^n}
 \right\}
=\sigma_2(\mathbf{d}),
$$
où $\sigma_2(\mathbf{d})$ est définie en (\ref{sigmadeux}) et $\omega_{\mathbf{d}}(2)$ en (\ref{omegadeux}).

\subsection{Le cas de la densité archimédienne}

Enfin, pour terminer le traitement de la constante, il reste à regarder la densité archimédienne. On remarque pour commencer que
$$
\omega_{\mathcal{R}}(\infty)=2^6\omega_{\mathcal{R}}^+(\infty)
$$
où $\omega_{\mathcal{R}}^+(\infty)$ est défini de la même façon que $\omega_{\mathcal{R}}(\infty)$ avec les conditions supplémentaires $s_i>0$ et $t_i>0$ pour tout $1 \leqslant i \leqslant 3$. On utilise la forme de Leray en paramétrant par les $t_i$. La forme de Leray est par conséquent ici donnée par
$$
(-2^3t_1t_2t_3)^{-1}\mbox{d}s_1\mbox{d}s_2\mbox{d}s_3\mbox{d}x_1\mbox{d}x_2.
$$
En effet, la variété est définie comme le lieu des zéros des polynômes définis par
$$
f_i(\mathbf{x},\mathbf{s},\mathbf{t})=\frac{F_i(\mathbf{x})}{d_i}-(s_i^2+t_i^2)
$$
et on a 
$$
\det\begin{pmatrix}
 \frac{\partial f_1}{\partial t_1}& \frac{\partial f_2}{\partial t_1} & \frac{\partial f_3}{\partial t_1}\\[2mm]
 \frac{\partial f_1}{\partial t_2} & \frac{\partial f_2}{\partial t_2} & \frac{\partial f_3}{\partial t_2}\\[2mm]
 \frac{\partial f_1}{\partial t_3} & \frac{\partial f_2}{\partial t_3} & \frac{\partial f_3}{\partial t_3}\\[2mm]
\end{pmatrix}=
\det\begin{pmatrix}
 -2t_1& 0 & 0\\
 0 & -2t_2 & 0\\
 0 & 0 & -2t_3
 \end{pmatrix}=-2^3t_1t_2t_3.
$$
On utilise ici le calcul de l'intégrale suivante
\begin{eqnarray}
\int_0^{\sqrt{A}} \frac{\mbox{d}s}{\sqrt{A-s^2}}=\frac{\pi}{2}.
\label{pi2}
\end{eqnarray}
En substituant $t_i=\sqrt{d_i^{-1}F_i(\mathbf{x})-s_i^2}$, on obtient finalement
$$
\omega_{\mathcal{R}}^+(\infty)=2^{-3} \int_{\mathbf{x} \in \mathcal{R}} \left( \prod_{1 \leqslant i \leqslant 3} \int_0^{\sqrt{d_i^{-1}F_i(\mathbf{x})}} \frac{\mbox{d}s_i}{\sqrt{d_i^{-1}F_i(\mathbf{x})-s_i^2}} \right)\mbox{d}x_1\mbox{d}x_2
$$
soit
$$
\omega_{\mathcal{R}}(\infty)=\pi^3 \rm{vol}(\mathcal{R}).
$$
Cela achève la preuve du Théorème \ref{theor4}.

\section{Démonstration du Théorème \ref{theor1}}
\subsection{Passage aux torseurs et reformulation du problème de comptage}

On note dans la suite $Z^m$ l'ensemble des vecteurs de $\mathbb{Z}^m$ premiers entre eux dans leur ensemble. On introduit également la norme de $\mathbb{R}^5$ suivante
\begin{eqnarray}
\parallel \mathbf{x}\parallel =\max\left\{|x_0|,|x_1|,|x_2|,\delta^{-1}|x_3|,\delta^{-1}|x_4|\right\},
\label{norm}
\end{eqnarray}
où 
\begin{eqnarray}
\delta=\sqrt{(|a_1|+|b_1|)(|a_2|+|b_2|)(|a_3|+|b_3|+|c_3|)}.
\label{delt}
\end{eqnarray}
On introduit la hauteur exponentielle sur $\mathbb{P}^4(\mathbb{Q})$, associée à la norme définie en (\ref{norm}), définie par
$$
H_4:
\left\{
\begin{array}{lll}
 \mathbb{P}^4(\mathbb{Q}) & \longrightarrow & \mathbb{R}^+_{\ast}\\
  {[}x_0:x_1:x_2:x_3:x_4] & \longmapsto & \parallel ( x_{0} ,x_1,x_2,x_3,x_4)\parallel 
\end{array}
\right.
$$
avec $[x_0:x_1:x_2:x_3:x_4]$ choisi de telle sorte que les $x_i$ soient des entiers premiers entre eux. On introduit aussi $\mathcal{T}_{{\rm spl}} \subset \mathbb{A}^5_{\mathbb{Q}}=\mbox{Spec}\left(\mathbb{Q}[y,z,t,u,v]\right)$ le sous-schéma définie par l'équation
\begin{eqnarray}
y^2+z^2=t^2F_1(u,v)F_2(u,v)F_3(u,v)
\label{T}
\end{eqnarray}
avec les conditions $(y,z,t) \neq 0$ et $(u,v) \neq 0$. Il s'agit d'un $\mathbb{G}_m^2$-torseur pour $S$ \cite[Définition~4.1]{35}. On pose 
$$
\mathcal{D}=\{ d \in \mathbb{N} \quad | \quad p|d \Rightarrow p\equiv 1\Mod{4} \}.
$$
Lorsque $d_0 \in \mathcal{D}$, alors tous les diviseurs de $d_0$ sont aussi dans l'ensemble $\mathcal{D}$. On introduit ensuite
$$
r(n;m)=\#\{ (a,b) \in \mathbb{Z}^2 \quad | \quad n=a^2+b^2, \quad (m,a,b)=1 \}.
$$
On a alors clairement $r(n;1)=r(n)$ et on remarque que $r(y^2n;y)=0$ si $y \not \in \mathcal{D}$.\\
Utilisant une inversion de Möbius pour traiter la condition de coprimalité, on aboutit, pour~$y \in \mathcal{D}$, à la formule suivante
\begin{eqnarray}
r(y^2n;y)=\sum_{k|y} \mu(k)r\left( \frac{y^2n}{k^2} \right).
\label{mobius}
\end{eqnarray}
Enfin, considérant l'ensemble
\begin{eqnarray}
\Sigma:=\left\{ \boldsymbol{\varepsilon} \in \{-1,+1\}^3  \quad | \quad \varepsilon_1\varepsilon_2 \varepsilon_3=1, \quad \varepsilon_1=1\right\}
\label{sigma}
\end{eqnarray}
et pour $\boldsymbol{\varepsilon} \in \Sigma$ et $T \geqslant 1$ la région
$$
R^{\boldsymbol{\varepsilon}}(T)=\left\{ (u,v) \in \mathbb{R}^2 \quad \Bigg| \quad \begin{array}{l}
         |u|,|v| \leqslant \sqrt{T},\\                                                                                           
         \varepsilon_i F_i(u,v)>0                                                                                          \end{array}
     \right\},
$$
on en déduit le lemme suivant.
\begin{lemme}
On a $N(B)=N_1(B)+O(B)$ où
$$
N_1(B)=\sum_{k \in \mathcal{D}}\mu(k) \sum_{\ell \leqslant B/k \atop \ell \in \mathcal{D}} \sum_{\boldsymbol{\varepsilon} \in \Sigma} \sum_{(u,v) \in Z^2 \cap R^{\boldsymbol{\varepsilon}}(B/k\ell)}r\left(\ell^2F_1^+(u,v)F_2^+(u,v)F_3^+(u,v) \right),
$$
où l'on note $F_i^+=\varepsilon_i F_i$.
\label{lemme18}
\end{lemme}
\noindent
\textit{Démonstration}-- Le lemme 2 de \cite{Brow} nous garantit que $N(B)=\frac{1}{4}T(B)$ avec
$$
T(B)=\# \left\{(y,z,t;u,v) \in Z^3 \times Z^2 \quad \Bigg| \quad \begin{array}{l}
                                                                  \parallel  (v^2t,uvt,u^2t,y,z)) \parallel  \leqslant B,\\
                                                                  y^2+z^2=t^2F(u,v)
                                                                 \end{array}
    \right\}.
$$
Puisque pour un $(y,z,t;u,v)$ compté, on a 
$$
\parallel  (v^2t,uvt,u^2t,y,z)) \parallel =\max\{u^2,v^2\}|t|,
$$
on en déduit, par symétrie sur le signe de $t$, que
\begin{eqnarray}
N(B)=\frac{1}{2}\# \left\{(y,z,t;u,v) \in (Z^3 \times Z^2) \cap \mathcal{T}_{{\rm spl}} \quad | \quad  0<\max\{u^2,v^2\}t \leqslant B  \right\}.
\label{N}
\label{formule}
\end{eqnarray}
Un raisonnement élémentaire montre que la contribution des $(u,v)$ tels que $$F_1(u,v)F_2(u,v)F_3(u,v)=0$$ est $O(1)$ ce qui permet d'écrire $N(B)=N_1(B)+O(B)$ avec
$$
N_1(B)=\frac{1}{2}\sum_{\substack{t \leqslant B \\ t \in \mathcal{D}}} \sum_{\substack{\varepsilon_i \in \{-1,+1\} \\ \varepsilon_1 \varepsilon_2 \varepsilon_3=1}} \sum_{(u,v) \in Z^2 \cap R^{\boldsymbol{\varepsilon}}(B/t)}r(t^2F_1(u,v)F_2(u,v)F_3(u,v);t).
$$
En posant $k\ell=t$, la formule (\ref{mobius}) fournit par conséquent
$$
N_1(B)=\frac{1}{2}\sum_{k \in \mathcal{D}}\mu(k) \sum_{\ell \leqslant B/k \atop \ell \in \mathcal{D}} \sum_{\varepsilon_i \in \{-1,+1\} \atop \varepsilon_1 \varepsilon_2 \varepsilon_3=1} \sum_{(u,v) \in Z^2 \cap R^{\boldsymbol{\varepsilon}}(B/k\ell)}r\left(\ell^2F_1^+(u,v)F_2^+(u,v)F_3^+(u,v) \right).
$$
On en déduit alors immédiatement la formule donnée dans l'énoncé du lemme avec $\varepsilon_1=1$. Le passage à la somme sur $\boldsymbol{\varepsilon} \in \Sigma$ sera utile en section 7 puisque l'ensemble $\Sigma$ sera utilisé pour décrire certaines classes d'isomorphie de torseurs.
\hfill
$\square$\\
\newline
Dans la suite, on notera pour éviter d'alourdir les notations $\omega(\gcd(a_0,\ldots,a_n))=\omega(a_0,\ldots,a_n)$. L'idée est maintenant de passer d'un problème de comptage sur $\mathcal{T}_{{\rm spl}}$ à un problème de comptage sur des variétés affines de $\mathbb{A}^8$ de la forme (\ref{var}), autrement dit d'exprimer le problème de comptage sur certains torseurs qui seront explicités en section 7. On introduit pour ce faire, lorsque $\mathbf{d}=(d_1,d_2,d_3)$, $d=d_1d_2d_3$,  $\mathbf{n}=(n_1,n_2,n_3)$ et $\mathbf{d}'=(d'_1,d'_2,d'_3)$, les notations
$$
c(d, \mathbf{n})=
3^{\omega(d,n_1,n_2,n_3)}\times
\prod_{\{i,j,k\}=\{1,2,3\}\atop i<j}2^{\omega(d,n_i,n_j)-\omega(d,n_i,n_j,n_k)}
$$
et
$$
c'(\mathbf{d}, \mathbf{d}', \mathbf{n})=
c(d, \mathbf{n})\times 2^{\omega(d'_1d'_2,n_1/d_1,n_2/d_2)} .
$$
On a alors le lemme suivant où on rappelle qu'on note $r_0=\frac{r}{4}$.
\begin{lemme}
Soient $n_0$, $n_1$, $n_2$ et $n_3$ quatre entiers supérieurs ou égaux à 1 avec $n_0 \in \mathcal{D}$. On a alors la formule suivante
$$
r_0(n_0 n_1 n_2 n_3)=
\begin{aligned}[t]
&
\sum_{\substack{\mathbf{d}, \mathbf{d}' \in \mathcal{D}^3 \\ d_1d_2d_3 |n_0, d_id'_jd'_k | n_i}}\sum_{m_{i}|C(n_{i})} \frac{\mu(d_1d_2d_3)\mu(d'_1d'_2)\mu(d'_3) }{c'(\mathbf{d}, \mathbf{d}', \mathbf{n})}r_0\left( \frac{n_0}{d_1d_2d_3} \right)\\
& \times  r_0\left( \frac{n_1}{d_1d'_2d'_3m_{1}} \right)r_0\left( \frac{n_2}{d_2d'_1d'_3m_{2}} \right)r_0\left( \frac{n_3}{d_3d'_1d'_2m_{3}} \right),
\end{aligned}
$$
où $\{i,j,k\}$ parcourt l'ensemble des permutations de $\{1,2,3\}$ et où
\begin{eqnarray}
C(n_i)=\prod_{\substack{p\equiv 3\Mod{4}\\ \nu_p(n_i)\equiv 1 \Mod{2}\\ \nu_p(n_1n_2n_3)\equiv 0 \Mod{2}}} p.
\label{nni}
\end{eqnarray}
\label{lemme19}
\end{lemme}
\noindent
\textit{Démonstration}-- Utilisant la formule d'éclatement pour deux entiers (voir par exemple \cite[lemme 10]{30}), on obtient
$$
r_0(n_0 n_1 n_2 n_3)=\sum_{d|(n_1 n_2 n_3,n_0)} \mu(d) r_0\left( \frac{n_0}{d} \right)r_0\left( \frac{n_1 n_2 n_3}{d} \right).
$$
Dans la somme, les entiers $d$ sont sans facteur carré et dans $\mathcal{D}$ lorsque $n_0\in \mathcal{D}$. On compte alors le nombre de décompositions $d=d_1 d_2 d_3$ avec $d_1|n_1$, $d_2|n_2$ et $d_3|n_3$.  Posant
$$
N(d,\mathbf{n})=\#\left\{ (d_1,d_2,d_3) \in \mathbb{N}^3 \quad | \quad d_i | n_i, \quad \mbox{et} \quad d=d_1d_2d_3 \right\}
$$
le nombre de telles décompositions, 
on a l'égalité $c(d, \mathbf{n})=N(d,\mathbf{n})$, lorsque $d$ est sans facteur carré. En effet, $c(\hspace{1mm}\cdot\hspace{1mm},\mathbf{n})$ et $N(\hspace{1mm}\cdot\hspace{1mm},\mathbf{n})$ sont deux fonctions multiplicatives. Il suffit donc de montrer qu'elles coïncident sur les nombres premiers. Soit alors $p$ un nombre premier. Dans le cas où $p | n_1n_2n_3$, on a clairement $$N(p,\mathbf{n})=\#\{ i \hspace{2mm} | \hspace{2mm} p | n_i \}=c(p,\mathbf{n}).$$
Cela permet bien de conclure à l'égalité souhaitée. On a par conséquent
$$
r_0(n_0 n_1 n_2 n_3)=
\sum_{d_1d_2d_3|n_0 \atop d_1 | n_1, d_2 |n_2, d_3|n_3} \frac{\mu(d_1d_2d_3 )}{c(d, \mathbf{n})}
r_0\left( \frac{n_0}{d_1d_2d_3} \right)r_0\left( \frac{n_1}{d_1}\frac{n_2}{d_2}\frac{n_3}{d_3} \right).
$$
Posant
$$
n_i^{(1)}=\prod_{p\equiv 1\Mod{4}}p^{\nu_p(n_i)}, \quad n_i^{(3)}=\prod_{p\equiv 3 \Mod{4}}p^{\nu_p(n_i)},
$$
on a par multiplicativité,
\begin{eqnarray}
r_0\left( \frac{n_1}{d_1}\frac{n_2}{d_2}\frac{n_3}{d_3} \right)=r_0\left( \frac{n_1^{(1)}}{d_1}\frac{n_2^{(1)}}{d_2}\frac{n_3^{(1)}}{d_3}\right)r_0\left( n_1^{(3)}n_2^{(3)}n_3^{(3)}\right).
\label{new}
\end{eqnarray}
On utilise alors la formule d'éclatement pour trois entiers donnée par le lemme 9 de \cite{35} pour le premier terme du membre de droite de (\ref{new}):
\begin{small}
$$
r_0\left( \frac{n_1^{(1)}}{d_1}\frac{n_2^{(1)}}{d_2}\frac{n_3^{(1)}}{d_3}\right)=\sum_{d'_id'_j|n_k/d_k\atop d'_i \in \mathcal{D}} \frac{\mu(d'_1)\mu(d'_2d'_3)}{2^{\omega(d'_2,n_2/d_2)+\omega(d'_3,n_3/d_3)}} r_0\left(\frac{n_1^{(1)}}{d_1d'_2d'_3} \right)r_0\left( \frac{n_2^{(1)}}{d_2d'_1d'_3}\right)r_0\left( \frac{n_3^{(1)}}{d_3d'_1d'_2}\right).
$$
\end{small}
Concernant le deuxième terme du membre de droite de (\ref{new}), en remarquant que $r_0$ ne prend que les valeurs 0 ou 1 sur des entiers n'ayant que des facteurs premiers congrus à 3 modulo 4, on a les égalités
$$
r_0\left( n_1^{(3)}n_2^{(3)}n_3^{(3)}\right)
\begin{aligned}[t]
&=r_0\left( \frac{n_1^{(3)}}{C(n_{1})}\right)r_0\left(\frac{n_2^{(3)}}{C(n_{2})}\right)r_0\left(\frac{n_3^{(3)}}{C(n_{3})}\right)\\
&=\sum_{m_{i}|C(n_{i})}r_0\left( \frac{n_1^{(3)}}{m_{1}}\right)r_0\left(\frac{n_2^{(3)}}{m_{2}}\right)r_0\left(\frac{n_3^{(3)}}{m_{3}}\right).
\end{aligned}
$$
La seconde égalité découle du fait que pour tout $m_{i}$ divisant strictement $C(n_{i})$, on a
$$
r_0\left( \frac{n_1^{(3)}}{m_{1}}\right)r_0\left(\frac{n_2^{(3)}}{m_{2}}\right)r_0\left(\frac{n_3^{(3)}}{m_{3}}\right)=0.
$$
En regroupant les facteurs par multiplicativité, on obtient bien la formule souhaitée.
\hfill
$\square$\\
\newline
On note dans la suite
$$
F_{i,e}=eF_i^+=e\varepsilon_i F_i \quad \mbox{et} \quad F_{3,e}=e^2F_3^+=e^2\varepsilon_3F_3,
$$
pour tout entier naturel $e$. On introduit également
\begin{eqnarray}
\Delta_{ij}^{(1)}=\prod_{p | \Delta_{ij} \atop p \equiv 1 \Mod{4}}p \quad \mbox{et} \quad \Delta_{ij}^{(3)}=\prod_{p | \Delta_{ij} \atop p \equiv 3 \Mod{4}}p
\label{deelta}
\end{eqnarray}
et l'ensemble
\begin{eqnarray}
M=\left\{ \mathbf{m}\in \mathbb{N}^3 \quad \Bigg| \quad  
\begin{array}{l}
m_{i} \bigg| \left[\Delta_{ij}^{(3)},\Delta_{ik}^{(3)}\right], \quad (m_{i},m_{j}) \big| \Delta_{ij}^{(3)}\\[0.3cm]
 \mu^2(m_{i})=1, \quad  \sqrt{m_{1}m_{2}m_{3}} \in \mathbb{N}
\end{array}, \quad \{i,j,k\}=\{1,2,3\}
\right\}.
\label{m}
\end{eqnarray}
Posant $\mathcal{R}=R^{\boldsymbol{\varepsilon}}(1)$ et pour $d \in \mathbb{N}$ fixé,
\begin{eqnarray}
f_d(n)=\sum_{sk=n}\mu(k)r(ds^2),
\label{fd}
\end{eqnarray}
on déduit alors des Lemmes \ref{lemme18} et \ref{lemme19} le lemme suivant.
\begin{lemme}
Avec les notations (\ref{Ss}) et (\ref{sigma}), on a
$$
N_1(B)=
\begin{aligned}[t]
&
\frac{1}{2^6}\underset{e\geqslant 1 \atop d \in \mathcal{D}}{\sum}\mu(e)\mu(d) \underset{\substack{n \leqslant N \\ n \in \mathcal{D}}}{\sum}f_d(n) \underset{\boldsymbol{\varepsilon} \in \Sigma}{\sum}\hspace{1mm}\underset{\substack{k_4k_1k'_1|\gcd(\Delta_{23},d) \\ k_4k_2k'_2|\gcd(\Delta_{13},d)\\ k_4k_3k'_3|\gcd(\Delta_{12},d)}}{\sum} \sum_{\substack{ k_4k'_4|\gcd(\Delta_{12},\Delta_{13},\Delta_{23},d)\\ k_5k'_5|\Delta_{12}}} 
\\
& 
\frac{\mu(k'_1)\mu(k'_2)\mu(k'_3)\mu(k'_4)\mu(k'_5)}{3^{\omega(k_4)}2^{\omega(k_5)+\omega(k_1)+\omega(k_2)+\omega(k_3)}}
 \sum_{\substack{\mathbf{d}, \mathbf{d}' \in \mathcal{D}^3 \\ d=d_1d_2d_3, d'_i |\Delta_{jk} \\ k_5k'_5 | d'_1d'_2}}\sum_{\mathbf{m} \in M}\mu(d'_1d'_2)\mu(d'_3) S\left(\sqrt{\frac{B}{de^2n}},\mathbf{e},\mathbf{E}\right),
\end{aligned}
$$
%$$
%N_1(B)=\frac{1}{2^6}\underset{e\in \mathbb{N}_{>0} \atop d \in \mathcal{D}}{\sum}\mu(e)\mu(d) \underset{\substack{n \leqslant N \\ n \in \mathcal{D}}}{\sum}f_d(n) \underset{\boldsymbol{\varepsilon} \in \Sigma}{\sum}\underset{\substack{k_4k_1k'_1|\gcd(\Delta_{23},d) \\ k_4k_2k'_2|\gcd(\Delta_{13},d)\\ k_4k_3k'_3|\gcd(\Delta_{12},d)}}{\sum} 
%\hspace{6cm}
%$$
%\begin{footnotesize}
%$$
%\times
%\sum_{\substack{ k_4k'_4|\gcd(\Delta_{12},\Delta_{13},\Delta_{23},d)\\ k_5k'_5|\gcd(\Delta_{12},d'_1d'_2)}} \frac{\mu(k'_1)\mu(k'_2)\mu(k'_3)\mu(k'_4)\mu(k'_5)}{3^{\omega(k_4)}2^{\omega(k_5)+\omega(k_1)+\omega(k_2)+\omega(k_3)}} \sum_{\substack{\mathbf{d}, \mathbf{d}' \in \mathbb{N}^3 \\ d=d_1d_2d_3, d'_j| \Delta_{ik}}}(\chi\mu)(d'_1d'_2)(\chi\mu)(d'_3) S\left(\sqrt{\frac{B}{de^2n}},\mathbf{e},\mathbf{E}\right),
%$$
%\end{footnotesize}
\noindent
où $N=\frac{\delta B}{d^{5/4} e}$ et
$$
S(e^{-1}\sqrt{T},\mathbf{e},\mathbf{E})=S(e^{-1}\sqrt{T},\mathbf{e},\mathbf{E};R^{\boldsymbol{\varepsilon}}, F_{1,e},F_{2,e},F_{3,e})
$$
pour $T \geqslant 1$ et avec
\begin{eqnarray}
\begin{split}
&e_1=d_1 d'_2 d'_3 m_{1}, \quad e_2=d_2d'_1d'_3m_{2}, \quad e_3=d_3d'_1d'_2m_{3},\\
&E_1=[e_1,k_4k_2k'_2,k_4k_3k'_3,k_4k'_4,d_1k_5k'_5],\\
 &E_2=[e_{2},k_4k_1k'_1,k_4k_3k'_3,k_4k'_4,d_2k_5k'_5],\\
&E_3=[e_{3},k_4k_1k'_1,k_4k_2k'_2,k_4k'_4].
\label{E}
\end{split}
\end{eqnarray}
\label{lemme20}
\end{lemme}
\noindent
\textit{Démonstration.}-- Pour $(u,v)$ premiers entre eux, on a 
\begin{eqnarray}
\gcd(F_i(u,v),F_j(u,v)) | \Delta_{ij}.
\label{dprimes}
\end{eqnarray}
Or, grâce au Lemme \ref{lemme19}, on a l'égalité
\begin{small}
$$
\begin{aligned}[t]
& r\left(\ell^2F_1^+(u,v)F_2^+(u,v)F_3^+(u,v) \right)=\frac{1}{2^6}\sum_{\substack{\mathbf{d}, \mathbf{d}' \in \mathcal{D}^3 \\ d_1d_2d_3 |\ell^2, d_id'_jd'_k | F_i(u,v)}}\sum_{\mathbf{m}\in M \atop m_{i} | F_i(u,v)} \\
& \frac{\mu(d_1d_2d_3)\mu(d'_1d'_2)\mu(d'_3)}{c'(\mathbf{d}, \mathbf{d}',F_1^+(u,v),F_2^+(u,v),F_3^+(u,v))}r\left( \frac{\ell^2}{d_1d_2d_3} \right)r\left( \frac{F_1^+(u,v)}{e_1} \right)r\left( \frac{F_2^+(u,v)}{e_{2}} \right)r\left( \frac{F_3^+(u,v)}{e_{3}} \right),
\end{aligned}
$$
\end{small}
En effet, par (\ref{dprimes}), un triplet $\mathbf{m}$ tel que $m_i | C(F_i(u,v))$ avec la notation (\ref{nni}), est dans~$M$. Par ailleurs, pour un triplet $\mathbf{m} \in M$ tel que $m_i | F_i(u,v)$ pour tout $1 \leqslant i \leqslant 3$ mais $m_j \nmid C(F_j(u,v))$ pour un certain $j$, alors
$$
r\left( \frac{F_1^+(u,v)}{d_1d'_2d'_3m_{1}} \right)r\left( \frac{F_2^+(u,v)}{d_2d'_1d'_3m_{2}} \right)r\left( \frac{F_3^+(u,v)}{d_3d'_1d'_2m_{3}} \right)=0.
$$
L'entier $d=d_1d_2d_3$ divise $\ell$ et est donc nécessairement dans $\mathcal{D}$. On écrit alors $\ell=ds$ et on obtient grâce au Lemme \ref{lemme18}
$$
N_1(B)=\frac{1}{2^6} \sum_{\substack{dk \leqslant B \\ k,d \in \mathcal{D}}} \mu(d)\mu(k) \sum_{\substack{s \leqslant \frac{B}{dk} \\ s \in \mathcal{D}}} r(ds^2)\sum_{\substack{\mathbf{d}, \mathbf{d}' \in \mathcal{D}^3 \\ d=d_1d_2d_3}}\sum_{\mathbf{m}\in M}\mu(d'_1d'_2)\mu(d'_3) \mathcal{S}_{\mathbf{d}, \mathbf{d}',\mathbf{m}}\left( \frac{B}{dsk} \right)
$$
où, pour $T \geqslant 1$, on a posé
$$
 \mathcal{S}_{\mathbf{d}, \mathbf{d}',\mathbf{m}}(T)=\sum_{\boldsymbol{\varepsilon} \in \Sigma} \sum_{\substack{(u,v) \in Z^2 \cap R^{\boldsymbol{\varepsilon}}(T) \\ e_{1} | F_1(u,v), e_{2} | F_2(u,v),\\ e_{3} | F_3(u,v)}} \frac{r\left( \frac{F_1^+(u,v)}{e_{1}} \right)r\left( \frac{F_2^+(u,v)}{e_{2}} \right)r\left( \frac{F_3^+(u,v)}{e_{3}} \right)}{c'(\mathbf{d}, \mathbf{d}',F_1^+(u,v),F_2^+(u,v),F_3^+(u,v))}.
$$
Avec la notation (\ref{fd}) et grâce à (\ref{dprimes}), on aboutit à la formule
$$
N_1(B)=\frac{1}{2^6} \sum_{dn \leqslant B \atop dn \in \mathcal{D}} \mu(d)f_d(n) \sum_{\mathbf{d}, \mathbf{d}' \in \mathcal{D}^3 \atop d=d_1d_2d_3, d'_j|\Delta_{ik}}\sum_{\mathbf{m}\in M}\mu(d'_1d'_2)\mu(d'_3) \mathcal{S}_{\mathbf{d}, \mathbf{d}',\mathbf{m}}\left( \frac{B}{dn} \right).
$$
On obtient alors en utilisant à nouveau (\ref{dprimes})
$$
\mathcal{S}_{\mathbf{d}, \mathbf{d}',\mathbf{m}}(T)=
\begin{aligned}[t]
& 
\sum_{\boldsymbol{\varepsilon} \in \Sigma}\sum_{\substack{k_1k_4|\gcd(\Delta_{23},d) \\ k_2k_4|\gcd(\Delta_{13},d) \\k_3k_4|\gcd(\Delta_{12},d) }}\sum_{\substack{ k_4|\gcd(\Delta_{12},\Delta_{13},\Delta_{23},d)\\ k_5|\gcd(\Delta_{12},d'_1d'_2)}} \frac{1}{3^{\omega(k_4)}2^{\omega(k_5)+\omega(k_1)+\omega(k_2)+\omega(k_3)}}\\
& \times\sum_{\substack{(u,v) \in Z^2 \cap R^{\boldsymbol{\varepsilon}}(T) \\ e_{1} | F_1(u,v), e_{2} | F_2(u,v) \\ e_{3} | F_3(u,v)}}
r\left( \frac{F_1^+(u,v)}{e_{1}} \right)r\left( \frac{F_2^+(u,v)}{e_{2}} \right)r\left( \frac{F_3^+(u,v)}{e_{3}} \right),
\end{aligned}
$$
où la somme intérieure porte sur les couples $(u,v)$ tels que
\begin{eqnarray}
\left\{
\begin{array}{l}
k_1=\gcd(d,F_2(u,v),F_3(u,v))/\gcd(d,F_1(u,v),F_2(u,v),F_3(u,v)),\\[2mm] k_2=\gcd(d,F_1(u,v),F_3(u,v))/\gcd(d,F_1(u,v),F_2(u,v),F_3(u,v)),\\[2mm] k_3=\gcd(d,F_1(u,v),F_2(u,v))/\gcd(d,F_1(u,v),F_2(u,v),F_3(u,v)),\\[2mm] k_4=\gcd(d,F_1(u,v),F_2(u,v),F_3(u,v)),\\[2mm]
k_5=\gcd(d'_1d'_2,F_1(u,v)/d_1,F_2(u,v)/d_2).\\
\end{array}
\right.
\label{conditions}
\end{eqnarray}
Plusieurs inversions de Möbius fournissent alors
$$
\mathcal{S}_{\mathbf{d}, \mathbf{d}',\mathbf{m}}(T)=
\begin{aligned}[t]
&
\sum_{\boldsymbol{\varepsilon} \in \Sigma}\sum_{\substack{k_4k_1k'_1|\gcd(\Delta_{23},d) \\ k_4k_2k'_2|\gcd(\Delta_{13},d)\\ k_4k_3k'_3|\gcd(\Delta_{12},d)}}\sum_{\substack{k_4k'_4|\gcd(\Delta_{12},\Delta_{13},\Delta_{23},d)\\ k_5k'_5|\gcd(\Delta_{12},d'_1d'_2) }} \frac{\mu(k'_1)\mu(k'_2)\mu(k'_3)\mu(k'_4)\mu(k'_5)}{3^{\omega(k_4)}2^{\omega(k_5)+\omega(k_1)+\omega(k_2)+\omega(k_3)}} \\
&  \times\sum_{\substack{(u,v) \in Z^2 \cap R^{\boldsymbol{\varepsilon}}(T) \\ e_{1} | F_1(u,v), e_{2} | F_2(u,v) \\ e_{3} | F_3(u,v)}}
r\left( \frac{F_1^+(u,v)}{e_{1}} \right)r\left( \frac{F_2^+(u,v)}{e_{2}} \right)r\left( \frac{F_3^+(u,v)}{e_{3}} \right),
\end{aligned}
$$
où la somme intérieure porte désormais sur les couples $(u,v)$ tels que
\begin{eqnarray}
\left\{
\begin{array}{l}
k_4k_1k'_1|\gcd(d,F_2(u,v),F_3(u,v)),\\[2mm] 
k_4k_2k'_2|\gcd(d,F_1(u,v),F_3(u,v)),\\[2mm] 
k_4k_3k'_3|\gcd(d,F_1(u,v),F_2(u,v)),\\[2mm]
k_4k'_4|\gcd(d,F_1(u,v),F_2(u,v),F_3(u,v)),\\[2mm]
k_5k'_5|\gcd(d'_1d'_2,F_1(u,v)/d_1,F_2(u,v)/d_2).\\
\end{array}
\right.
\label{conditions2}
\end{eqnarray}
On a $(\mathbf{e},\mathbf{E}) \in \mathfrak{D}$ et on peut alors réécrire cette somme sous la forme
\begin{eqnarray}
\begin{split}
\mathcal{S}_{\mathbf{d}, \mathbf{d}',\mathbf{m}}(T)=
& \sum_{\boldsymbol{\varepsilon} \in \Sigma}\sum_{\substack{k_4k_1k'_1|\gcd(\Delta_{23},d) \\ k_4k_2k'_2|\gcd(\Delta_{13},d)\\ k_4k_3k'_3|\gcd(\Delta_{12},d)}}\sum_{\substack{k_4k'_4|\gcd(\Delta_{12},\Delta_{13},\Delta_{23},d)\\ k_5k'_5|\gcd(\Delta_{12},d'_1d'_2)  }} \frac{\mu(k'_1)\mu(k'_2)\mu(k'_3)\mu(k'_4)\mu(k'_5)}{3^{\omega(k_4)}2^{\omega(k_5)+\omega(k_1)+\omega(k_2)+\omega(k_3)}}\\
& \times
\sum_{(u,v) \in Z^2 \cap \sqrt{T}R^{\boldsymbol{\varepsilon}}(1) \cap \Lambda(\mathbf{E})}
r\left( \frac{F_1^+(u,v)}{e_1} \right)r\left( \frac{F_2^+(u,v)}{e_2} \right)r\left( \frac{F_3^+(u,v)}{e_3} \right).
\label{sommme}
\end{split}
\end{eqnarray}
Reste encore à enlever la condition de coprimalité sur les couples $(u,v)$ au moyen d'une dernière inversion de Möbius.  Avec la notation (\ref{Ss}), la somme intérieure de (\ref{sommme}) est alors égale à
$$
\sum_{e =1}^{+\infty} \mu(e)S(e^{-1}\sqrt{T},\mathbf{e},\mathbf{E}).
$$
On obtient ainsi le lemme avec $N=\frac{B}{de^2}$. Mais on peut remarquer que $S\left(\sqrt{\frac{B}{de^2n}},\mathbf{e},\mathbf{E}\right)$ est nulle si l'on n'a pas
$$
d_1 \leqslant (|a_1|+|b_1|)\sqrt{\frac{B}{dn}}, \quad d_2 \leqslant (|a_2|+|b_2|)\sqrt{\frac{B}{dn}} \quad \mbox{et} \quad d_3 \leqslant (|a_3|+|b_3|+|c_3|)\frac{B}{dn}.
$$
En effet, on a par exemple
$$
F_{1,e}(u,v) \leqslant (|a_1|+|b_1|)e \max\{|u|,|v|\} \leqslant (|a_1|+|b_1|)\sqrt{\frac{B}{dn}}
$$
dans la région considérée. Avec la notation (\ref{delt}), on a donc en réalité que
\begin{eqnarray}
d^{\frac{3}{4}}n^{\frac{1}{2}} \leqslant \delta B^{\frac{1}{2}}.
\label{inegal}
\end{eqnarray}
Mais pour que $S\left(\sqrt{\frac{B}{de^2n}},\mathbf{e},\mathbf{E}\right)$ soit non nulle, il faut aussi imposer
$
de^2n \leqslant B
$
et par conséquent,
$
d^{\frac{1}{2}}en^{\frac{1}{2}} \leqslant  B^{\frac{1}{2}},
$
ce qui, combiné avec l'inégalité (\ref{inegal}), implique l'inégalité
$$
d^{\frac{5}{4}}en \leqslant \delta B
$$
et démontre bien le lemme énoncé.
\hfill
$\square$

\subsection{Fin de la preuve du Théorème \ref{theor1}}

On obtient alors le théorème suivant.
\begin{theor}
Lorsque $B \rightarrow +\infty$, on a
$$
N(B) =c_0 B \log(B)(1+o(1)),
$$
avec
$$
c_0=
\begin{aligned}[t]
& \frac{\pi^2}{2^4}\sum_{e=1}^{+\infty} \frac{\mu(e)}{e^2} \sum_{d \in \mathcal{D}} \frac{\mu(d)r_0(d)\varphi^{\dagger}(d)}{d}\sum_{\boldsymbol{\varepsilon} \in \Sigma} {\rm{vol}}\left( R^{\boldsymbol{\varepsilon}}(1) \right)\sum_{\substack{\mathbf{d}, \mathbf{d}' \in \mathcal{D}^3 \\ d=d_1d_2d_3, d'_i| \Delta_{jk}}}\sum_{\mathbf{m}\in M}\mu(d'_1d'_2)\mu(d'_3)\\
&\times
\sum_{\substack{k_4k_1k'_1|\gcd(\Delta_{23},d) \\ k_4k_2k'_2|\gcd(\Delta_{13},d) \\k_4k_3k'_3|\gcd(\Delta_{12},d) }}\sum_{\substack{ k_4k'_4|\gcd(\Delta_{12},\Delta_{13},\Delta_{23},d)\\ k_5k'_5|\gcd(\Delta_{12},d'_1d'_2)}} \frac{\mu(k'_1)\mu(k'_2)\mu(k'_3)\mu(k'_4)\mu(k'_5)}{3^{\omega(k_4)}2^{\omega(k_5)+\omega(k_1)+\omega(k_2)+\omega(k_3)}} \sigma^{\boldsymbol{\varepsilon}}(\mathbf{e},\mathbf{E},e),
\end{aligned}
$$
avec
\begin{eqnarray}
\sigma^{\boldsymbol{\varepsilon}}(\mathbf{e},\mathbf{E},e):=\prod_{p} \sigma_p^{\boldsymbol{\varepsilon}}(\mathbf{e},\mathbf{E},e) ,
\label{jesaisplus}
\end{eqnarray}
où
$$
 \sigma_p^{\boldsymbol{\varepsilon}}(\mathbf{e},\mathbf{E},e)=\left( 1-\frac{\chi(p)}{p} \right)^3 \sum_{\boldsymbol{\nu} \in \mathbb{N}^3} \frac{\chi(p)^{\nu_1+\nu_2+\nu_3}\rho\left( p^{N_1},p^{N_2},p^{N_3} \right)}{p^{2(N_1+N_2+N_3)}},
$$
$N_i=\max\{\nu_p(E_i),\nu_i+\nu_p(e_i)\}$
et
$$
\sigma_2^{\boldsymbol{\varepsilon}}(\mathbf{e},e)=8\lim_{n \rightarrow +\infty} 2^{-2n} \#\left\{ 
 \mathbf{x} \in \left(\mathbb{Z}/2^n\mathbb{Z}\right)^2 \quad \bigg| \quad 
  e^{{\rm{deg}}(F_i)} F_i(\mathbf{x}) \in  \varepsilon_i m_i\mathcal{E}_{2^n}
 \right\}.
$$
\label{theor5}
\end{theor}
La fin de cette section est consacrée à la preuve de ce théorème. Il permet de démontrer le Théorème~\ref{theor1} modulo le fait que la constante $c_0$ soit bien la constante conjecturée par Peyre, ce qui fera l'objet de la section suivante. On utilise évidemment le Théorème \ref{theor3} pour estimer la somme 
$$
S(\sqrt{T},\mathbf{e},\mathbf{E})=S(\sqrt{T},\mathbf{e},\mathbf{E};\mathcal{R}, F_{1,e},F_{2,e},F_{3,e})
$$
et obtenir le lemme suivant.
\begin{lemme}
Soit $\varepsilon>0$ tel que $r'(\sqrt{T})^{1-\varepsilon} \geqslant 1$. On a alors avec les notations du Théorème \ref{theor5}
$$
S(\sqrt{T},\mathbf{e},\mathbf{E})=\pi^3 {\rm{vol}}(R^{\boldsymbol{\varepsilon}}(1)) T \prod_{p} \sigma_p^{\boldsymbol{\varepsilon}}(\mathbf{e},\mathbf{E},e)+
O_{\varepsilon,F_{i,e}}\left( \frac{T}{\left(\log(T)\right)^{\eta-\varepsilon}}\right).
$$
De plus, on a 
$$
\sigma^{\boldsymbol{\varepsilon}}(\mathbf{e},\mathbf{E},e) \ll (de)^{\varepsilon}a'(\mathbf{E},\mathbf{\Delta}),
$$
uniformément en tous les paramètres $\mathbf{d}$, $\mathbf{d}'$, $\mathbf{m}$, $\mathbf{k}$, $\mathbf{k}'$, $\varepsilon_i$, $e$, $\mathbf{e}$, $\mathbf{E}$ et $d$.
\label{lemme21}
\end{lemme}
On raisonne alors comme dans la preuve du lemme 9 de \cite{30} afin d'obtenir une majoration uniforme de $S(\sqrt{T},\mathbf{e},\mathbf{E})$ mis à part qu'on utilise une majoration différente pour la quantité~$\det(G_f(\mathcal{A}))$ qui y apparaît. On remarque ici qu'il faut bien faire attention au fait qu'on travaille avec les formes $F_{i,e}$ et qu'on a notamment une dépendance des coefficients en la variable $e$. On a ainsi le lemme suivant.
\begin{lemme}
Pour $T\geqslant 1$ et $(\mathbf{e},\mathbf{E}) \in \mathfrak{D}$, on a
$$
\begin{aligned}
S(\sqrt{T},\mathbf{e},\mathbf{E})\ll& c(F_1)c(F_2)c(F_3)(\Delta^{\ast})^3(de L_{\infty})^{\varepsilon}\gcd(d,e)a'(\mathbf{E},\mathbf{\Delta})\\
&\times
\left( r_{\infty}(R^{\boldsymbol{\varepsilon}}(1))^2\frac{T}{d}+r_{\infty}(R^{\boldsymbol{\varepsilon}}(1))^{1+\varepsilon}T^{\frac{1}{2}+\varepsilon}\right),
\end{aligned}
$$
où l'on note $c(P)$ le contenu d'un polynôme à coefficients entiers.
\label{lemme22}
\end{lemme}
\noindent
\textit{Démonstration.}-- On raisonne ici comme dans la section 6 de \cite{30}. Avec les notations (\ref{not}), on pose
$$
E''_i=\frac{E_i}{\gcd(E'_i,b)\gcd(E_i,e\ell_i)} \quad \mbox{et} \quad E''_3=\frac{E_3}{\gcd(E'_3,b^2)\gcd(E_3,e^2\ell_3)}.
$$
Des calculs élémentaires garantissent qu'on a
$$
E''=E''_1E''_2E''_3 \geqslant \frac{E}{\gcd(E_1,be\ell_1)\gcd(E_2,be\ell_2)\gcd(E_3,b^2e^2\ell_3)}.\hspace{1.9cm}
$$
Mais revenant à la définition de $\frac{E_1}{\gcd(E_1,be\ell_1)}$, on obtient avec la notation du Lemme \ref{lemme6}
$$
\frac{E_1}{\gcd(E_1,be\ell_1)} \geqslant \frac{1}{\Delta^{\ast}}\frac{d_1}{\gcd(d_1,be\ell_1)}.
$$
En procédant de la même façon avec les indices 2 et 3, on aboutit finalement à l'inégalité
$$
E'' \geqslant \frac{1}{(\Delta^{\ast})^3}\frac{d}{\gcd(d_1,be\ell_1)\gcd(d_2,be\ell_2)\gcd(d_3,b^2e^2\ell_3)}.
$$
On obtient alors
$$
E'' \geqslant \frac{1}{c(F_1)c(F_2)c(F_3)(\Delta^{\ast})^3}\frac{d}{\gcd(d_1,be)\gcd(d_2,be)\gcd(d_3,b^2e^2)}.
$$
Puisque $d=d_1d_2d_3$ est sans facteur carré, $d_3$ l'est aussi et donc
$$
E'' \geqslant \frac{1}{c(F_1)c(F_2)c(F_3)(\Delta^{\ast})^3}\frac{d}{\gcd(d_1,be)\gcd(d_2,be)\gcd(d_3,be)}.
$$
Pour finir, on utilise le fait que
$$
\gcd(d_1,be)\gcd(d_2,be)\gcd(d_3,be)=\gcd(d_1d_2d_3,be)=\gcd(d,be).
$$
pour conclure à l'inégalité
$$
E'' \geqslant \frac{1}{c(F_1)c(F_2)c(F_3)(\Delta^{\ast})^3}\frac{d}{\gcd(d,be)}.
$$
Cela permet de finir la démonstration du lemme exactement comme dans le lemme 9 de \cite{30}.
\hfill
$\square$\\
\newline
On a ensuite besoin de quelques résultats issus de \cite[lemme 14]{36} concernant la fonction~$f_d$ définie en (\ref{fd}) pour $d \in \mathcal{D}$.
Si $d \in \mathcal{D}$ est sans facteur carré, alors, pour $x \geqslant 2$, on a
\begin{eqnarray}
\sum_{n \leqslant x \atop n \in \mathcal{D}} \frac{f_d(n)}{n}=\frac{4r_0(d)\varphi^{\dagger}(d)}{\pi}\left( \log(x)+O\left(  \log^3(2+\omega(d)) \right) \right),
\label{abcd}
\end{eqnarray}
où
$$
\varphi^{\dagger}(d)=\prod_{p|d}\left( 1+\frac{1}{p} \right)^{-1}.
$$
Pour $\varepsilon>0$ et $0<\theta \leqslant 1$, on déduit de ces résultats la majoration
\begin{eqnarray}
\sum_{n \leqslant x \atop n \in \mathcal{D}} \frac{|f_d(n)|}{n^{\theta}} \leqslant x^{1-\theta}\sum_{n \leqslant x \atop n \in \mathcal{D}} \frac{|f_d(n)|}{n} \ll d^{\varepsilon}x^{1-\theta}\log(x).
\label{abcde}
\end{eqnarray}
\newline
\indent
On en déduit à présent la conjecture de Manin pour les surfaces de Châtelet considérées. On utilise la formule asymptotique obtenue grâce au Théorème \ref{theor3} de $S(\sqrt{T},\mathbf{e},\mathbf{E})$ dans l'expression de $N_1(B)$ que l'on a obtenue dans le Lemme \ref{lemme21}. Pour pallier la non uniformité de cette estimation, on pose
$$
\mathcal{U}(B)=\mathcal{U}_{\mathbf{e},\mathbf{E},e}^{\boldsymbol{\varepsilon}}(B)=\sum_{n \leqslant N \atop n \in \mathcal{D}} f_d(n) S\left(\sqrt{\frac{B}{de^2n}},\mathbf{e},\mathbf{E}\right),
$$
avec $N$ défini lors du Lemme \ref{lemme20}, de telle sorte que
$$
N_1(B)=
\begin{aligned}[t]
& \frac{1}{2^6}\sum_{e \geqslant 1 \atop d \in \mathcal{D}}\mu(e) \mu(d)  \sum_{\boldsymbol{\varepsilon} \in \Sigma}\sum_{\substack{\mathbf{d}, \mathbf{d}' \in \mathcal{D}^3 \\ d=d_1d_2d_3, d'_i| \Delta_{jk}}}\sum_{\mathbf{m}\in M}\mu(d'_1d'_2)\mu(d'_3)\\
& \times
\sum_{\substack{k_4k_1k'_1|\gcd(\Delta_{23},d) \\ k_4k_2k'_2|\gcd(\Delta_{13},d) \\k_4k_3k'_3|\gcd(\Delta_{12},d) }}\sum_{\substack{ k_4k'_4|\gcd(\Delta_{12},\Delta_{13},\Delta_{23},d) \\ k_5k'_5|\gcd(\Delta_{12},d'_1d'_2)}} \frac{\mu(k'_1)\mu(k'_2)\mu(k'_3)\mu(k'_4)\mu(k'_5)}{3^{\omega(k_4)}2^{\omega(k_5)+\omega(k_1)+\omega(k_2)+\omega(k_3)}}\mathcal{U}(B).
\end{aligned}
$$
Avec la notation (\ref{jesaisplus}), on introduit alors la quantité
$$
E^{\boldsymbol{\varepsilon}}(B;\mathbf{e},\mathbf{E},e)=
\frac{1}{B\log(B)}\left|\mathcal{U}(B)- 4\pi^2 {\rm{vol}}(R^{\boldsymbol{\varepsilon}}(1))\sigma^{\boldsymbol{\varepsilon}}(\mathbf{e},\mathbf{E},e)r_0(d)\varphi^{\dagger}(d)\frac{B \log(B)}{de^2}  \right|.
$$
Les résultats (\ref{abcd}) et (\ref{abcde}) impliquent que
$$
E^{\boldsymbol{\varepsilon}}(B;\mathbf{e},\mathbf{E},e) \underset{B \rightarrow +\infty}{\longrightarrow} 0,
$$
à $\mathbf{e}$, $\mathbf{E}$, $\mathbf{d}$, $\mathbf{d'}$, $\mathbf{m}$, $\mathbf{k}$, $\mathbf{k}'$, $d$ et $e$ fixés. Pour conclure, on souhaite appliquer un théorème de convergence dominée et il suffit donc de montrer la majoration
$$
\sum_{e \geqslant 1 \atop d \in \mathcal{D}}\sum_{\boldsymbol{\varepsilon} \in \Sigma}\sum_{\substack{\mathbf{d}, \mathbf{d}' \in \mathcal{D}^3 \\ d=d_1d_2d_3, d'_i| \Delta_{jk}}}\sum_{\mathbf{m}\in M}\sum_{\substack{k_4k_1k'_1|\gcd(\Delta_{23},d) \\ k_4k_2k'_2|\gcd(\Delta_{13},d)\\ k_4k_3k'_3|\gcd(\Delta_{12},d) }}\sum_{\substack{ k_4k'_4|\gcd(\Delta_{12},\Delta_{13},\Delta_{23},d) \\k_5k'_5|\gcd(\Delta_{12},d'_1d'_2) }}E^{\boldsymbol{\varepsilon}}(B;\mathbf{e},\mathbf{E},e) \ll 1.
$$
Pour ce faire, on a besoin du lemme suivant.
\begin{lemme}
Pour tout $\varepsilon>0$, on a
$$
 \left|\sigma^{\boldsymbol{\varepsilon}}(\mathbf{e},\mathbf{E},e)\right| \ll d^{-\frac{1}{6}+\varepsilon}e^{\varepsilon}.
$$
\label{lemme24}
\end{lemme}
\noindent
\textit{Démonstration.}-- On procède en majorant chaque facteur eulérien dans la définition de~$\sigma^{\boldsymbol{\varepsilon}}(\mathbf{e},\mathbf{E},e)$ en (\ref{jesaisplus}). Pour $p=2$, on peut majorer $\sigma_2^{\boldsymbol{\varepsilon}}(\mathbf{e},\mathbf{E},e)$ par 8. Pour les $p$ impairs, on utilise la majoration donnée par le point f) du Lemme \ref{lemme2}. On a alors
$$
\left|\sigma_p^{\boldsymbol{\varepsilon}}(\mathbf{e},\mathbf{E},e)\right|\ll \sum_{\boldsymbol{\nu}\in \mathbb{N}^3} (N_3+1)\frac{1}{p^{\frac{N_1}{3}+\frac{N_2}{3}+\frac{N_3}{6}}},
$$
où les $N_i=\max\{\nu_p(E_i),\nu_i+\nu_p(e_i)\}$ ont été définis dans le Lemme \ref{lemme21}. On tire alors parti du fait que $d$ soit sans facteur carré. On note $d=d_1d_2d_3$ et $\delta_i=\nu_p(d_i)$ de sorte qu'on peut supposer $\delta_1+\delta_2+\delta_3=1$ lorsque $p|d$. De plus, par définition de $E_i$ et $N_i$, on voit que~$N_i \geqslant \nu_i+\delta_i$ et donc on a
$$
\prod_{p|d} \left|\sigma_p^{\boldsymbol{\varepsilon}}(\mathbf{e},\mathbf{E},e)\right| \ll \prod_{p|d} p^{-\frac{\delta_1+\delta_2+\delta_3}{6}} \prod_{p|d}\sum_{\boldsymbol{\nu}\in \mathbb{N}^3} (N_3+1)\frac{1}{p^{\frac{\nu_1}{3}+\frac{\nu_2}{3}+\frac{\nu_3}{6}}}.
$$
On a alors
$$
\prod_{p|d} p^{-\frac{\delta_1+\delta_2+\delta_3}{6}}=\prod_{p|d} p^{-\frac{1}{6}}=d^{-\frac{1}{6}}
$$
puisque $d$ est sans facteur carré. D'autre part, 
$$
\prod_{p|d} \sum_{\boldsymbol{\nu}\in \mathbb{N}^3} (N_3+1)\frac{1}{p^{\frac{\nu_1}{3}+\frac{\nu_2}{3}+\frac{\nu_3}{6}}} \ll d^{\varepsilon} \prod_{p|d} (1+\nu_p(d))\ll d^{\varepsilon}
$$
et finalement
$$
\prod_{p|d} \left|\sigma_p^{\boldsymbol{\varepsilon}}(\mathbf{e},\mathbf{E},e)\right| \ll d^{-\frac{1}{6}+\varepsilon}.
$$
\hfill
$\square$\\
\newline
En raisonnant alors comme dans la section 7 de \cite{36} et en utilisant le Lemme \ref{lemme22}, on obtient 
$$
E^{\boldsymbol{\varepsilon}}(B;\mathbf{e},\mathbf{E},e) \ll
(de)^{\varepsilon}\gcd(d,e)a'(\mathbf{E},\mathbf{\Delta})  \left( \frac{1}{d^2e^2}+\frac{1}{d^{\frac{9}{8}}e^{\frac{3}{2}}}+\frac{1}{d^{\frac{5}{4}}e^2} \right),
$$
soit
$$
E^{\boldsymbol{\varepsilon}}(B;\mathbf{e},\mathbf{E},e) \ll
(de)^{\varepsilon}\gcd(d,e)a'(\mathbf{E},\mathbf{\Delta})\frac{1}{d^{\frac{9}{8}}e^{\frac{3}{2}}}.
$$
Ensuite, puisque 
$
a'(\mathbf{E},\mathbf{\Delta})  
$
concerne les formes primitives, on n'a en réalité pas de dépendance en $e$ et par définition de cette quantité en (\ref{a}), on obtient que 
$$
a'(\mathbf{E},\mathbf{\Delta}) \ll 1
$$
où la constante est indépendante des paramètres de sommation. Enfin, le fait que
$$
\sum_{e=1}^{+\infty} \sum_{d \in \mathcal{D}} (de)^{\varepsilon}\gcd(d,e)\frac{1}{d^{\frac{9}{8}}e^{\frac{3}{2}}} \ll 1.
$$
pour $\varepsilon$ assez petit nous autorise à appliquer le théorème de convergence dominée pour obtenir le Théorème \ref{theor5}.

\section{La constante de Peyre}

Pour conclure la démonstration du Théorème \ref{theor1}, il reste à montrer que la constante $c_0$ obtenue dans le Théorème \ref{theor5} est en accord avec la conjecture de Peyre \cite[formule 5.1]{P03}, puisqu'on peut montrer que les surfaces de Châtelet sont des variétés ``presque de Fano" au sens de \cite[Définition 3.1]{P03}. On note $c_S$ la constante conjecturée par Peyre. Suivant \cite{P95} et~\cite{Sal}, on a
$$
c_S=\alpha(S) \beta(S) \omega_H\left(S(\mathbb{A}_{\mathbb{Q}})^{{\rm{Br}}(S)}\right),
$$
où
$$
\beta(S)=\#H^1(\mbox{Gal}(\overline{\mathbb{Q}},\mathbb{Q}),\mbox{Pic}(\overline{S}))=\mbox{Coker}\left( \mbox{Br}(\mathbb{Q}) \rightarrow \mbox{Br}(S) \right),
$$
$\alpha(S)$ est le volume d'un certain polytope dans le dual du cône effectif et $\omega_H\left(S(\mathbb{A}_{\mathbb{Q}})^{{\rm{Br}}(S)}\right)$ est un nombre de Tamagawa.

\subsection{Les facteurs $\alpha(S)$ et $\beta(S)$}

Pour commencer, il est nécessaire de rappeler quelques éléments relatifs à la géométrie de $S$ tirés de \cite{CoSSWD1} et \cite{CoSSWD2}. On considère $K=\mathbb{Q}(\sqrt{\Delta})$ le corps de décomposition de $F$ ainsi que le corps biquadratique $L=\mathbb{Q}\left(\sqrt{\Delta}, i\right)$. Si l'on note $\mathcal{G}=\mbox{Gal}\left(L / \mathbb{Q} \right)$, on a clairement que $\mathcal{G} \cong \left( \mathbb{Z}/2\mathbb{Z}\right)^2$ engendré par la conjugaison complexe~$\sigma$ et par~$\tau$, la conjugaison dans $K$. Enfin, pour toute extension $k$ de $\mathbb{Q}$, on notera $S_{k}=S \times_{\mathbb{Q}} \mbox{Spec}(k)$.
\par
Sur $\overline{\mathbb{Q}}$, on écrit également
$
F_3(u,v)=(a'_3u+b'_3 v)(a'_4u+b'_4 v)$, avec $a'_3$, $a'_4$, $b'_3$ et $b'_4$ dans $K$ tels que $\tau(a'_3)=a'_4$ et $\tau(b'_3)=b'_4$. En particulier, si l'on définit
\begin{eqnarray}
\Delta_{34}=a'_3b'_4-a'_4b'_3\neq 0,
\label{d34}
\end{eqnarray}
on a $\tau(\Delta_{34})=-\Delta_{34}$. Sur $L$, la surface $S_{L}$ admet 10 diviseurs exceptionnels qui sont des courbes d'auto-intersection négative. Huit d'entre elles sont d'auto-intersection $-1$ et sont données pour $i \in \{1,2 \}$ par
$$
D_i^+: \quad u=-b_i, \quad x+iy=0; \qquad D_i^-: \quad u=-b_i, \quad x-iy=0
$$
et pour $i \in \{3,4 \}$ par
$$
D_i^+: \quad u=-b'_i, \quad x+iy=0; \qquad D_i^-: \quad u=-b'_i, \quad x-iy=0
$$
tandis que les deux dernières sont d'auto-intersection $-2$ et sont données par
$$
E^+: \quad t=0, \quad x+iy=0; \qquad
E^-: \quad t=0, \quad x-iy=0.
$$
On a également
$$
\mbox{Pic}(S_{L})=\langle [E^{\pm}],[D_i^{\pm}]\mid i\in \{1,2,3,4 \}  \rangle=\langle [E^{+}],[D_1^{+}],[D_2^{+}],[D_3^{+}],[D_4^{+}],[D_1^{-}]  \rangle \cong \mathbb{Z}^6
$$
avec les relations
\begin{eqnarray}
[D_i^{+}]+[D_i^{-}]=[D_j^{+}]+[D_j^{-}]
\label{relations1}
\end{eqnarray}
pour $i,j \in \{1,2,3,4 \}$ et 
\begin{eqnarray}
[E^+]+[D_i^{+}]+[D_j^{+}]=[E^-]+[D_{\ell}^{-}]+[D_m^{-}]
\label{relations2}
\end{eqnarray}
pour $\{i,j,\ell,m\}=\{1,2,3,4\}$ et on rappelle que la formule d'adjonction fournit la relation
\begin{eqnarray}
\omega_S^{-1}=2E^++\sum_{i=1}^4D_i^+=2E^-+\sum_{i=1}^4D_i^-.
\label{antic}
\end{eqnarray}
Or, $\mbox{Pic}(S_{\overline{\mathbb{Q}}})\cong \mbox{Pic}(S_{L})$ et $\mbox{Pic}(S)=\left(\mbox{Pic}(S_{L})\right)^{\mathcal{G}}$ puisque $S(\mathbb{Q}) \neq \emptyset$. Décrivons alors l'action de $\mathcal{G}$ sur le groupe de Picard géométrique. On a
$$
\sigma(E^+)=E^-; \qquad 
\sigma(D_i^+)=D_i^-
$$
pour tout $i \in \{1,2,3,4\}$ lorsque $\Delta>0$ tandis qu'on a
$$
\sigma(E^+)=E^-; \qquad 
\forall i \in \{1,2\} \quad \sigma(D_i^+)=D_i^-; \qquad
\sigma(D_3^+)=D_4^- \quad \mbox{et} \quad \sigma(D_3^-)=D_4^+
$$
 et
$$
\tau(D_3^+)=D_4^+;\qquad
\tau(D_3^-)=D_4^-
$$
les autres élèments étant fixés. On en déduit les lemmes suivants.
\begin{lemme}
On  a $\mbox{Pic}(S)=\langle [\omega_S^{-1}], [D_1^+]+[D_1^-]\rangle$.
\label{lemmepic}
\end{lemme}
\noindent
\textit{Démonstration.}-- Soit 
$
D=a[D_1^+]+b[D_2^+]+c[D_3^+]+d[D_4^+]+e[E^+]+f[D_1^-] \in \mbox{Pic}(S_{\overline{\mathbb{Q}}}).
$
On a que $D \in \mbox{Pic}(S)$ si, et seulement si $\sigma(D)=\tau(D)=D$, soit si, et seulement si, 
$$
\begin{aligned}
D&=(b+c+d-e)([D_1^+]+[D_1^-])+f[D_1^+]+(a-e)[D_1^-]\\
&+(e-b)[D_2^+]+(e-c)[D_3^+]+(e-d)[D_4^+]+e[E^+]
\end{aligned}
$$
et
$
D=a[D_1^+]+b[D_2^+]+d[D_3^+]+c[D_4^+]+e[E^+]+f[D_1^-]
$ (ce qui implique notamment $c=d$).
Comme $[D_1^+]+[D_1^-] \in\mbox{Pic}(S)$, on en déduit que 
$D \in \mbox{Pic}(S)$ si, et seulement si $$
\begin{aligned}
\sigma\left(D-(b+c+d-e)([D_1^+]+[D_1^-])\right)&=\tau\left(D-(b+c+d-e)([D_1^+]+[D_1^-])\right)\\
&=D-(b+c+d-e)([D_1^+]+[D_1^-]).\\
\end{aligned}
$$
Cela équivaut à $e=2b=2c=2d$ et $a=b+f$ si bien que
$$
D-(b+c+d-e)([D_1^+]+[D_1^-])=(a-2b)([D_1^+]+[D_1^-])+b\left[\omega_S^{-1}\right]
$$
et
$$
D-(b+c+d-e)([D_1^+]+[D_1^-]) \in \mathbb{Z}\left([D_1^+]+[D_1^-]\right)\oplus\mathbb{Z}\omega_S^{-1}.
$$
\hfill
$\square$\\
\begin{lemme}
On a $\alpha(S)=\frac{1}{2}.$
\label{lemmealpha}
\end{lemme}
\noindent
\textit{Démonstration.}-- Posons $e_1=\omega_S^{-1}$ et $e_2=[D_1^+]+[D_1^-]$. On sait que dans ce cas, le cône effectif~$\Lambda_{\rm{eff}}(S)$ est engendré par les sommes d'éléments d'orbites des courbes d'auto-intersections négatives. On a par conséquent ici
$$
\Lambda_{\rm{eff}}(S)=\langle [E^+]+[E^-],[D_1^+]+[D_1^-],[D_3^+]+[D_4^+]+[D_3^-]+[D_4^-] \rangle=\langle [E^+]+[E^-],[D_1^+]+[D_1^-]\rangle.
$$
On utilise alors la définition suivante de la constante $\alpha(S)$ donnée dans \cite{P95}
$$
\alpha(S)=\mbox{Vol}\left\{ x \in \Lambda_{\rm{eff}}(S)^{\vee} \quad | \quad \langle\omega_S^{-1},x\rangle=1 \right\},
$$
où la mesure sur l'hyperplan
$$
\mathcal{H}=\left\{ x \in \mbox{Pic}(S)^{\vee}\otimes_{\mathbb{Z}} \mathbb{R} \quad | \quad \langle\omega_S^{-1},x\rangle=1 \right\}
$$
est définie dans \cite{P95}. Le cône $\Lambda_{\rm{eff}}(S)$ est donc engendré par
$
e_1-2e_2$ et $e_2$
si bien que la constante $\alpha(S)$ est donnée par le volume de la région
$$
\{ (x,y) \in \mathbb{R}^2 \mid x=1, \quad y>0, \quad x-2y>0\}.
$$
Autrement dit, on obtient la longueur du segment $\left[0,\frac{1}{2}\right]$, soit $\alpha(S)=\frac{1}{2}$.
\hfill
$\square$\\
\newline
\indent
On aurait également pu montrer que le groupe de Galois $G=\mbox{Gal}(\overline{\mathbb{Q}}/\mathbb{Q})$ ne fixe aucune $(-1)$-courbe et que le cardinal du groupe de Weyl associé au type de singularité $2\mathbf{A}_1$ où les deux singularités sont conjuguées est égal à 2. La table 2 de \cite{DEJ} ainsi que \cite{DJT} nous permettait alors d'en conclure que $\alpha(S)=\frac{1}{2}$.\\
\indent
Pour calculer $\beta(S)$, on utilise la proposition 7.1.2 de \cite{Sk}, qui fournit $\beta(S)=2$ dans le cas considéré dans cet article.

\subsection{Torseurs versels}

D'après \cite[proposition 8.3]{Pey} ou \cite[proposition 2.1]{De}, on déduit que l'ensemble des classes d'isomorphisme de torseurs versels au-dessus de $S$ possédant au moins un point rationnel est fini, ce qui permet d'exhiber une partition finie de l'ensemble des points rationnels de $S$, indexée par toute famille de représentants de ces classes d'isomorphisme. Contrairement à \cite{35}, il est plus délicat dans le cas de cet article de déterminer explicitement un tel système de représentants.\\
\par
On considère l'ensemble
$$
\mathcal{B}=\left\{ \boldsymbol{\beta} \in \mathbb{Z}^2 \times \mathbb{Z}\left[ \Delta\right]^2 \hspace{1mm} \bigg| \hspace{1mm} 
\begin{array}{c}
\tau(\beta_3)=\beta_4, \quad \exists (\beta'_3,n) \in \mathbb{Z}^2 \quad \mbox{tels que} \quad \beta_3\beta_4=\beta'_3n\\[2mm]
 \mbox{avec} \quad   \sqrt{\beta_1\beta_2\beta'_3} \in \mathbb{Z} \quad \mbox{et} \quad n \in N_{\mathbb{Q}[i]/\mathbb{Q}}\left( \mathbb{Z}[i] \right)\\
\end{array}
\right\}
$$
ainsi que le sous-ensemble $\mathcal{B}_{M}^{\Sigma}$ des $\boldsymbol{\beta} \in \mathcal{B}$ pour lesquels il existe $(\boldsymbol{\varepsilon},\mathbf{m}) \in \Sigma \times M$ tels que
$$
\beta_1=\varepsilon_1m_1; \qquad \beta_2=\varepsilon_2m_2; \qquad \beta'_3=\varepsilon_3m_3. 
$$
Il est bon de noter que l'ensemble $\boldsymbol{\beta} \in \mathcal{B}_M^{\Sigma}$ est infini. Pour $\boldsymbol{\beta}\in \mathcal{B}$, on pose $\mathcal{T}_{\boldsymbol{\beta}}$ le sous-ensemble constructible de $\mathbf{A}^{10}_{\mathbb{Q}}=\mbox{Spec}\left(\mathbb{Q}\left[X_i,Y_i \mid 0 \leqslant i \leqslant 4 \right] \right)$ défini par les deux équations quadratiques invariantes sous le groupe de Galois $\mathcal{G}$ suivantes
$$
\phi_1^{\boldsymbol{\beta}}=\phi_2^{\boldsymbol{\beta}}=0
$$
avec
$$
\left\{
\begin{aligned}
\phi_1^{\boldsymbol{\beta}}=&a_2\beta_1\left(X_1^2+Y_1^2\right)-a_1\beta_2\left(X_2^2+Y_2^2\right)-\\[2mm]
&\frac{\Delta_{12}}{\Delta_{34}}\Bigg(a_4\beta_3\left( X_3^2+Y_3^2+\Delta\left(X_4^2+Y_4^2\right)+2\sqrt{\Delta}(X_3X_4+Y_3Y_4)\right)+\\
&a_3\beta_4\left( X_3^2+Y_3^2+\Delta\left(X_4^2+Y_4^2\right)-2\sqrt{\Delta}(X_3X_4+Y_3Y_4)\right) \Bigg) \\
\phi_2^{\boldsymbol{\beta}}=&b_2\beta_1\left(X_1^2+Y_1^2\right)-b_1\beta_2\left(X_2^2+Y_2^2\right)-\\[2mm]
&\frac{\Delta_{12}}{\Delta_{34}}\Bigg(b_4\beta_3\left( X_3^2+Y_3^2+\Delta\left(X_4^2+Y_4^2\right)+2\sqrt{\Delta}(X_3X_4+Y_3Y_4)\right)+\\
&b_3\beta_4\left( X_3^2+Y_3^2+\Delta\left(X_4^2+Y_4^2\right)-2\sqrt{\Delta}(X_3X_4+Y_3Y_4)\right) \Bigg) ,\\
\end{aligned}
\right.
$$
et les inégalités
$$
\forall i \neq j \in \llbracket 0,2 \rrbracket^2, \qquad\left((X_{i},Y_{i}),(X_j,Y_j)\right)\neq((0,0),(0,0))
$$
et
$$
\forall i \in \llbracket 0,2 \rrbracket, \qquad \left((X_{i},Y_{i}),(X_3,Y_3,X_4,Y_4)\right)\neq((0,0),(0,0,0,0)).
$$
Définissons à présent un morphisme $\pi_{\boldsymbol{\beta}}:\mathcal{T}_{\boldsymbol{\beta}} \rightarrow S$. Pour ce faire, comme dans \cite[section~4]{35}, il suffit de définir un morphisme $\hat{\pi_{\boldsymbol{\beta}}}:\mathcal{T}_{\boldsymbol{\beta}} \rightarrow \mathcal{T}_{{\rm{spl}}}$. Considérons alors une extension finie $k$ de~$\mathbb{Q}$ et $((x_i,y_i))_{0 \leqslant i \leqslant 4}$ dans $\mathcal{T}_{\boldsymbol{\beta}}(k)$. Il existe alors un couple $(u,v) \in k^2\smallsetminus\{(0,0)\}$ tel que
$$
\begin{aligned}
u=&\frac{1}{\Delta_{12}}\left(b_2\beta_1(x_1^2+y_1^2)-b_1\beta_2(x_2^2+y_2^2)\right)\\[2mm]
=&\frac{1}{\Delta_{34}}\Bigg(a_4\beta_3\left( x_3^2+y_3^2+\Delta\left(x_4^2+y_4^2\right)+2\sqrt{\Delta}(x_3x_4+y_3y_4)\right)+\\
&a_3\beta_4\left( x_3^2+y_3^2+\Delta\left(x_4^2+y_4^2\right)-2\sqrt{\Delta}(x_3x_4+y_3y_4)\right) \Bigg) \\
\end{aligned}
$$
et
$$
\begin{aligned}
v=&\frac{1}{\Delta_{12}}\left(a_2\beta_1(x_1^2+y_1^2)-a_1\beta_2(x_2^2+y_2^2)\right)\\[2mm]
=&\frac{1}{\Delta_{34}}\Bigg(b_4\beta_3\left( x_3^2+y_3^2+\Delta\left(x_4^2+y_4^2\right)+2\sqrt{\Delta}(x_3x_4+y_3y_4)\right)+\\
&b_3\beta_4\left( x_3^2+y_3^2+\Delta\left(x_4^2+y_4^2\right)-2\sqrt{\Delta}(x_3x_4+y_3y_4)\right) \Bigg) \\
\end{aligned}
$$
tels que
$$
F_1(u,v)=\beta_1(x_1^2+y_1^2), \qquad F_2(u,v)=\beta_2(x_2^2+y_2^2)
$$
et
$$
 F_3(u,v)=\beta_3\beta_4 \left[\left(x_3^2-y_3^2-\Delta(x_4^2-y_4^2)\right)^2+\left( x_3y_3-\Delta y_3y_4\right)^2\right].
$$
On remarque en particulier que la quantité $F_3(u,v)/(\beta_3\beta_4)$ est bien une somme de deux carrés mais une somme de deux carrés particulière puisqu'il s'agit d'une norme de $L$ sur $\mathbb{Q}$. On note alors $\alpha_{\boldsymbol{\beta}}$ la racine carrée positive de $\beta_1 \beta_2 \beta'_3$ et $n=z_{\boldsymbol{\beta}}\overline{z_{\boldsymbol{\beta}}}$ et on considère 
$$
\left\{
\begin{array}{l}
x+iy=z_{\boldsymbol{\beta}}\alpha_{\boldsymbol{\beta}}\left(z_0\right)^2\prod_{j=1}^4z_j\\[2mm]
x-iy=\overline{z_{\boldsymbol{\beta}}}\alpha_{\boldsymbol{\beta}}\left(\overline{z_0}\right)^2\prod_{j=1}^4\overline{z_j}\\[2mm]
t=z_0\overline{z_0}\\
\end{array}
\right.
$$
où l'on a posé $z_j=x_j+iy_j$ pour $j \in \{0,1,2\}$ et
$$
z_3=\left(x_3^2-y_3^2-\Delta(x_4^2-y_4^2)\right)+i\left( x_3y_3-\Delta y_3y_4\right).
$$
On a alors $x^2+y^2=t^2F(u,v)$ et $(x,y,t,u,v) \in \mathcal{T}_{{\rm{spl}}}(k)$, ce qui permet de définir $\pi_{\boldsymbol{\beta}}$.

\begin{lemme}
Pour tout $\boldsymbol{\beta} \in \mathcal{B}$, la variété $\mathcal{T}_{\boldsymbol{\beta}}$ équipée du morphisme $\pi_{\boldsymbol{\beta}}$ et de l'action naturelle de $T_{{\rm{NS}}}$ définie de la même façon que dans \cite[section 4]{35} est un torseur versel pour~$S$ et
$$
S(\mathbb{Q})=\bigsqcup_{\boldsymbol{\beta} \in \mathcal{B}_M^{\Sigma}}\pi_{\boldsymbol{\beta}}\left(\mathcal{T}_{\boldsymbol{\beta}}(\mathbb{Q})\right).
$$
\label{lemme73}
\end{lemme}
\noindent
\textit{Démonstration}--
On utilisera dans la suite les notations $\Delta_{ij}=\mbox{Res}(F_i,F_j)$ et $\Delta_{43}=-\Delta_{34}$ qui généralisent (\ref{disc}) et (\ref{d34}). En particulier, on a $\Delta_{ij}=-\Delta_{ji}$. D'après \cite[section 4]{35} et en utilisant le formalisme développé dans \cite[Example 2.4.3]{Pi}, on sait que sur $\overline{\mathbb{Q}}$, un anneau de Cox de $S$ est donné par
\begin{eqnarray}
\overline{R}=\overline{\mathbb{Q}}[Z_i^+,Z_i^- \mid 0 \leqslant i \leqslant 4]/\left(\Delta_{jk}Z_i^+Z_i^-+\Delta_{ki}Z_j^+Z_j^-+\Delta_{ij}Z_{k}^+Z_{k}^-\right)_{1 \leqslant i<j<k\leqslant 4}
\label{coxch}
\end{eqnarray}
où $\mbox{div}(Z_0^{\pm})=E^{\pm}$ et $\mbox{div}(Z_i^{\pm})=D_i^{\pm}$ pour $ i \in \llbracket 1,4 \rrbracket$. On a plus précisément que
$$
\left(\Delta_{jk}Z_i^+Z_i^-+\Delta_{ki}Z_j^+Z_j^-+\Delta_{ij}Z_{k}^+Z_{k}^-\right)_{1 \leqslant i<j<k\leqslant 4}=\left(P_{1,2,3},P_{1,2,4} \right)
$$
si l'on note $P_{i,j,k}$ la forme quadratique $\Delta_{jk}Z_i^+Z_i^-+\Delta_{ki}Z_j^+Z_j^-+\Delta_{ij}Z_{k}^+Z_{k}^-$. Considérons alors les variables invariantes sous $\mathcal{G}$ suivantes
$$
X_k=\frac{Z_k^++Z_k^-}{2} \quad \mbox{et} \quad Y_k=\frac{Z_k^+-Z_k^-}{2i}
$$
pour $k \in \{ 0,1,2\}$ et
$$
X_{3}=\frac{Z_3^++Z_4^++Z_3^-+Z_4^-}{4} \quad \mbox{et} \quad Y_{3}=\frac{Z_3^+-Z_3^-+Z_4^+-Z_4^-}{4i}
$$
et 
$$
X_{4}=\frac{Z_3^+-Z_4^++Z_3^--Z_4^-}{4\sqrt{\Delta}} \quad \mbox{et} \quad Y_{4}=\frac{Z_3^++Z_4^--Z_3^--Z_4^+}{4\sqrt{\Delta}i}.
$$
Par exemple, lorsque $a'_3\neq 0$, puisqu'on a 
$$\phi_1^{(1,1,1,1)}=a'_3P_{1,2,4}-a'_4P_{1,2,3} \quad \mbox{et} \quad \phi_2^{(1,1,1,1)}=\frac{1}{a'_3}\left(\Delta_{34} P_{1,2,3}+b'_3\phi_1^{(1,1,1,1)} \right),
$$ 
une descente galoisienne garantit qu'un anneau de Cox pour $S$ sur $\mathbb{Q}$ est donné par
$$
R=\mathbb{Q}[X_i,Y_i \mid 0 \leqslant i \leqslant 4]/(\phi_1^{(1,1,1,1)},\phi_2^{(1,1,1,1)}).
$$
Par \cite[corollary 2.3.9]{Sk}, un torseur versel est unique à twist près par un élément de
 $H^1(\mathbb{Q},\mbox{Pic}(S_{\overline{\mathbb{Q}}}))$. On peut alors montrer en adaptant la preuve de \cite[proposition 2.69]{Pi} ou en utilisant \cite[theorem 7.1]{CoSSWD2} que pour tout cocycle $c \in H^1(\mathbb{Q},\mbox{Pic}(S_{\overline{\mathbb{Q}}}))$, on obtient que l'anneau de Cox tordu par $c$ est de la forme
$$
R^{c}=\mathbb{Q}[X_i,Y_i \mid 0 \leqslant i \leqslant 4]/(\phi_1^{\boldsymbol{\beta}},\phi_2^{\boldsymbol{\beta}})
$$
pour un certain $\boldsymbol{\beta} \in \mathcal{B}$.\par
Finalement, il reste à établir le dernier point:
$$
S(\mathbb{Q})=\bigsqcup_{\boldsymbol{\beta} \in \mathcal{B}_M^{\Sigma}}\pi_{\boldsymbol{\beta}}\left(\mathcal{T}_{\boldsymbol{\beta}}(\mathbb{Q})\right).
$$
Considérons un point $P \in S(\mathbb{Q})$. Il existe alors un unique $(y,z,t,u,v) \in \mathbb{Z}^4$ tel que 
$$
\left\{
\begin{array}{l}
(y,z,t)=(u,v)=1\\   t>0\\   F_1(u,v) \geqslant 0\\
t^2F_1(u,v)F_2(u,v)F_3(u,v)=y^2+z^2
\end{array}
\right.
$$
et $\pi_{{\rm{spl}}}((y,z,t,u,v))=P$ avec $\pi_{{\rm{spl}}}$ l'application $\pi_{{\rm{spl}}}:\mathcal{T}_{{\rm{spl}}} \rightarrow S$ définie dans \cite[definition~4.1]{35}. Lorsque $F_1(u,v)F_2(u,v)F_3(u,v) > 0$,  il existe un unique triplet $(\varepsilon_1,\varepsilon_2,\varepsilon_3) \in \{-1,+1\}^3$ avec~$\varepsilon_1=1$ tels que $\varepsilon_1 \varepsilon_2 \varepsilon_3=1$ et
$$
\varepsilon_i F_i(u,v) >0.
$$
Le fait que $F_1(u,v)F_2(u,v)F_3(u,v)$ soit une somme de deux carrés implique que, lorsque $p \equiv 3\Mod{4}$, on a
$$
\nu_p(F_1(u,v)F_2(u,v)F_3(u,v))=\nu_p(F_1(u,v))+\nu_p(F_2(u,v))+\nu_p(F_3(u,v)) \equiv 0 \Mod{2}.
$$
On pose alors
$$
m_{1}=\prod_{p \equiv 3\Mod{4} \atop \nu_p(F_1(u,v)) \equiv 1\Mod{2}}p, \quad m_{2}=\prod_{p \equiv 3\Mod{4} \atop \nu_p(F_2(u,v)) \equiv 1\Mod{2}}p, \quad m_{3}=\prod_{p \equiv 3\Mod{4} \atop \nu_p(F_3(u,v)) \equiv 1\Mod{2}}p,
$$
de sorte que $m_{i} | F_i(u,v)$. On en déduit que si $p|m_{1}$, exactement une seule des deux valeurs $\nu_p(F_2(u,v))$ et $\nu_p(F_3(u,v))$ est impaire. On a donc que $p|m_{1}$ et $p|m_{2}$ par exemple et puisque~$(u,v)=1$, cela implique que $p| \Delta_{12}^{(3)}$, si bien qu'on en déduit les relations
$$
m_{1} \bigg| \left[\Delta_{12}^{(3)},\Delta_{13}^{(3)}\right], \quad m_{2} \bigg| \left[\Delta_{12}^{(3)},\Delta_{23}^{(3)}\right], \quad m_{3} \bigg| \left[\Delta_{13}^{(3)},\Delta_{23}^{(3)}\right],
$$
$$
(m_{i},m_{j}) \big| \Delta_{ij}^{(3)},
$$
et
$$
\nu_p\left( \frac{F_i(u,v)}{m_{i}} \right) \equiv 0\Mod{2},
$$
pour tous les nombres premiers $p \equiv 3\Mod{4}$. De plus, $m_{1}m_{2}m_{3}$ est un carré. Si par exemple $F_1(u,v)=0$, on pose $\varepsilon_1m_{1}$ comme étant l'unique entier sans facteur carré tel que~$\varepsilon_1\varepsilon_2\varepsilon_3m_{1}m_{2}m_{3}$ soit un carré.\par
On introduit alors $\beta_1=\varepsilon_1 m_1$ et $\beta_2=\varepsilon_2$. On pose alors $g$ le plus petit diviseur de $F_3(u,v)$ tel que
$
\frac{\varepsilon_3 F_3(u,v)}{g}
$
soit une norme de $L$. Nécessairement, en posant $\beta'_3=\varepsilon_3 m_3$, on a $\beta'_3 \mid g$ donc on peut écrire $g=\beta'_3 n$. On constate alors que, par définition de $\beta'_3$, que $n$ est une norme sur~$\mathbb{Q}[i]$. De plus, $g$ est une norme sur $K$, on peut donc écrire $g=\beta_3\beta_4=\beta'_3n$ avec~$\tau(\beta_3)=\beta_4$. Ainsi, $\boldsymbol{\beta} \in \mathcal{B}_{M}^{\Sigma}$ et $(y,z,t,u,v) \in \hat{\pi_{\boldsymbol{\beta}}}\left(\mathcal{T}_{\boldsymbol{\beta}} \right)$. On peut montrer que ce $\boldsymbol{\beta} \in \mathcal{B}_{M}^{\Sigma}$ est unique comme dans la preuve de \cite[Proposition 4.7]{35}.
\hfill
$\square$\\
\newline
On constate donc, en notant $\mathcal{X}_{\boldsymbol{\beta}}$ le sous-schéma de $\mathbb{A}^8_{\mathbb{Q}}=\mbox{Spec}\left( \mathbb{Q}[X_i,Y_i \mid 1\leqslant i \leqslant 4]\right)$ défini par les équations $\phi_1^{\boldsymbol{\beta}}$ et $\phi_2^{\boldsymbol{\beta}}$, que $\mathcal{T}_{\boldsymbol{\beta}}$ est égal au produit $\mathcal{X}_{\boldsymbol{\beta}}\times \mathbb{A}^2_{\mathbb{Q}}$. En notant $\mathcal{X}_{\boldsymbol{\beta}}^{\circ}$ le complémentaire de l'origine dans $\mathcal{X}_{\boldsymbol{\beta}}$, on a un isomorphisme entre l'intersection complète de $\mathbb{A}^{10}_{\mathbb{Q}}\smallsetminus\{0\}$ donnée par les équations
$$
F_1(u,v)=\beta_1(X_1^2+Y_1^2), \quad F_2(u,v)=\beta_2(X_2^2+Y_2^2)
$$
et
$$
F_3(u,v)=\beta_3\beta_4\left[ \left(X_{3}^2-Y_{3}^2-\Delta(X_4^2-Y_{4}^2)\right)^2+\left(X_{3}Y_{3}-\Delta X_{4}Y_{4}\right)^2 \right],
$$
et le schéma $\mathcal{X}_{\boldsymbol{\beta}}^{\circ}$. Or, lors de la preuve de la conjecture de Manin, on s'est ramené à un problème de comptage sur certaines variétés de la forme (\ref{var})
$$
F_i(u,v)=\beta_i(X_i^2+Y_i^2).
$$
On n'a par conséquent pas utilisé une descente sur les torseurs versels mais sur des torseurs d'un type différent que l'on explicite dans la section suivante lors de notre preuve de la conjecture de Manin. Cette description facilitera grandement le traitement de la constante afin d'établir que celle-ci correspond à la prédiction de Peyre. De plus, il apparaît difficile \textit{a priori} d'obtenir un système de représentant des classes d'isomorphisme de torseurs versels et de décrire précisément les ensembles du type $\pi_{\boldsymbol{\beta}}\left(\mathcal{T}_{\boldsymbol{\beta}}(\mathbb{Q})\right)$. En effet, pour ce faire, il faut caractériser les normes de l'extension biquadratique $L$.\par
Il est à noter plus généralement que dans les cas des articles \cite{35}, \cite{36} et \cite{39}, le même type de calculs fournissent que le seul cas de la conjecture de Manin pour les surfaces de Châtelet où l'on utilise une méthode de descente sur les torseurs versels est celui pour lequel $F$ est scindé. Dans tous les autres cas, on utilise une descente sur des torseurs d'un type différent dont la construction est analogue à celle de la section qui suit.

\subsection{Les torseurs utilisés dans la preuve de la conjecture de Manin}

On décrit explicitement dans cette section les torseurs qui sont utilisés dans la preuve du Théorème~\ref{theor1}, c'est-à-dire les torseurs qui correspondent aux variétés de la forme (\ref{var}). Cette description repose sur le formalisme développé dans \cite{Pi} et \cite{DeP}. \\
\par
Soit $T$ le tore algébrique dont le groupe des caractères est donné par
\begin{eqnarray}
\hat{T}=\left[E^+\right]\mathbb{Z} \oplus \left[D_1^+\right]\mathbb{Z} \oplus \left[D_2^+\right]\mathbb{Z} \oplus \left[D_3^++D_4^+\right]\mathbb{Z}\oplus \left[D_1^-\right]\mathbb{Z}
\label{torus}
\end{eqnarray}
ainsi que l'injection $\lambda: \hat{T} \hookrightarrow \mbox{Pic}(S_{\overline{\mathbb{Q}}})$.

\begin{lemme}
 La~$\mathbb{Q}-$algèbre
$
R':=\mathbb{Q}\left[X_i,Y_i \mid 0 \leqslant i \leqslant 3\right]/(f_1),
$
avec
$$
f_1=X_3^2+Y_3^2-\frac{1}{\Delta_{12}^2}\left( \Delta_{23}(X_1^2+Y_1^2)+\Delta_{31}(X_2^2+Y_2^2)  \right)\left( \Delta_{24}(X_1^2+Y_1^2)+\Delta_{41}(X_2^2+Y_2^2)  \right),
$$
est un anneau de Cox de $S$ de type $\lambda$ défini en (\ref{torus}).
\label{lemmecox}
\end{lemme}
\noindent
\textit{Démonstration}--
D'après (\ref{antic}), on a $\left[ \omega_S^{-1}\right] \in \lambda(\hat{T})$ et ainsi l'image de $\lambda$ contient la classe d'un diviseur ample. Suivant la preuve de \cite[proposition 2.71]{Pi},  un anneau de Cox $R'$ de type $\lambda$ est donné par l'anneau des invariants sous le groupe de Galois $\mathcal{G}$ de
$$
\bigoplus_{m \in \hat{T}} \overline{R}_m,
$$
où $\overline{R}$ a été défini (\ref{coxch}) et $\overline{R}_m$ correspond aux éléments homogènes de degré $m$ de $\overline{R}$. Supposons $m \in \hat{T}$ donné tel que
$$
m=\left[ a_0 E^++a_1D_1^++a_2D_2^++a_3(D_3^++D_4^+)+a_4 D_1^-\right] 
$$
avec $a_i \in \mathbb{Z}$. Pour déterminer $\overline{R}_m$, on cherche à résoudre le système linéaire donné par
\begin{small}
$$
\left[e_0^+E^++e_0^-E^-+\sum_{j=1}^4 \left(e_j^+ D_j^++e_j^-D_j^-\right) \right]=\left[ a_0 E^++a_1D_1^++a_2D_2^++a_3(D_3^++D_4^+)+a_4 D_1^-\right],
$$
\end{small}
où les $e_j^{\pm} \geqslant 0$. Grâce aux relations (\ref{relations1}) et (\ref{relations2}), ce dernier est équivalent à 
$$
\left\{
\begin{array}{l}
a_0={e_0^+}+{e_0^-}\\[1mm]
a_1=\sum_{j=2}^{4}{e_j^-}+{e_1^+}-{e_0^-}\\[1mm]
a_2=e_2^+-e_2^-+e_0^-\\
a_3=e_3^+-e_3^-+e_0^-=e_4^+-e_4^-+e_0^-\\[1mm]
a_4=\sum_{j=1}^{4}e_j^--2e_0^-\\
\end{array}
\right.
$$
En résolvant ce système en les $e_j^{\pm}$ pour $0 \leqslant j \leqslant 4$, il vient que $\overline{R'}$ est isomorphe au sous-anneau de $\overline{R}$ engendré par les variables
$$
Z_0^{\pm} ;\quad
 Z_1^{\pm} ;\quad
 Z_2^{\pm};\quad
Z_3^+Z_4^+;\quad
 Z_3^-Z_4^- 
$$ 
vérifiant la relation suivante (d'après (\ref{coxch}))
$$
\begin{aligned}
\left(Z_3^+Z_4^+\right)\left(Z_3^-Z_4^-\right)&=\left(Z_3^+Z_3^-\right)\left(Z_4^+Z_4^-\right)\\
&=\frac{1}{\Delta_{12}^2}\left( \Delta_{23}Z_1^+Z_1^-+\Delta_{31}Z_{2}^+Z_{2}^-  \right)\left( \Delta_{24}Z_1^+Z_1^-+\Delta_{41}Z_{2}^+Z_{2}^-  \right).\\
\end{aligned}
$$
Cette dernière relation est bien invariante par le groupe de Galois $\mathcal{G}$. En considérant les variables $\mathcal{G}-$invariantes suivantes
$$
X_k=\frac{Z_k^++Z_k^-}{2} \quad \mbox{et} \quad Y_k=\frac{Z_k^+-Z_k^-}{2i}
$$
pour $k \in \{0,1,2\}$ et
$$
X_3=\frac{Z_3^+Z_4^++Z_3^-Z_4^-}{2} \quad \mbox{et} \quad Y_3=\frac{Z_3^+Z_4^+-Z_3^-Z_4^-}{2i},
$$
il vient alors bien
$$
R' =\mathbb{Q}[X_i,Y_i \mid 0 \leqslant i \leqslant 3]/(f_1).
$$
\hfill
$\square$\\
\newline
\`A nouveau d'après \cite[corollary 2.3.9]{Sk}, un torseur de type $\lambda$ est unique à twist près par un élément de $H^1(\mathbb{Q},\hat{T})$. Le lemme suivant explicite ce groupe de cohomologie $H^1(\mathbb{Q},\hat{T})$.

\begin{lemme}
On a la série d'isomorphismes
$
H^1(\mathbb{Q},\hat{T}) \cong H^1\left(\mbox{Gal}\left(\mathbb{Q}[i]/\mathbb{Q}\right),\hat{T}\right) \cong  \mathbb{Z}/2\mathbb{Z}.
$
\label{lemmecoh}
\end{lemme}
\noindent
\textit{Démonstration}-- Le point clé concernant les torseurs de type $\lambda$ est que l'action de la conjugaison $\tau$ dans $K$ est triviale si bien que $H^1\left(\mbox{Gal}\left(\overline{\mathbb{Q}}/\mathbb{Q}[i]\right),\hat{T}\right)=\{0\}$. Ainsi, la suite exatce de restriction-inflation
$$
\xymatrix{
0 \ar[r] & H^1\left(\mbox{Gal}\left(\mathbb{Q}[i]/\mathbb{Q}\right),\hat{T}\right) \ar[r] & H^1(\mathbb{Q},\hat{T}) \ar[r] & H^1\left(\mbox{Gal}\left(\overline{\mathbb{Q}}/\mathbb{Q}[i]\right),\hat{T}\right)
}
$$
fournit l'isomorphisme $H^1(\mathbb{Q},\hat{T}) \cong H^1\left(\mbox{Gal}\left(\mathbb{Q}[i]/\mathbb{Q}\right),\hat{T}\right)$. Puisque le groupe $\mbox{Gal}\left(\mathbb{Q}[i]/\mathbb{Q}\right)$ est cyclique d'ordre 2 engendré par la conjugaison complexe $\sigma$, le groupe de cohomologie $H^1\left(\mbox{Gal}\left(\mathbb{Q}[i]/\mathbb{Q}\right),\hat{T}\right)$ coïncide avec l'homologie du complexe
$$
\xymatrix{
 \hat{T} \ar[r]^{{\rm{Id}}+\sigma} & \hat{T} \ar[r]^{{\rm{Id}}-\sigma} & \hat{T}.
}
$$
Autrement dit, 
$$
H^1\left(\mbox{Gal}\left(\mathbb{Q}[i]/\mathbb{Q}\right),\hat{T}\right) \cong \mbox{Ker}(\sigma+{\rm{Id}})/\mbox{Im}({\rm{Id}}-\sigma).
$$
On obtient alors aisément le fait que $\mbox{Ker}(\sigma+{\rm{Id}})$ est engendré par 
$$
\left[D_1^+\right]-\left[D_1^-\right]; \quad \left[D_1^+\right]-\left[D_2^+\right]; \quad 2\left[D_1^+\right]-\left[D_3^++D_4^+\right]
$$
et que $\mbox{Im}({\rm{Id}}-\sigma)$ est engendrée par
$$
\left[D_1^+\right]-\left[D_1^-\right]; \quad 2\left[D_2^+\right]-\left[D_1^+\right]-\left[D_1^-\right]; \quad 2\left[D_3^++D_4^+\right]-2\left[D_1^+\right]-2\left[D_1^-\right].
$$
Il vient alors immédiatement que $H^1\left(\mbox{Gal}\left(\mathbb{Q}[i]/\mathbb{Q}\right),\hat{T}\right) \cong \mathbb{Z}/2\mathbb{Z}$ engendré par la classe de~$\left[D_1^+\right]-\left[D_2^+\right]$.
\hfill
$\square$\\
\newline
On est désormais en mesure de décrire tous les anneaux de Cox de $S$ de type $\lambda$. On s'inspire ici des deux preuves de \cite[propositions 2.69-2.70]{Pi}. Tout élément de $H^1(\mathbb{Q},\hat{T})$ est représenté par un cocycle $c: \mbox{Gal}\left(\mathbb{Q}[i]/\mathbb{Q}\right) \rightarrow \hat{T}$, lui-même déterminé par l'image $c_{\sigma} \in \mbox{Hom}_{\mathbb{Z}}\left(T,\mathbb{Q}[i]^{\times}\right)$ de la conjugaison complexe.

\begin{lemme}
Tout anneau de Cox pour $S$ de type $\lambda$ est isomorphe à une $\mathbb{Q}-$algèbre de la forme
$$
R'_{\mathbf{n}}:=\mathbb{Q}\left[X_i,Y_i \mid 0 \leqslant i \leqslant 3\right]/(f_1^{\mathbf{n}})
$$
avec
\begin{eqnarray}
\begin{aligned}
f_1^{\mathbf{n}}&=n_3\left(X_3^2+Y_3^2\right)-\frac{1}{\Delta_{12}^2}\left( \Delta_{23}n_1(X_1^2+Y_1^2)+\Delta_{31}n_2(X_2^2+Y_2^2)  \right)\\
&\times\left( \Delta_{24}n_1(X_1^2+Y_1^2)+\Delta_{41}n_2(X_2^2+Y_2^2)  \right)
\end{aligned}
\label{fn}
\end{eqnarray}
et $\mathbf{n}=(n_1,n_2,n_3) \in \left(\mathbb{Z}\smallsetminus \{0\}\right)^3$ tels que $n_1n_2n_3$ soit une somme de deux carrés.
\label{lemmetordu}
\end{lemme}
\noindent
\textit{Démonstration}-- On sait que, pour tout anneau de Cox $R'_{\lambda}$ pour $S$ de type $\lambda$, il existe un cocycle $c: \mbox{Gal}\left(\mathbb{Q}[i]/\mathbb{Q}\right) \rightarrow \hat{T}$ tel que $R'_{\lambda}$ soit un twist de $R'$ défini en (\ref{lemmecox}). D'après \cite[proposition 2.41]{Pi}, l'action de la conjugaison complexe sur $\overline{R'}$ est la suivante:
$$
\sigma(Z_0^-)=c_{\sigma}\left( \left[E^+ \right]\right)Z_0^+; \quad \forall j \in \{1,2\} \quad \sigma(Z_j^-)=c_{\sigma}\left( \left[D_j^+ \right]\right)Z_j^+;$$
et
$$
\sigma(Z_3^-Z_4^-)=c_{\sigma}\left( \left[D_3^++D_4^+ \right]\right)Z_3^+Z_4^+.
$$
En prenant les nouvelles variables suivantes, invariantes sous l'action du groupe de Galois $\mbox{Gal}\left(\mathbb{Q}[i]/\mathbb{Q}\right)$
$$
X_0=\frac{c_{\sigma}\left( \left[E^+ \right]\right)Z_0^++Z_0^-}{2}; \quad Y_0=\frac{c_{\sigma}\left( \left[E^+ \right]\right)Z_0^+-Z_0^-}{2i}; $$
$$
X_k=\frac{c_{\sigma}\left( \left[D_k^+ \right]\right)Z_k^++Z_k^-}{2}; \quad Y_k=\frac{c_{\sigma}\left( \left[D_k^+ \right]\right)Z_k^+-Z_k^-}{2i}
$$
pour $k \in \{ 1,2\}$ et
$$
X_3=\frac{c_{\sigma}\left( \left[D_3^++D_4^+ \right]\right)Z_3^+Z_4^++Z_3^-Z_4^-}{2} \quad \mbox{et} \quad Y_k=\frac{c_{\sigma}\left( \left[D_3^++D_4^+ \right]\right)Z_3^+Z_4^+-Z_3^-Z_4^-}{2i},
$$
il vient
$$
R'_{\lambda} \cong \mathbb{Q}\left[X_i,Y_i \mid 0 \leqslant i \leqslant 3\right]/(f_1^{c})
$$
avec
$$
\begin{aligned}
f_1^{c}&=n_{31}^c\left(X_3^2+Y_3^2\right)-\frac{1}{\Delta_{12}^2}\left( \Delta_{23}(X_1^2+Y_1^2)+\Delta_{31}n_{21}^c(X_2^2+Y_2^2)  \right)\\
&\times\left( \Delta_{24}(X_1^2+Y_1^2)+\Delta_{41}n_{21}^c(X_2^2+Y_2^2)  \right)
\end{aligned}
$$
et $n_{21}^c=c_{\sigma}\left( \left[D_1^+-D_2^+ \right]\right)$ et $n_{31}^c=c_{\sigma}\left( \left[2D_1^+-D_3^+-D_4^+ \right]\right)$. Le fait que $c$ soit un cocycle et les relations (\ref{relations1}) et (\ref{relations2}) impliquent que $\sigma(n_{i1})=n_{i1}$ si bien $n_{i1} \in \mathbb{Q}^{\times}$ pour $i \in \{2,3\}$. En écrivant $n_{21}^c=n_2/n_1$ et $n_{31}^c=n_3/{n_1}^2$ pour trois entiers $n_1$, $n_2$ et $n_3$ non nuls, on obtient bien
$$
R'_{\lambda} \cong \mathbb{Q}\left[X_i,Y_i \mid 0 \leqslant i \leqslant 3\right]/(f_1^{\mathbf{n}}).
$$
Pour conclure, il reste à établir le fait que $n_1n_2n_3$ est une somme de deux carrés. On a
$$
n_1n_2n_3=n_1^4 n_{21}^cn_{31}^c.
$$
Puisque les conditions de cocycle fournissent
$$
c_{\sigma}\left( \left[E^+ \right]\right) \sigma \left(c_{\sigma}\left( \left[E^- \right]\right) \right)=1; \quad \forall k \in \{1,2\} \quad c_{\sigma}\left( \left[D_k^+ \right]\right) \sigma \left(c_{\sigma}\left( \left[D_k^- \right]\right) \right)=1
$$
et
$$
c_{\sigma}\left( \left[D_3^++D_4^+ \right]\right) \sigma \left(c_{\sigma}\left( \left[D_3^-+D_4^- \right]\right) \right)=1,
$$
les relations (\ref{relations1}) et (\ref{relations2}) fournissent
$$
n_1n_2n_3=n_1^4c_{\sigma}\left( \left[E^++2D_1^+ \right]\right) \sigma \left(c_{\sigma}\left( \left[E^++2D_1^+ \right]\right) \right)
$$
si bien que $n_1n_2n_3$ est bien une somme de deux carrés. On peut établir plus précisément, en suivant \cite[proposition 2.70]{Pi} qu'étant donné  $(n_1,n_2,n_3) \in \mathbb{Q}^{\times}$, il existe un cocycle $c: \mbox{Gal}\left(\mathbb{Q}[i]/\mathbb{Q}\right) \rightarrow \hat{T}$ si, et seulement si, le produit $n_1n_2n_3$ est une somme de deux carrés.
\hfill
$\square$\\
\par
De façon analogue au Lemme \ref{lemme73}, on voit que pour $\mathbf{n}$ tel que $n_1n_2n_3$ soit une somme de deux carrés, le sous-ensemble constructible $\mathcal{T}_{\mathbf{n}}$ de $\mathbf{A}^{8}_{\mathbb{Q}}=\mbox{Spec}\left(\mathbb{Q}\left[X_i,Y_i \mid 0 \leqslant i \leqslant 3 \right] \right)$ défini par l'équation $f^{\mathbf{n}}$ donnée en (\ref{fn}) et les inégalités
$$
\forall i \neq j \in \llbracket 0, 3 \rrbracket^4, \quad \left((X_i,Y_i),(X_j,Y_j) \right)\neq \left((0,0),(0,0) \right)
$$
est un torseur de type $\lambda$ au-dessus de $S$. Pour toute extension finie $k$ de $\mathbb{Q}$ et tout $(x_i,y_i)_{0\leqslant i \leqslant 3}$ dans $\mathcal{T}_{\mathbf{n}}(k)$, à l'aide de
$$
u=\frac{1}{\Delta_{12}}\left(b_2n_1(x_1^2+y_1^2)-b_1n_2(x_2^2+y_2^2) \right)
$$
et
$$
v=\frac{1}{\Delta_{12}}\left(-a_2n_1(x_1^2+y_1^2)+a_1n_2(x_2^2+y_2^2) \right)
$$
on peut définir comme en section précédente un morphisme $\pi_{\mathbf{n}}:\mathcal{T}_{\mathbf{n}} \rightarrow S$. Lorsque $\mathbf{n}$ est de la forme
$$
\forall i \in \llbracket 1,3 \rrbracket, \quad n_i=\varepsilon_i m_i
$$
pour un certain $(\boldsymbol{\varepsilon},\mathbf{m}) \in \Sigma \times M$, on note $\mathcal{T}_{\mathbf{n}}=\mathcal{T}_{\mathbf{m},\boldsymbol{\varepsilon}}$ et $\pi_{\mathbf{n}}=\pi_{\mathbf{m},\boldsymbol{\varepsilon}}$. On peut alors montrer comme lors de la preuve du Lemme \ref{lemme73} qu'on a l'égalité
\begin{eqnarray}
S(\mathbb{Q})=\bigsqcup_{(\boldsymbol{\varepsilon},\mathbf{m}) \in \Sigma\times M}\pi_{\mathbf{m},\boldsymbol{\varepsilon}}\left(\mathcal{T}_{\mathbf{m},\boldsymbol{\varepsilon}}(\mathbb{Q})\right).
\label{uniondis}
\end{eqnarray}
De plus, on a
\begin{small}
\begin{eqnarray}
\begin{aligned}
\pi_{\mathbf{m},\boldsymbol{\varepsilon}}\left(  \mathcal{T}_{\mathbf{m},\boldsymbol{\varepsilon}}(\mathbb{Q}) \right)=&\left\{ [tu^2:tuv:tv^2:x_3:x_4] \in \mathbb{P}^4_{\mathbb{Q}} \left|
\begin{array}{c}
(t,u,v,x_3,x_4) \in \mathbb{Z}^5, (t,u,v)=(x_3,x_4)=1\\
 t\geqslant 0, x_3^2+x_4^2=t^2F(u,v), \\
\varepsilon_i F_i(u,v)>0\\
\nu_p(F_i(u,v))-\mu_i\equiv 0\Mod{2}\\
F_i(u,v)\in \varepsilon_i m_i \mathcal{E}
\end{array}
\right.
 \right\}\\[2mm]
\bigsqcup&\left\{ [tu^2:tuv:tv^2:x_3:x_4] \in \mathbb{P}^4_{\mathbb{Q}} \left|
\begin{array}{c}
(t,u,v,x_3,x_4) \in \mathbb{Z}^5, (t,u,v)=(x_3,x_4)=1\\
 t\geqslant 0, x_3^2+x_4^2=t^2F(u,v), \\
\varepsilon_i F_i(-u,-v)>0\\
\nu_p(F_i(u,v))-\mu_i\equiv 0\Mod{2}\\
F_i(-u,-v)\in \varepsilon_i m_i \mathcal{E}
\end{array}
\right.
 \right\}.\\
\end{aligned}
\label{formulemoche}
\end{eqnarray}
\end{small}
En effet, un point de $S$ admet exactement deux représentants de la forme $[tu^2:tuv:tv^2:x_3:x_4]$ avec $(t,u,v,x_3,x_4) \in \mathbb{Z}^5$ et  $(t,u,v)=(x_3,x_4)=1$.\par
\`A nouveau de la même façon qu'en section précédente, on montre, en notant $\mathcal{Y}_{\mathbf{n}}$ le sous-schéma de $\mathbf{A}^{8}_{\mathbb{Q}}=\mbox{Spec}\left(\mathbb{Q}\left[X_i,Y_i \mid 1 \leqslant i \leqslant 3 \right] \right)$ défini par l'équation $f^{\mathbf{n}}$ donnée en (\ref{fn}), que~$\mathcal{T}_{\mathbf{n}}$ est égal au produit $\mathcal{Y}_{\mathbf{n}} \times \mathbf{A}^{2}_{\mathbb{Q}}$. Si $\mathcal{Y}_{\mathbf{n}} ^{\circ}$ le complémentaire de l'origine dans $\mathcal{Y}_{\mathbf{n}} $, on a alors un isomorphisme entre l'intersection complète de $\mathbf{A}^{8}_{\mathbb{Q}}\smallsetminus \{0\}$ donnée par les équations
$$
F_i(u,v)=n_i(X_i^2+Y_i^2).
$$
On retrouve ainsi une variété de la forme (\ref{var}) et il s'agit donc bel et bien des torseurs utilisés dans la preuve de la conjecture de Manin en section 6.\par
On donne alors deux derniers lemmes qui permettent d'exprimer la constante conjecturée par Peyre de manière adéquate à notre traitement du problème de comptage.

\begin{lemme}
On a
$
{\mbox{\cyr SH}}^1(\mathbb{Q},T)=\{0\}
$
où
$$
{\mbox{\cyr SH}}^1(\mathbb{Q},T)=\mbox{Ker}\left( H^1(\mathbb{Q},T)   \longrightarrow H^1(\mathbb{R},T)\prod_p H^1(\mathbb{Q}_p,T) \right).
$$
\label{lemmesha}
\end{lemme}
\noindent
\textit{Démonstration}-- 
Le théorème de Tate fournit l'accouplement naturel non dégénéré suivant
$$
{\mbox{\cyr SH}}^1(\mathbb{Q},T) \times {\mbox{\cyr SH}}^2(\mathbb{Q},\hat{T}) \longrightarrow \mathbb{Q}/\mathbb{Z}
$$
si bien qu'il suffit d'établir que ${\mbox{\cyr SH}}^2(\mathbb{Q},\hat{T})=\{0\}$. Puisque
$
H^1\left(\mbox{Gal}\left(\overline{\mathbb{Q}}/\mathbb{Q}[i]\right),\hat{T}\right)=\{0\},
$
la suite de restriction-inflation d'ordre 2 fournit
$$
\xymatrix{
0 \ar[r] & H^2\left(\mbox{Gal}\left(\mathbb{Q}[i]/\mathbb{Q}\right),\hat{T}\right) \ar[r] & H^2(\mathbb{Q},\hat{T}) \ar[r] & H^2\left(\mbox{Gal}\left(\overline{\mathbb{Q}}/\mathbb{Q}[i]\right),\hat{T}\right)
}
$$
Montrons à présent que $H^2\left(\mbox{Gal}\left(\mathbb{Q}[i]/\mathbb{Q}\right),\hat{T}\right)=\{0\}$. Le groupe $\mbox{Gal}\left(\mathbb{Q}[i]/\mathbb{Q}\right)$ étant cyclique engendré par la conjugaison complexe $\sigma$, il vient que
$$
H^2\left(\mbox{Gal}\left(\mathbb{Q}[i]/\mathbb{Q}\right),\hat{T}\right)=\hat{T}^{\langle \sigma \rangle}/\mbox{Im}(\mbox{Id}+\sigma)=\mbox{Ker}(\sigma-\mbox{Id})/\mbox{Im}(\mbox{Id}+\sigma).
$$
On établit comme lors de la preuve du Lemme \ref{lemmecoh} que $\mbox{Ker}(\sigma-\mbox{Id})$ est de rang 2 engendré par
$
\left[D_1^+ \right]+\left[D_1^- \right]$ et $\left[\omega_S^{-1}\right].
$
Or, 
$$
\left[D_1^+ \right]+\left[D_1^- \right]=(\mbox{Id}+\sigma)\left( \left[D_1^+ \right]\right)\quad \mbox{et}  \quad \left[\omega_S^{-1}\right]=(\mbox{Id}+\sigma)\left(\left[E^+ \right] +2\left[D_1^+ \right]\right)
$$
si bien que $\mbox{Ker}(\sigma-\mbox{Id})=\mbox{Im}(\mbox{Id}+\sigma)$ et $H^2\left(\mbox{Gal}\left(\mathbb{Q}[i]/\mathbb{Q}\right),\hat{T}\right)=\{0\}$. On a donc le diagramme commutatif suivant
$$
\xymatrix{
&{\mbox{\cyr SH}}^2(\mathbb{Q},\hat{T}) \ar[d]&{\mbox{\cyr SH}}^2(\mathbb{Q}[i],\hat{T}) \ar[d]\\
0 \ar[r] &  H^2(\mathbb{Q},\hat{T})  \ar[r] \ar[d] & H^2\left(\mbox{Gal}\left(\overline{\mathbb{Q}}/\mathbb{Q}[i]\right),\hat{T}\right) \ar[d] \\
&  \prod_{v \in{\rm{val}}(\mathbb{Q})} H^2(\mathbb{Q}_v,\hat{T})  \ar[r] & \prod_{v \in{\rm{val}}(\mathbb{Q}[i])} H^2\left(\mathbb{Q}[i]_v,\hat{T}\right)
}
$$
qui implique que tout élément de ${\mbox{\cyr SH}}^2(\mathbb{Q},\hat{T})$ appartient à ${\mbox{\cyr SH}}^2(\mathbb{Q}[i],\hat{T})$. Or, sur $\mathbb{Q}[i]$, on a~$\hat{T} \cong \mathbb{Z}^5$ avec action triviale si bien que ${\mbox{\cyr SH}}^2(\mathbb{Q}[i],\hat{T})=\{0\}$ d'après \cite[lemme 1.9]{San}, ce qui permet de conclure la preuve.
\hfill
$\square$\\

\begin{lemme}
On a l'égalité
$$
S\left( \mathbb{A}_{\mathbb{Q}}\right)^{{\rm{Br}}(S)}=\bigsqcup_{(\mathbf{m},\boldsymbol{\varepsilon})\in M\times\Sigma}\pi_{\mathbf{m},\boldsymbol{\varepsilon}}\left(\mathcal{T}_{\mathbf{m},\boldsymbol{\varepsilon}}\left( \mathbb{A}_{\mathbb{Q}}\right)\right),
$$
où $\mathbb{A}_{\mathbb{Q}}$ désigne l'anneau des adèles de $\mathbb{Q}$.
\label{lemme77}
\end{lemme}
\noindent
\textit{Démonstration}-- 
Tout d'abord, d'après \cite[theorem 6.1.2]{Sk} et le Lemme \ref{lemmetordu}, il vient
$$
S\left( \mathbb{A}_{\mathbb{Q}}\right)^{{\rm{Br}}_{\lambda}(S)}=\bigcup_{\mathbf{n}\in \mathbb{N}^3 \atop n_1n_2n_3=\square+\square}\pi_{\mathbf{n}}\left(\mathcal{T}_{\mathbf{n}}\left( \mathbb{A}_{\mathbb{Q}}\right)\right),
$$
où ${\rm{Br}}_{\lambda}(S)=r^{-1}\left( \lambda_{\ast}\left( H^1(\mathbb{Q},\hat{T})\right)\right)$, $\lambda_{\ast}:H^1(\mathbb{Q},\hat{T}) \rightarrow H^1\left(\mathbb{Q},\mbox{Pic}(S_{\overline{\mathbb{Q}}})\right)$ est associée à l'injection~$\lambda$ et $r: {\rm{Br}}(S) \rightarrow H^1\left(\mathbb{Q},\mbox{Pic}(S_{\overline{\mathbb{Q}}})\right)$ est l'application canonique issue de la suite spectrale de Hochschild-Serre $H^p\left( \mathbb{Q},H^q\left( S_{\overline{\mathbb{Q}}},\mathbb{G}_m\right)\right)\Rightarrow H^{p+q}\left(S, \mathbb{G}_m\right)$. On a ici utilisé le fait que la surface~$S$ est géométriquement rationnelle auquel cas ${\rm{Br}}_1(S)={\rm{Br}}(S)$.\par
De la suite exacte courte
$$
\xymatrix{
 0 \ar[r] & \hat{T} \ar[r] & 
\mbox{Pic}(S_{\overline{\mathbb{Q}}}) 
\ar[r] & \mbox{Pic}(S_{\overline{\mathbb{Q}}})/\hat{T} 
\ar[r] & 0
}
$$
on tire la suite exacte
$$
\xymatrix{
 \left(\mbox{Pic}(S_{\overline{\mathbb{Q}}})/\hat{T} \right)^{\mathcal{G}}
\ar[r] & H^1(\mathbb{Q},\hat{T}) \ar[r]^{\lambda_{\ast} \hspace{5mm}} & H^1\left(\mathbb{Q},\mbox{Pic}(S_{\overline{\mathbb{Q}}})\right)
}
$$
Or, on montre facilement que $\mbox{Pic}(S_{\overline{\mathbb{Q}}})/\hat{T}$ est de rang 1 engendré par la classe de $\left[D_3^+\right]$ si bien que
$$
\left(\mbox{Pic}(S_{\overline{\mathbb{Q}}})/\hat{T} \right)^{\mathcal{G}}=\{0\}
$$
et $\lambda_{\ast}$ est injective. Le Lemme \ref{lemmecoh} et la valeur de $\beta(S)$ impliquent par conséquent 
$$
\lambda_{\ast}\left( H^1(\mathbb{Q},\hat{T})\right)=H^1\left(\mathbb{Q},\mbox{Pic}(S_{\overline{\mathbb{Q}}})\right)
$$
de sorte que ${\rm{Br}}_{\lambda}(S)={\rm{Br}}(S)$.\par
On a ainsi
$$
S\left( \mathbb{A}_{\mathbb{Q}}\right)^{{\rm{Br}}(S)}=\bigcup_{\mathbf{n}\in \mathbb{N}^3 \atop n_1n_2n_3=\square+\square}\pi_{\mathbf{n}}\left(\mathcal{T}_{\mathbf{n}}\left( \mathbb{A}_{\mathbb{Q}}\right)\right)
$$
et en particulier
$$
\bigcup_{\mathbf{m},\boldsymbol{\varepsilon}\in M\times\Sigma}\pi_{\mathbf{m},\boldsymbol{\varepsilon}}\left(\mathcal{T}_{\mathbf{m},\boldsymbol{\varepsilon}}\left( \mathbb{A}_{\mathbb{Q}}\right)\right) \subset S\left( \mathbb{A}_{\mathbb{Q}}\right)^{{\rm{Br}}(S)}.
$$
D'autre part, d'après \cite{CoSSWD1}, il vient $\overline{S(\mathbb{Q})}=S\left( \mathbb{A}_{\mathbb{Q}}\right)^{{\rm{Br}}(S)}$ et par continuité de l'application $\mathcal{T}_{\mathbf{m},\boldsymbol{\varepsilon}}(\mathbb{Q}) \rightarrow S(\mathbb{Q})$ induite par $\pi_{\mathbf{m},\boldsymbol{\varepsilon}}$ ainsi que (\ref{uniondis}), il s'ensuit
$$
\overline{S(\mathbb{Q})}=S\left( \mathbb{A}_{\mathbb{Q}}\right)^{{\rm{Br}}(S)}=\bigcup_{(\mathbf{m},\boldsymbol{\varepsilon})\in M\times\Sigma}\pi_{\mathbf{m},\boldsymbol{\varepsilon}}\left(\overline{\mathcal{T}_{\mathbf{m},\boldsymbol{\varepsilon}}\left( {\mathbb{Q}}\right)}\right) \subset\bigcup_{(\mathbf{m},\boldsymbol{\varepsilon})\in M\times\Sigma}\pi_{\mathbf{m},\boldsymbol{\varepsilon}}\left(\mathcal{T}_{\mathbf{m},\boldsymbol{\varepsilon}}\left( \mathbb{A}_{\mathbb{Q}}\right)\right)
$$
si bien qu'on peut conclure à l'égalité
$$
S\left( \mathbb{A}_{\mathbb{Q}}\right)^{{\rm{Br}}(S)}=\bigcup_{(\mathbf{m},\boldsymbol{\varepsilon})\in M\times\Sigma}\pi_{\mathbf{m},\boldsymbol{\varepsilon}}\left(\mathcal{T}_{\mathbf{m},\boldsymbol{\varepsilon}}\left( \mathbb{A}_{\mathbb{Q}}\right)\right).
$$
Suivant alors l'argument utilisé dans \cite{Sal} lors de la preuve du lemme 6.17, on a que tout point de $S\left( \mathbb{A}_{\mathbb{Q}}\right)^{{\rm{Br}}(S)}$ appartient à exactement $\#{\mbox{\cyr SH}}^1(\mathbb{Q},T)$ ensembles de la forme $\pi_{\mathbf{m},\boldsymbol{\varepsilon}}\left(\mathcal{T}_{\mathbf{m},\boldsymbol{\varepsilon}}\left( \mathbb{A}_{\mathbb{Q}}\right)\right)$ avec $(\mathbf{m},\boldsymbol{\varepsilon})\in M\times\Sigma$. Le Lemme \ref{lemmesha} permet alors d'en déduire que l'union est disjointe et de conclure la preuve.
\hfill
$\square$\\

\subsection{Expression de la constante de Peyre} Le Lemme \ref{lemme77} fournit
$$
c_S=\alpha(S)\beta(S) \sum_{\boldsymbol{\varepsilon} \in \Sigma \atop \mathbf{m} \in M}\omega_H\big(\pi_{\boldsymbol{\varepsilon},\mathbf{m}} (\mathcal{T}_{\boldsymbol{\varepsilon},\mathbf{m}}(\mathbb{A}_{\mathbb{Q}}))\big).
$$
On a $\alpha(S)\beta(S)=1$ et on écrit
$$
\omega_H(\pi_{\boldsymbol{\varepsilon},\mathbf{m}} (\mathcal{T}_{\boldsymbol{\varepsilon},\mathbf{m}}(\mathbb{A}_{\mathbb{Q}})))=\omega_{\infty}(\boldsymbol{\varepsilon},\mathbf{m}) \prod_p \omega_p(\boldsymbol{\varepsilon},\mathbf{m}),
$$
où, en passant comme dans \cite{36} aux densités sur le torseur intermédiaire $\mathcal{T}_{{\rm{spl}}}$ défini par l'équation (\ref{T}) et au vu de (\ref{formulemoche}), l'on a les égalités:
$$
\omega_{p}(\boldsymbol{\varepsilon},\mathbf{m})=\lim_{n \rightarrow +\infty}\frac{1}{p^{4n}} \# \left\{  (u,v,y,z,t) \in \left( \mathbb{Z}/p^n \mathbb{Z}\right)^5 \quad \Bigg| \quad 
\begin{array}{l}
 F(u,v)\equiv y^2+z^2\Mod{p^n}\\
 p\nmid (u,v), \quad p \nmid (y,z,t)\\
2|\nu_p(F_i(u,v))-\mu_{i}
\end{array}
\right\}
$$
pour $p\neq2$.
Pour tout entier $d$ et $F \in \mathbb{Z}[x_1,x_2]$, on introduit alors
$$
\rho^{\ast}(d)=\rho^{\ast}(d;F)=\#\{\mathbf{x} \in [0,d[^2  \quad | \quad d |F(\mathbf{x}), \quad (x_1,x_2,d)=1\}.
$$
Exactement le même type de raisonnement que dans la section 8.2 de \cite{39} conduit, lorsque $p \nmid \Delta_{12}^{(3)}\Delta_{13}^{(3)}\Delta_{23}^{(3)}$, aux expressions
$$
\omega_{p}(\boldsymbol{\varepsilon},\mathbf{m})=\left(1-\frac{1}{p^2}\right) \left(1-\frac{1}{p^2} +\sum_{ \nu \geqslant 1} \frac{\chi(p^{\nu})\rho^{\ast}(p^{\nu};F)}{p^{2\nu}} \right)
$$
lorsque $p \equiv 3\Mod{4}$ et
$$
\omega_{p}(\boldsymbol{\varepsilon},\mathbf{m})=\left( 1-\frac{1}{p^2} \right)^2+\left(1-\frac{1}{p}\right)^2\sum_{\nu \geqslant 1}\frac{\chi(p^{\nu})\rho^{\ast}(p^{\nu};F)}{p^{2\nu}}
$$
lorsque $p\equiv 1\Mod{4}$. En particulier, $\omega_{p}(\boldsymbol{\varepsilon},\mathbf{m})=\omega_p$ ne dépend ni de $\boldsymbol{\varepsilon}$ ni de $\mathbf{m}$. Dans le cas $p|\Delta_{12}^{(3)}\Delta_{13}^{(3)}\Delta_{23}^{(3)}$ et $p \equiv 3 \Mod{4}$, on a $(t,p)=1$ et donc on obtient 
$$
\omega_{p}(\boldsymbol{\varepsilon},\mathbf{m})=\lim_{n \rightarrow +\infty}\frac{1-\frac{1}{p}}{p^{3n}} \# \left\{  (u,v,y,z) \in \left( \mathbb{Z}/p^n \mathbb{Z}\right)^4 \quad \Bigg| \quad 
\begin{array}{l}
 F(u,v)\equiv y^2+z^2\Mod{p^n}\\
 p\nmid (u,v)\\
 2|\nu_p(F_i(u,v))-\mu_{i}
\end{array}
\right\}.
$$
Pour $p=2$, on a
$$
\begin{aligned}
\omega_{2}(\boldsymbol{\varepsilon},\mathbf{m})=&\lim_{n \rightarrow +\infty}\frac{1}{2^{4n}} \# \left\{  (u,v,y,z,t) \in \left( \mathbb{Z}/2^n \mathbb{Z}\right)^5 \quad \Bigg| \quad 
\begin{array}{l}
 F(u,v)\equiv y^2+z^2\Mod{2^n}\\
 2\nmid (u,v), \quad 2 \nmid (y,z,t)\\
F_i(u,v)\in \varepsilon_i m_i \mathcal{E}_{2^n}
\end{array}
\right\}\\[1cm]
&+\lim_{n \rightarrow +\infty}\frac{1}{2^{4n}} \# \left\{  (u,v,y,z,t) \in \left( \mathbb{Z}/2^n \mathbb{Z}\right)^5 \quad \Bigg| \quad 
\begin{array}{l}
 F(u,v)\equiv y^2+z^2\Mod{2^n}\\
 2\nmid (u,v), \quad 2 \nmid (y,z,t)\\
F_i(-u,-v)\in \varepsilon_i m_i \mathcal{E}_{2^n}
\end{array}
\right\}.\\
\end{aligned}
$$
On a également $(2,t)=1$ et donc, grâce au Lemme \ref{lemme16}, on aboutit à l'expression
$$
\omega_{2}(\boldsymbol{\varepsilon},\mathbf{m})=2\lim_{n \rightarrow +\infty}2^{-2n} \# \left\{  (u,v) \in \left( \mathbb{Z}/2^n \mathbb{Z}\right)^2 \quad \Bigg| \quad 
\begin{array}{l}
 2\nmid (u,v)\\
 F_i(u,v) \in \varepsilon_im_{i}\mathcal{E}_{2^n}
\end{array}
\right\}.
$$
Noter ici qu'on a une dépendance en $\boldsymbol{\varepsilon}$ et en $\mathbf{m}$. Enfin on traite le cas de la densité archimédienne. Grâce à la symétrie du problème, on peut se restreindre à $y>0$ et $z>0$. En utilisant une forme de Leray en paramétrant en $z$, on obtient
$$
\omega_{\infty}(\boldsymbol{\varepsilon},\mathbf{m})=4\lim_{B\rightarrow+\infty}\frac{1}{B\log(B)}\int_{D} \frac{\mbox{d}u\mbox{d}v \mbox{d}t\mbox{d}y}{2\sqrt{t^2F(u,v)-y^2}}
$$
avec
$$
D=\left\{ (u,v,y,t) \in \mathbb{R}^4 \quad \Bigg| \quad 
\begin{array}{c}
(u,v) \in \mathcal{R}_{(\varepsilon_1,\varepsilon_2,\varepsilon_3)}\left(\sqrt{B/t}\right)\bigsqcup \mathcal{R}_{(-\varepsilon_1,-\varepsilon_2,\varepsilon_3)}\left(\sqrt{B/t}\right), \\[3mm]
 0<y<t\sqrt{F(u,v)}, \quad 1<t<B\\
\end{array}
 \right\}.
$$
La formule (\ref{pi2}) fournit alors immédiatement
$$
\omega_{\infty}(\boldsymbol{\varepsilon},\mathbf{m})=\pi\left({\rm Vol}(\mathcal{R}_{(\varepsilon_1,\varepsilon_2,\varepsilon_3)}(1))+{\rm Vol}(\mathcal{R}_{(-\varepsilon_1,-\varepsilon_2,\varepsilon_3)}(1))\right)=2\pi\mbox{vol}(\mathcal{R}^{\boldsymbol{\varepsilon}}(1)),
$$
qui ne dépend que de $\boldsymbol{\varepsilon}$ mais pas de $\mathbf{m}$. On obtient ainsi l'égalité
\begin{eqnarray}
c_S=\sum_{\boldsymbol{\varepsilon} \in \Sigma \atop \mathbf{m} \in M} 2\pi\mbox{vol}(\mathcal{R}^{\boldsymbol{\varepsilon}}(1))\prod_p\omega_{p}(\boldsymbol{\varepsilon},\mathbf{m}) .
\label{const2}
\end{eqnarray}

\subsection{Transformation de la constante $c_0$}

On revient dans cette section à l'expression de la constante $c_0$ obtenue grâce au Théorème \ref{theor5} que l'on met sous une forme similaire à celle de $c_S$ en (\ref{const2}). On réécrit cette constante $c_0$ sous la forme suivante
\begin{eqnarray}
\begin{aligned}
c_0=& \frac{\pi^2}{2^4} \sum_{d \in \mathcal{D}} \frac{\mu(d)r_0(d)\varphi^{\dagger}(d)}{d}\sum_{\boldsymbol{\varepsilon} \in \Sigma} \mbox{vol}\left( R^{\boldsymbol{\varepsilon}}(1) \right)\sum_{\substack{\mathbf{d}, \mathbf{d}' \in \mathcal{D}^3 \\ d=d_1d_2d_3, d'_i| \Delta_{jk}}}\sum_{\mathbf{m}\in M}\mu(d'_1d'_2)\mu(d'_3)\\
& \times
\sum_{\substack{k_4k_1k'_1|\gcd(\Delta_{23},d) \\ k_4k_2k'_2|\gcd(\Delta_{13},d)\\k_4k_3k'_3|\gcd(\Delta_{12},d)}}\sum_{\substack{  k_4k'_4|\gcd(\Delta_{12},\Delta_{13},\Delta_{23},d)\\ k_5k'_5|\gcd(\Delta_{12},d'_1d'_2) }} \frac{\mu(k'_1)\mu(k'_2)\mu(k'_3)\mu(k'_4)\mu(k'_5)}{3^{\omega(k_4)}2^{\omega(k_5)+\omega(k_1)+\omega(k_2)+\omega(k_3)}} \sigma^{\boldsymbol{\varepsilon}}_{\ast}(\mathbf{e},\mathbf{E}),
\end{aligned}
\label{czero}
\end{eqnarray}
où
\begin{eqnarray}
\sigma^{\boldsymbol{\varepsilon}}_{\ast}(\mathbf{e},\mathbf{E})=\sum_{e=1}^{+\infty}\frac{\mu(e)}{e^2}\sigma^{\boldsymbol{\varepsilon}}(\mathbf{e},\mathbf{E},e).
\label{pl}
\end{eqnarray}
Par multiplicativité en $e$ de la quantité $\sigma^{\boldsymbol{\varepsilon}}(\mathbf{e},\mathbf{E},e)$, on peut écrire $\sigma^{\boldsymbol{\varepsilon}}_{\ast}(\mathbf{e},\mathbf{E})$ sous forme d'un produit eulérien
\begin{eqnarray}
\sigma^{\boldsymbol{\varepsilon}}_{\ast}(\mathbf{e},\mathbf{E})=\prod_p \sigma^{\boldsymbol{\varepsilon}}_{\ast,p}(\mathbf{e},\mathbf{E}).
\label{sigmaetoile}
\end{eqnarray}
Pour $p>2$, le facteur eulérien est donné par
$$
\sigma^{\boldsymbol{\varepsilon}}_{\ast,p}(\mathbf{e},\mathbf{E})=\left(1-\frac{\chi(p)}{p}\right)^3 \sum_{\boldsymbol{\nu} \in \mathbb{N}^3} \frac{\chi(p)^{\nu_1+\nu_2+\nu_3}\tilde{\rho}(p^{N_1},p^{N_2},p^{N_3};F_{1},F_{2},F_3)}{p^{2(N_1+N_2+N_3+1)}}
$$
où les $N_i$ sont définis lors du Lemme \ref{lemme21} et avec
$$
\tilde{\rho}(p^{N_1},p^{N_2},p^{N_3};F_{1},F_{2},F_3)=\#\{ \mathbf{x}\in \left( \mathbb{Z}/p^{N_1+N_2+N_3+1}\mathbb{Z} \right)^2 \quad | \quad p^{N_i}|F_i(\mathbf{x}), \quad p\nmid\mathbf{x}\},
$$
tandis que pour $p=2$, on obtient un facteur eulérien
\begin{eqnarray}
\sigma^{\boldsymbol{\varepsilon}}_{\ast,2}(\mathbf{e},\mathbf{E})=8\lim_{n \rightarrow +\infty} 2^{-2n}\#\left\{ 
 \mathbf{x} \in \left(\mathbb{Z}/2^n\mathbb{Z}\right)^2 \quad \Bigg| \quad 
 \begin{array}{l}
  F_i(\mathbf{x}) \in e_i\varepsilon_i\mathcal{E}_{2^n}\\
  2 \nmid\mathbf{x}
 \end{array}
 \right\}.
\label{sigma22}
\end{eqnarray}
On écrit alors $c_0$ sous la forme
$$
c_0=\sum_{\boldsymbol{\varepsilon} \in \Sigma \atop \mathbf{m} \in M} \mbox{vol}(\mathcal{R}_{\boldsymbol{\varepsilon}}(1)) c_0(\boldsymbol{\varepsilon},\mathbf{m})
$$
avec
$$
c_0(\boldsymbol{\varepsilon},\mathbf{m})=
\begin{aligned}[t]
& \frac{\pi^2}{2^4} \sum_{d \in \mathcal{D}} \frac{\mu(d)r_0(d)\varphi^{\dagger}(d)}{d}\sum_{\mathbf{d},\mathbf{d}' \in \mathcal{D}^3 \atop d=d_1d_2d_3, d'_i| \Delta_{jk}} 
\mu(d'_1d'_2)\mu(d'_3)\\
& \times
\sum_{\substack{k_4k_1k'_1|\gcd(\Delta_{23},d) \\ k_4k_2k'_2|\gcd(\Delta_{13},d)\\ k_4k_3k'_3|\gcd(\Delta_{12},d) }}\sum_{\substack{  k_4k'_4|\gcd(\Delta_{12},\Delta_{13},\Delta_{23},d)\\ k_5k'_5|\gcd(\Delta_{12},d'_1d'_2)  }} \frac{\mu(k'_1)\mu(k'_2)\mu(k'_3)\mu(k'_4)\mu(k'_5)}{3^{\omega(k_4)}2^{\omega(k_5)+\omega(k_1)+\omega(k_2)+\omega(k_3)}} \sigma^{\boldsymbol{\varepsilon}}_{\ast}(\mathbf{e},\mathbf{E}), 
\end{aligned}
$$
où
$\sigma^{\boldsymbol{\varepsilon}}_{\ast,p}(\mathbf{e},\mathbf{E})
$
a été défini en (\ref{pl}). De plus, en utilisant les notations (\ref{E}) et en posant
$$
e'_1=d_1d'_2d'_3 \quad e'_2=d_2d'_1d'_3 ,\quad e'_3=d_3d'_1d'_2,
$$
$$
E'_1=[e'_1,k_4k_2k'_2,k_4k_3k'_3,k_4k'_4,d_1k_5k'_5], \quad E'_2=[e'_2,k_4k_1k'_1,k_4k_3k'_3,k_4k'_4,d_2k_5k'_5],
$$
$$
E'_3=[e'_3,k_4k_1k'_1,k_4k_2k'_2,k_4k'_4],
$$
et
$$
e''_1=E''_1=m_{1} \quad e''_2=E''_2=m_{2}, \quad e''_3=E''_3=m_{3},
$$
on a les décompositions $e_i=e'_ie''_i$ et $E_i=E'_iE''_i$. On a alors $N_i=\max\{\nu_p(e'_i)+\nu_i,\nu_p(E'_i)\}$ lorsque $p \equiv 1 \Mod{4}$ et $N_i=\max\{\nu_p(e''_i)+\nu_i,\nu_p(E''_i)\}$ lorsque $p \equiv 3 \Mod{4}$.\\
Enfin, on obtient
$$
\sigma^{\boldsymbol{\varepsilon}}_{\ast,2}(\mathbf{e},\mathbf{E})=\sigma^{\boldsymbol{\varepsilon}}_{\ast,2}(\mathbf{m})
=8\lim_{n \rightarrow +\infty} 2^{-2n}\#\left\{ 
 \mathbf{x} \in \left(\mathbb{Z}/2^n\mathbb{Z}\right)^2 \quad \Bigg| \quad 
 \begin{array}{l}
  F_i(\mathbf{x}) \in m_{i}\varepsilon_i\mathcal{E}_{2^n}\\
  2 \nmid\mathbf{x}
 \end{array}
 \right\}.
$$
On pose alors
\begin{small}
$$
V_1(\mathbf{d},\mathbf{d}')=
\sum_{\substack{k_4k_1k'_1|\gcd(\Delta_{23},d) \\ k_4k_2k'_2|\gcd(\Delta_{13},d)\\  k_4k_3k'_3|\gcd(\Delta_{12},d) }}\sum_{\substack{   k_4k'_4|\gcd(\Delta_{12},\Delta_{13},\Delta_{23},d)\\ k_5k'_5|\gcd(\Delta_{12},d'_1d'_2) }} \frac{\mu(k'_1)\mu(k'_2)\mu(k'_3)\mu(k'_4)\mu(k'_5)}{3^{\omega(k_4)}2^{\omega(k_5)+\omega(k_1)+\omega(k_2)+\omega(k_3)}}\prod_{p \equiv 1\Mod{4}}\sigma_{\ast,p}(\mathbf{e}',\mathbf{E}'),
$$
\end{small}
puis
$$
V_3(\mathbf{m})=\prod_{p \equiv 3\Mod{4}}\sigma_{\ast,p}(\mathbf{e}'',\mathbf{E}'')
$$
et enfin
$$
V_2(\boldsymbol{\varepsilon},\mathbf{m})=\frac{1}{4}\sigma^{\boldsymbol{\varepsilon}}_{\ast,2}(\mathbf{m})
$$
de sorte qu'on ait
$$
c_0(\boldsymbol{\varepsilon},\mathbf{m})=2\pi\prod_p  \left(1+\frac{\chi(p)}{p}\right)\left(1-\frac{\chi(p)}{p}\right)^{-2}c_1 \times V_2(\boldsymbol{\varepsilon},\mathbf{m})V_3(\mathbf{m})
$$
avec
$$
c_1=
\sum_{d \in \mathcal{D}} \frac{\mu(d)r_0(d)\varphi^{\dagger}(d)}{d}\sum_{\mathbf{d}, \mathbf{d}' \in \mathcal{D}^3 \atop d=d_1d_2d_3, d'_i| \Delta_{jk}} \mu(d'_1d'_2)\mu(d'_3)V_1(\mathbf{d},\mathbf{d}').
$$
On a ici utilisé la décomposition en produit eulérien
$$
\frac{\pi}{8}=\prod_p \left(1+\frac{\chi(p)}{p}\right)\left(1-\frac{\chi(p)}{p}\right)^{-2}.
$$
La stratégie pour conclure à la validation de la conjecture de Peyre est alors, suivant la section 9.5 de \cite{39}, de décomposer $ c_0(\boldsymbol{\varepsilon},\mathbf{m})$ sous la forme
$$
c_0(\boldsymbol{\varepsilon},\mathbf{m})=2\pi\prod_p c_p(\boldsymbol{\varepsilon},\mathbf{m})
$$
et d'établir que pour tout nombre premier $p$, on a $c_p(\boldsymbol{\varepsilon},\mathbf{m})=\omega_p(\boldsymbol{\varepsilon},\mathbf{m}).$\\

\subsection{Fin de la preuve de la conjecture de Peyre}

On pose $c_2(\boldsymbol{\varepsilon},\mathbf{m})=V_2(\boldsymbol{\varepsilon},\mathbf{m})$ de sorte que $c_2(\boldsymbol{\varepsilon},\mathbf{m})=\omega_2(\boldsymbol{\varepsilon},\mathbf{m})$. Dans le cas $p \equiv 1 \Mod{4}$, on pose $\delta'_i=\nu_p(d'_{i})$
et
$$
N'_1=\nu_1+\delta_1+\delta'_2+\delta'_3, \quad 
N_1=\max\{N'_1,\kappa_4+\kappa_2+\kappa'_2,\kappa_4+\kappa_3+\kappa'_3,\kappa_4+\kappa'_4,\kappa_5+\kappa'_5\},
$$
et
$$
N'_2=\nu_2+\delta_2+\delta'_1+\delta'_3, \quad 
N_2=\max\{N'_2,\kappa_4+\kappa_1+\kappa'_1,\kappa_4+\kappa_3+\kappa'_3,\kappa_4+\kappa'_4,\kappa_5+\kappa'_5\}
$$
et enfin
$$
N'_3=\nu_3+\delta_3+\delta'_1+\delta'_2, \quad 
N_3=\max\{N'_3,\kappa_4+\kappa_1+\kappa'_1,\kappa_4+\kappa_2+\kappa'_2,\kappa_4+\kappa'_4\}.
$$
Puisque $r_0(p)=2$ et par définition de $\varphi^{\dagger}$, on a alors
$$
c_p(\boldsymbol{\varepsilon},\mathbf{m})=\left(1-\frac{1}{p^2}\right) \sum_{0 \leqslant \delta \leqslant 1}\left(-\frac{2}{p+1} \right)^{\delta} \sum_{0 \leqslant \delta_1,\delta_2, \delta_3 \leqslant 1 \atop \delta_1+\delta_2+\delta_3=\delta} f_p(\boldsymbol{\delta})
$$
avec
\begin{small}
\begin{eqnarray}
\begin{split}
f_p(\boldsymbol{\delta})=& \sum_{\substack{0 \leqslant \delta'_i \leqslant \min\{\nu_p(\Delta_{jk}),1 \}\\ \delta'_1\delta'_2=0}} \!\!(-1)^{\delta'_{1}+\delta'_{2}+\delta'_{3}} 
\sum_{\kappa_4,\kappa'_4 \geqslant 0 \atop \kappa_4+\kappa'_4 \leqslant \min\{\delta,\nu_p(\Delta_{12}),\nu_p(\Delta_{13}),\nu_p(\Delta_{23})\}}\!\! \frac{(-1)^{\kappa'_4}}{3^{\kappa_4}}\\
& \times
\sum_{i=1}^3\sum_{\kappa_i,\kappa'_i \geqslant 0 \atop \kappa_4+\kappa_i+\kappa'_i \leqslant \min\{\delta,\nu_p(\Delta_{jk})\}} \frac{(-1)^{\kappa'_i}}{2^{\kappa_i}}
\sum_{\kappa_5,\kappa'_5 \geqslant 0 \atop \kappa_5+\kappa'_5 \leqslant \min\{\delta'_1+\delta'_2,\nu_p(\Delta_{12})\}} \frac{(-1)^{\kappa'_5}}{2^{\kappa_5}}\\
& \times \sum_{\boldsymbol{\nu} \in \mathbb{N}^3} \frac{\tilde{\rho}\left(p^{N_1},p^{N_2},p^{N_3};F_1,F_2,F_3 \right)}{p^{2(N_1+N_2+N_3+1)}}, 
\label{ccc}
\end{split}
\end{eqnarray}
\end{small}
avec $\{i,j,k\}=\{1,2,3\}$. Enfin, dans le cas $p \equiv 3\Mod{4}$, si l'on note $\mu_{i}=\nu_p(m_{i})$, on pose
$$
c_p(\boldsymbol{\varepsilon},\mathbf{m})=\left(1-\frac{1}{p^2}\right)g(\boldsymbol{\mu}),
$$
où
$$g(\boldsymbol{\mu})=
\sum_{\boldsymbol{\nu} \in \mathbb{N}^3} \frac{(-1)^{\nu_1+\nu_2+\nu_3}\tilde{\rho}\left(p^{N_1},p^{N_2},p^{N_3};F_1,F_2,F_3 \right)}{p^{2(N_1+N_2+N_3+1)}}, 
$$
et avec $N_1=\mu_{1}+\nu_1$, $N_2=\mu_{2}+\nu_2$ et $N_3=\mu_{3}+\nu_3$.
On a ainsi
$$
c_0(\boldsymbol{\varepsilon},\mathbf{m})=2\pi\prod_p c_p(\boldsymbol{\varepsilon},\mathbf{m}).
$$
\par
Il s'agit à présent de remonter la formule d'éclatement du Lemme \ref{lemme19} utilisée dans le passage aux torseurs. On distingue selon la congruence de $p$ modulo 4. On commence par le cas $p \equiv 1 \Mod{4}$ et on étudie la quantité
$$
c'_p(\boldsymbol{\varepsilon},\mathbf{m})=\left(1-\frac{1}{p^2}\right)^{-1}c_p(\boldsymbol{\varepsilon},\mathbf{m}).
$$
Pour ce faire, on ne considère dans un premier temps que la quantité
$$
f'_p(\boldsymbol{\delta},\boldsymbol{\delta}')=
\begin{aligned}[t]
&
\sum_{\kappa_4,\kappa'_4 \geqslant 0 \atop \kappa_4+\kappa'_4 \leqslant \min\{\delta,\nu_p(\Delta_{ij})\}} \frac{(-1)^{\kappa'_4}}{3^{\kappa_4}}\sum_{\kappa_1,\kappa'_1 \geqslant 0 \atop \kappa_4+\kappa_1+\kappa'_1 \leqslant \min\{\delta,\nu_p(\Delta_{23})\}} \frac{(-1)^{\kappa'_1}}{2^{\kappa_1}}\sum_{\kappa_2,\kappa'_2 \geqslant 0 \atop \kappa_4+\kappa_2+\kappa'_2 \leqslant \min\{\delta,\nu_p(\Delta_{13})\}} \frac{(-1)^{\kappa'_2}}{2^{\kappa_2}}\\
& \times\sum_{\kappa_3,\kappa'_3 \geqslant 0 \atop \kappa_4+\kappa_3+\kappa'_3 \leqslant \min\{\delta,\nu_p(\Delta_{12})\}} \frac{(-1)^{\kappa'_3}}{2^{\kappa_3}}\sum_{\kappa_5,\kappa'_5 \geqslant 0 \atop \kappa_5+\kappa'_5 \leqslant \min\{\delta'_1+\delta'_2,\nu_p(\Delta_{12})\}} \frac{(-1)^{\kappa'_5}}{2^{\kappa_5}}\\
&\times \sum_{\boldsymbol{\nu} \in \mathbb{N}^3} \frac{\tilde{\rho}\left(p^{N_1},p^{N_2},p^{N_3};F_1,F_2,F_3 \right)}{p^{2(N_1+N_2+N_3+1)}}.
\end{aligned}
$$
On raisonne comme dans la section 9.5 de \cite{39} et on utilise, pour $a$ et $b$ deux entiers naturels, la formule
$$
\sum_{\substack{\kappa,\kappa' \geqslant 0 \\ \kappa+\kappa' \leqslant \min\{a,b\}\\ 0 \leqslant \kappa' \leqslant 1}} \frac{(-1)^{\kappa'}}{z^{\kappa}}=\frac{1}{z^{\min\{a,b\}}}.
$$
On commence par se placer dans le cas d'un nombre premier $p \equiv 1 \Mod{4}$ qui divise tous les résultants $\Delta_{ij}$.
On obtient alors l'égalité
$$
f'_p(\boldsymbol{\delta},\boldsymbol{\delta}')=\sum_{\boldsymbol{\nu} \in \mathbb{N}^3} \frac{(\nu_1+1)(\nu_2+1)(\nu_3+1) \rho^{\dagger}(p^{N'_1},p^{N'_2},p^{N'_3};F_1,F_2,F_3)}{C(\delta,N'_1,N'_2,N'_3)2^{\min\{\delta'_1+\delta'_2,N'_1,N'_2,N'_3\}}p^{2(N'_1+N'_2+N'_3+1)}},
$$
où
$$
C(\delta,N'_1,N'_2,N'_3)=(3/8)^{\min\{\delta,N'_1,N'_2,N'_3\}}2^{\min\{\delta,N'_2,N'_3\}+\min\{\delta,N'_1,N'_3\}+\min\{\delta,N'_2,N'_3\}}
$$
et avec, lorsque $(a_1,a_2,a_3) \in \mathbb{N}$,
$$
\rho^{\dagger}(p^{a_1},p^{a_2},p^{a_3};F_1,F_2,F_3)=\#\left\{ (u,v) \in \left(\mathbb{Z}/p^{a_1+a_2+a_3+1}\mathbb{Z}\right)^2 \hspace{2mm} \big| \hspace{2mm} p^{a_i} \parallel  F_i(u,v), \hspace{2mm} p \nmid (u,v) \right\}.
$$
Il reste à voir sur quels $\boldsymbol{\delta}'=(\delta'_{1},\delta'_{2},\delta'_{3})$ on somme. D'après (\ref{ccc}), on somme sur les triplets suivants:
$$
(0,0,0), \quad (0,0,1), \quad (1,0,1), \quad (0,1,1), \quad (1,0,0), \quad (0,1,0),
$$
si bien que
$$
f_p(\boldsymbol{\delta})
\begin{aligned}[t]
& =\sum_{\boldsymbol{\nu} \in \mathbb{N}^3}\frac{\rho^{\dagger}(p^{N''_1},p^{N''_2},p^{N''_3};F_1,F_2,F_3)}{C(\delta,N''_1,N''_2,N''_3)p^{2(N''_1+N''_2+N''_3+1)}}\bigg((\nu_1+1)(\nu_2+1)(\nu_3+1)-\nu_1 \nu_2(\nu_3+1) \\
&\hspace{3.7cm}\left. +\frac{1}{2}\nu_1\nu_3(\nu_2-1)+\frac{1}{2}\nu_2\nu_3(\nu_1-1)-\frac{1}{2}\nu_2\nu_3(\nu_1+1)-\frac{1}{2}\nu_1\nu_3(\nu_2+1)       \right)\\
&=\sum_{\boldsymbol{\nu} \in \mathbb{N}^3}(\nu_1+\nu_2+\nu_3+1)\frac{\rho^{\dagger}(p^{N''_1},p^{N''_2},p^{N''_3};F_1,F_2,F_3)}{C(\delta,N''_1,N''_2,N''_3)p^{2(N''_1+N''_2+N''_3+1)}},
\end{aligned}
$$
avec $N''_i=\nu_i+\delta_i$. En remarquant que pour $N''_i$ fixés, il y a $C(\delta,N''_1,N''_2,N''_3)$ triplets $\boldsymbol{\delta}$ tels que $N''_i=\nu_i+\delta_i$, on peut conclure exactement comme dans \cite[section 9.5]{39} que l'on obtient la quantité adéquate.\par
Sans être complètement exhaustif, on traite ensuite un cas éloquent dont l'adaptation aux cas restants ne pose aucune difficulté. Tout d'abord, on se place dans un cas où $\nu_p(\Delta_{12})=0$. Supposons alors par exemple que $\nu_p(\Delta_{13})=0$ et $\nu_p(\Delta_{23})\geqslant 1$. Dans ce cas de figure, on a clairement
$$
\min\{\delta'_1+\delta'_2,\nu_p(\Delta_{12}),N'_1,N'_2,N'_3\}=\min\{\delta,\nu_p(\Delta_{12}),\nu_p(\Delta_{13}),\nu_p(\Delta_{23}),N'_1,N'_2,N'_3)\}=0
$$
et un examen des conditions sur les $\boldsymbol{\delta}'$ dans la formule (\ref{ccc}) montre que la somme sur les~$\boldsymbol{\delta}'$ porte sur les triplets $(0,0,0)$ et $(1,0,0)$, si bien qu'on obtient 
$$
f_p(\boldsymbol{\delta})=
\sum_{\boldsymbol{\nu} \in \mathbb{N}^3}\frac{\rho^{\dagger}(p^{N'_1},p^{N'_2},p^{N'_3};F_1,F_2,F_3)}{C'(\delta,N'_1,N'_2,N'_3)p^{2(N'_1+N'_2+N'_3+1)}}\left((\nu_1+1)(\nu_2+1)(\nu_3+1)-\nu_2 \nu_3(\nu_1+1)\right),
$$
où
$$
C'(\delta,N'_1,N'_2,N'_3)=2^{\min\{\delta,N'_2,N'_3\}+\min\{\delta,N'_1,N'_3\}+\min\{\delta,N'_2,N'_3\}}.
$$
Mais, le fait que $p$ ne divise ni $\Delta_{12}$ ni $\Delta_{13}$ implique que nécessairement $\nu_1\nu_2=\nu_1\nu_3=0$ de sorte que
$$
(\nu_1+1)(\nu_2+1)(\nu_3+1)-\nu_2 \nu_3(\nu_1+1)=\nu_1+\nu_2+\nu_3+1+\nu_1\nu_2+\nu_1\nu_3=\nu_1+\nu_2+\nu_3+1
$$
et on peut à nouveau conclure comme dans la section 9.5 de \cite{39}. 
\par
Pour finir, dans le cas $p \equiv 3\Mod{4}$, on a
$$
g(\boldsymbol{\mu})=\sum_{\boldsymbol{\nu} \in \mathbb{N}^3} \frac{\rho^{\dagger}\left( p^{2\nu_1+\mu_{1}},p^{2\nu_2+\mu_{2}},p^{2\nu_3+\mu_{3}}; F_1,F_2,F_3 \right)}{p^{2(2\nu_1+2\nu_2+2\nu_3+1)}}
$$
que l'on traite comme dans la section 9.5 de \cite{39}. Cela permet finalement de conclure à l'égalité
$$
c_0=\sum_{\boldsymbol{\varepsilon} \in \Sigma \atop \mathbf{m} \in M}\omega_{\infty}(\boldsymbol{\varepsilon},\mathbf{m}) \prod_p \omega_p(\boldsymbol{\varepsilon},\mathbf{m})
$$
et achève la preuve du fait que $c_0=c_S$.

\bibliographystyle{plain-fr}
\nocite{*}
\bibliography{bibliogr10}
\end{document}